# ASYMPTOTIC STABILITY OF STEADY STATES FOR THE COMPRESSIBLE NAVIER-STOKES-RIESZ SYSTEM IN THE PRESENCE OF VACUUM

JOSÉ A. CARRILLO, RENJUN DUAN, ANETA WRÓBLEWSKA-KAMIŃSKA, AND JUNHAO ZHANG

ABSTRACT. We consider a one-dimensional physical vacuum free boundary problem on the compressible Navier-Stokes-Riesz system for an attractive Riesz potential $|x|^{2s-1}/(2s-1)$ with $0 < s < 1/2$. It is proved that for the adiabatic constant $\gamma$ satisfying $2(1-s) < \gamma < 1 + 2s/3$ under the additional condition that $3/8 < s < 1/2$, there exists a unique global-in-time strong solution. Specifically, we establish the Lyapunov-type stability of the compactly supported steady states in the Lagrangian coordinates and we also obtain the time rate of convergence for the strong solution to steady states with the same mass in weighted Sobolev spaces where the weights indicate the behavior of solutions near the vacuum free boundary. The difficulties and challenges in the proof are caused not only by the degeneracy due to the vacuum free boundary but also by the non-local feature of the Riesz potential.

## 1. INTRODUCTION

The evolution of compressible viscous flow with a self-consistent force in $\mathbb{R}^n$ can be modeled by the following system of equations,

$$\partial_t \rho + \text{div}(\rho \mathbf{u}) = 0, \tag{1.1}$$

$$\partial_t(\rho \mathbf{u}) + \text{div}(\rho \mathbf{u} \otimes \mathbf{u}) + \text{div}\mathfrak{S} + \chi \rho \nabla \Psi_k = 0, \tag{1.2}$$

for $t > 0$, $\mathbf{x} \in \mathbb{R}^n$. Here, $\rho > 0$ and $\mathbf{u} \in \mathbb{R}^n$ are the density and velocity of the fluid, respectively. The stress tensor $\mathfrak{S}$ is given by

$$\mathfrak{S} := pI_n - \mu \left( \nabla \mathbf{u} + \nabla \mathbf{u}^{\mathsf{T}} - \frac{2}{3}(\text{div}\mathbf{u})I_n \right) - \lambda(\text{div}\mathbf{u})I_n,$$

where $I_n$ is the $n \times n$ identity matrix, $\mu > 0$ is the shear viscosity, $\lambda > 0$ is the bulk viscosity, and $p$ is the pressure that satisfies the equation of state,

$$p := K\rho^{\gamma},$$

with constant $K > 0$ and adiabatic exponent $\gamma > 1$. The interaction potential $\Psi_k$ is given by $\Psi_k := W_k * \rho$ where $k = 2s - n$ with $s \in (0, n/2)$ and

$$\nabla \Psi_k(\mathbf{x}) := \begin{cases} \nabla(W_k * \rho)(\mathbf{x}), & \text{for } 1 - n < k < 0, \\ \int_{\mathbb{R}^n} \nabla W_k(\mathbf{x} - \mathbf{y})(\rho(\mathbf{y}) - \rho(\mathbf{x})) \, \text{d}\mathbf{y}, & \text{for } -n < k \le 1 - n. \end{cases} \tag{1.3}$$

Date: January 2, 2026.

2020 Mathematics Subject Classification. 35Q30, 35R35, 35B40.

Key words and phrases. Navier-Stokes-Riesz system, free boundary problem, compactly supported stationary profile, stability, a priori estimates.

Notice that in the more singular range $-n < k \leq 1 - n$, we need to substract a term to make the integral convergent assuming suitable regularity of the density $\rho$. Alternatively, one could defined the integral using Cauchy principal value integrals. Here, $W_k$ is the Riesz kernel given by

$$W_k(\mathbf{x}) := \begin{cases} \dfrac{|\mathbf{x}|^k}{k}, & \text{for } -n < k < 0, \\ \log|\mathbf{x}|, & \text{for } k = 0. \end{cases}$$

It is also known that $\Psi_k$ is governed by a fractional diffusion process,

$$c_{n,s}(-\Delta)^s \Psi_k = \rho,$$

where

$$c_{n,s} = \begin{cases} (2s-n)\dfrac{\Gamma(n/2-s)}{\pi^{n/2}4^s\Gamma(s)}, & \text{for } k < \min\{0, 2-n\}, \\ 2\pi, & \text{for } k = 0 \text{ and } n = 2, \end{cases}$$

with $\Gamma$ being the gamma function. In (1.2), $\chi \in \mathbb{R}$ represents the strength of repulsion if $\chi < 0$ or attraction if $\chi > 0$. The system (1.1)-(1.2) is often called the compressible Navier-Stokes-Riesz system, which can be applied to describe Riesz gases with applications in plasma physics, star or galaxy dynamics, and swarming models in mathematical biology, see for instance [1, 6, 5, 27]. From now on, we will fix $K = \chi = 1$, since we are dealing only for the attractive case.

It is well-known that the system admits non-trivial steady solutions caused by two competing behaviors: the repulsion (pressure) between fluid particles and the attraction via a non-local self-consistent force. Indeed, by setting $\mathbf{u} = 0$, the system (1.1)-(1.2) is reduced to

$$\nabla(\bar{\rho}^\gamma) + \bar{\rho}\nabla\bar{\Psi}_k = 0, \tag{1.4}$$

where $\bar{\Psi}_k$ satisfies (1.3) with $\rho$ substituted by $\bar{\rho}$. We call (1.4) by the steady diffusion-aggregation equation. According to results in [4, 8, 9, 2] (see a summary in Appendix A), for $\gamma_c < \gamma < \gamma^*$, where $\gamma_c := 2(1 - s/n)$ and

$$\gamma^* := \begin{cases} \dfrac{k-1}{k}, & \text{for } n \geq 1 \text{ and } 0 < s < \tfrac{1}{2}, \\ +\infty, & \text{for } n \geq 2 \text{ and } \tfrac{1}{2} \leq s < \tfrac{n}{2}, \end{cases}$$

all stationary states $\bar{\rho} \geq 0$, see precise definition in Appendix A, satisfying (1.4), are radially symmetric decreasing about their center of mass. Moreover, they are unique among the class of radially decreasing stationary states, and they can be characterized as the unique global minimizer of the nonlinear functional $\mathcal{F}$ associated with (1.4),

$$\mathcal{F}[\rho] := \int_{\mathbb{R}^n} \frac{1}{\gamma - 1}\rho^\gamma + \frac{1}{2}\rho(W_k * \rho)\,\mathrm{d}\mathbf{x}, \tag{1.5}$$

in the class

$$\mathcal{Y} := \left\{ \rho \in L^1_+(\mathbb{R}^n) \cap L^\gamma(\mathbb{R}^n), \|\rho\|_{L^1} = 1, \int_{\mathbb{R}^n} \mathbf{x}\rho(\mathbf{x})\,\mathrm{d}\mathbf{x} = 0 \right\}.$$

Moreover, they are compactly supported $\mathrm{supp}(\bar{\rho}) = \bar{B}_{\bar{R}}$, $\bar{R} > 0$, and $\bar{\rho} \in W^{1,\infty}(\mathbb{R}^n)$ for $\gamma_c < \gamma \leq 2$. Therefore, $(\bar{\rho}, 0)$ is a steady solution to the system (1.1)-(1.2) and $(\bar{\rho}, 0)$ is also a global minimizer of the following energy functional $\mathcal{E}$ associated with (1.6)-(1.7),

$$\mathcal{E}[\rho, \mathbf{u}] := \int_{\mathbb{R}^n} \frac{1}{2}\rho|\mathbf{u}|^2 + \frac{1}{\gamma - 1}\rho^\gamma + \frac{1}{2}\rho W_k * \rho\,\mathrm{d}\mathbf{x}.$$

Indeed, for an arbitrary solution $(\rho, \mathbf{u})$ to (1.1)-(1.2), we have

$$\mathcal{E}[\rho, \mathbf{u}] \geq \mathcal{E}[\rho, 0] = \mathcal{F}[\rho] \geq \mathcal{F}[\bar{\rho}].$$

In this paper, our aim is to prove the global-in-time stability of the steady solution $(\bar{\rho}, 0)$ by considering a one-dimensional $(n = 1)$ vacuum free boundary problem for the compressible Navier-Stokes-Riesz system (1.1)-(1.2) with $2(1-s) < \gamma \leq 2$. To be more specific, we consider

$$\partial_\tau \rho + \partial_\xi (\rho u) = 0, \tag{1.6}$$

$$\rho(\partial_\tau u + u \partial_\xi u) + \partial_\xi p = \partial_{\xi\xi} u - \rho \partial_\xi \Psi_k, \tag{1.7}$$

for $\tau > 0$ and $0 < \xi < a(\tau)$ where $k = 2s-1$, and $a(\tau)$ is a free endpoint determined by

$$\begin{cases} \dfrac{\mathrm{d}a(\tau)}{\mathrm{d}\tau} = u(\tau, a(\tau)), \\ a(0) = a_0, \end{cases} \tag{1.8}$$

for some fixed $a_0 > 0$. Here, we set $4\mu/3 + \lambda$ to be unity for convenience. At the origin, we impose a no-flow boundary condition,

$$u(\tau, 0) = 0, \tag{1.9}$$

while at the free endpoint, we impose the Neumann condition,

$$\partial_\xi u(\tau, a(\tau)) = 0. \tag{1.10}$$

From now on, we will omit the dependence for $\Psi$ and $W$ on $k$ for the sake of simplicity.

The initial condition is given by

$$(\rho, u)(0, \xi) = (\rho_0, u_0)(\xi), \tag{1.11}$$

with $\rho_0$ satisfying

$$\int_0^{a_0} \rho_0(\xi)\,\mathrm{d}\xi = \int_0^{\bar{R}} \bar{\rho}(x)\,\mathrm{d}x, \quad \rho_0(a_0) = 0, \quad (\rho_0^{\gamma-1}(a_0))' = (\bar{\rho}^{\gamma-1}(\bar{R}))', \tag{1.12}$$

which indicates that, initially, the mass of the fluid equals the mass of the steady state, the fluid continuously connects to the vacuum across the free boundary, and the sound speed $c = \sqrt{\rho^{\gamma-1}}$ coincides with that of the steady state to the first-order derivative at the vacuum boundary. Moreover, by the conservation of mass (1.6) and (1.10), we have further $\rho(\tau, a(\tau)) = 0$ for all $\tau > 0$. It also should be noticed that, since $\bar{\rho} \in W^{1,\infty}$ when $2(1-s) < \gamma \leq 2$,

$$-\infty < (\rho_0^{\gamma-1}(a_0))' = (\bar{\rho}^{\gamma-1}(\bar{R}))' < 0,$$

which implies that the vacuum boundary is physical for the compressible Euler equations according to the definition in Jang-Masmoudi's work [16, 20].

Moreover, the free boundary problem (1.6)-(1.7) with boundary conditions (1.9), (1.10) and initial conditions (1.11) admits a symmetric extension, that is, if we define $(\hat{\rho}, \hat{u})$ by

$$(\hat{\rho}, \hat{u})(\tau, \xi) = \begin{cases} (\rho, u)(\tau, \xi), & \xi \geq 0, \\ (\rho, -u)(\tau, -\xi), & \xi < 0, \end{cases}$$

then $(\hat{\rho}, \hat{u})$ is a solution to (1.6)-(1.7) over the free boundary $(-a(\tau), a(\tau))$ with boundary conditions and initial conditions extended correspondingly from (1.10) and (1.11). Since $\hat{u}$ is an odd extension of $u$, we have $\hat{u}(t, 0) = 0$. Without loss of generality, we may focus our attention on $(\hat{\rho}, \hat{u})$. For the sake of simplicity, we still denote the extended solution by $(\rho, u)$.

Since $\Psi$ is defined over $\mathbb{R}$, we utilize a natural extension of $\rho$ denoted by $\tilde{\rho}$,

$$\tilde{\rho}(\xi) := \begin{cases} \rho(\xi), & -a(\tau) \leq \xi \leq a(\tau), \\ 0, & |\xi| > a(\tau). \end{cases}$$

Therefore by the equation (1.3), $\Psi$ has the following form,

$$\partial_\xi \Psi(\tau, \xi) = \int_{-\infty}^{\infty} W'(\xi - \eta)(\tilde{\rho}(\tau, \eta) - \tilde{\rho}(\tau, \xi)) \, \mathrm{d}\eta.$$

Our results can be stated as follows: assume that the initial data $(\rho_0, u_0, a_0)$ is a small perturbation of the steady solutions $(\bar{\rho}, 0, \bar{R})$ in a suitable sense (see Theorem 2.1 and Theorem 2.2) with the same total mass, then for the adiabatic exponent $2(1-s) < \gamma < 1 + 2s/3$ under the additional condition that $3/8 < s < 1/2$, there is a unique global-in-time strong solution $(\rho, u, a)$ to the free boundary problem of (1.6)-(1.7) with boundary conditions (1.9), (1.10) and initial conditions (1.11). Moreover, let $\eta(t, x)$ be defined by (2.1), (2.2), and (2.3). Then the steady solution $(\bar{\rho}, 0)$ is stable in the following sense

$$\|(\eta(t, x) - x, \rho(t, \eta(t, x)) - \bar{\rho}(x), u(t, \eta(t, x)))\|_{L_x^\infty}$$
$$\leq C \|(\eta_0(x) - x, \rho_0(\eta_0(x)) - \bar{\rho}(x), u_0(\eta_0(x)))\|_{L_x^\infty}, \quad t \geq 0,$$

for some constant $C > 0$ uniformly in time $t$. Decay rates can also be obtained under some suitably weighted $L^2$ norms. Here, the restriction $\gamma < 1 + 2s/3$ is a technical requirement for studying the problem in an appropriate functional space. For details, see Section 2.2.

The difficulties and challenges of the vacuum free boundary value problem for the compressible Navier-Stokes-Riesz system are not only caused by the degeneracy due to the vacuum free boundary, but also caused by the non-local nature of the Riesz potential, which forces us to take into account of the effect of the solution over the whole domain even if we focus on the local behavior of the solution. This makes a huge difference compared to cases involving the Newtonian potential.

Let us conclude this section by reviewing some related works. The physical vacuum free boundary problem is a long-standing research topic in the compressible hydrodynamics due to its physical importance and mathematical difficulties. In the theory of inviscid flows, significant results on the local-in-time existence have been investigated. For the compressible Euler system, Jang-Masmoudi [16] made the breakthrough in establishing the local-in-time well-posedness of vacuum free boundary problem in an one-dimensional case. Then Coustand-Shkoller [10] obtained the smoother local solutions all the way to the moving boundary in higher-order Sobolev spaces. Later in [11] they extended their method to three dimensional cases. A similar result was also obtained by Jang-Masmoudi [20] using a different approach. If the self-gravitation is taking into consideration, Luo-Xin-Zeng [22] proved a general uniqueness theorem and established a new local-in-time existence theory for spherically symmetric motions in a functional space involving less derivatives than those constructed before. Hadžić-Jang [13] even obtained nonlinear stability of some compactly supported expanding star-solutions of the mass-critical gravitational Euler-Poisson system, then they further constructed the global-in-time solutions without any symmetry assumptions in [14]. In the theory of viscous flows, not only the local-in-time results but also the global-in-time results have been found. Jang [17] established the local-in-times well-posedness of strong solutions. Later, Zeng [29] studied the vacuum free boundary problem of the compressible Navier-Stokes system and constructed the global-in-time strong solutions. Ou-Zeng [28] extended the study to more physical situation where the viscosity is assumed to be density-dependent and an externally gravity force exists. If the self-gravitation is considered, there exists an important class of solutions called Lane-Emden solutions which are spherically symmetric and stable, and minimize the energy among all possible solutions (cf. [12]). An important progress concerning the Lane-Emden solutions has been made by Luo-Xin-Zeng [24, 25]. They consider the free boundary problem of the compressible Navier-Stokes-Poisson system and proved that for the viscous gaseous star with $4/3 < \gamma < 2$, the Lane-Emden solution is asymptotically nonlinear stable, moreover, the global-in-time regularity uniformly up to the vacuum boundary of solutions is proved. Later, Luo-Wang-Zeng [26] generalized the previous results by considering white dwarfs and polytropes for all $\gamma > 4/3$. The extra dissipation induced by viscosity plays an important role in obtaining these global-in-time results. It has been noted that the dissipation induced by damping also guarantee the global-in-time results. For results in this direction, we refer the readers to the results [23, 30, 31, 32]. Despite the stability results, we refer the readers to [15, 18, 19] for results concerning the instability.

Very recently, Carrillo-Charles-Chen-Yuan [7] obtained a global weak solution with spherical symmetry and finite energy for the Cauchy problem of the compressible Euler-Riesz system, the nonlinear stability of steady states is also verified in a certain sense. The proof of their result relies on the vanishing viscosity method. As an intermediate step, they also studied the free boundary problem of the compressible Navier-Stokes-Riesz system. To solve the free boundary problem, the Lagrangian transformation is also incredibly useful in their analysis. Since the solution to this free boundary problem is an approximation to the solution (containing vacuum) to the Euler-Riesz system, they can lifted the free boundary problem to consider a non-vacuum case, which directly ensures that the Lagrangian transformation

$$(t, r) \mapsto (\tau := t, x := \int_a^r \rho(t, y) y^{n-1} \, \mathrm{d}y),$$

is a smooth bijective mapping. Indeed, we have

$$\nabla_{(t,r)} x = (-\rho u r^{n-1}, \rho r^{n-1}), \quad \nabla_{t,r} \tau = (1, 0), \tag{1.13}$$

which implies that the Lagrangian transformation is bijective by inverse function theorem. If vacuum appears, then (1.13) is not applicable to the inverse function theorem. Our result is important in the following two aspects: First, we succeed in justify that the Lagrangian transformation is a bijective mapping by considering the strong solution to the compressible Navier-Stokes-Riesz system; Second, we obtain a time-asymptotic behavior of the solution and a time decay estimate in some particular weighted $L^2$ norms where the weights also indicate the behavior of the solution close to the vacuum free boundary.

The paper's structure is organized as follows: In Section 2, we describe the problem under the Lagrangian transformation, then the main results are presented. In Section 3, we establish three fundamental estimates for the strong solution framework: (1) Lower-order estimates: Initial $L^2$ energy estimates derived via wave equation techniques, enhanced through viscous evolution analysis of spatial derivatives near vacuum. (2) Pointwise characterization: Bootstrap analysis of perturbation amplitudes using Taylor-expanded representation formula of $\partial_x \eta$. (3) Higher-order regularity: Extension of lower-order estimates with higher-order estimates in temporal regularity, as well as the higher-order estimates in second-order spatial derivatives. The principal technical challenge involves managing the Riesz potential's singular non-local structure, which is addressed through decompositions inspired by the results of lower-order energy and higher-order energy estimates. In addition, some time-decay estimates are also established. In Section 4, we provide a concise difference scheme formulation for demonstrating local existence and uniqueness of strong solutions, subsequently employing Section 3's results to prove **Theorem 2.1** and **Theorem 2.2**. Appendices include essential background material: existence theory for steady states of compressible Navier-Stokes-Riesz systems, and Hardy-type inequalities extensively utilized in our analysis.

## 2. Lagrangian formulation and main results

2.1. **Lagrangian formulation.** In order to fix the free boundary, we formulate the problem in Lagrangian variable. Let $x$ be the reference variable and define the Lagrangian variable $\eta(t, x)$ by

$$\frac{\mathrm{d}\eta(t, x)}{\mathrm{d}t} = u(t, \eta(t, x)), \tag{2.1}$$

for $t > 0$ and $-\bar{R} \leq x \leq \bar{R}$, with the initial condition

$$\eta(0, x) = \eta_0(x), \tag{2.2}$$

where $\eta_0(x)$ is the initial position which maps $[-\bar{R}, \bar{R}]$ to $[-a_0, a_0]$ satisfying

$$\int_{\eta_0(-x)}^{\eta_0(x)} \rho_0(y)\,\mathrm{d}y = \int_{-x}^{x} \bar{\rho}(y)\,\mathrm{d}y, \quad \eta_0(-x) = -\eta_0(x), \quad \eta_0'(\pm\bar{R}) = 1, \tag{2.3}$$

so that

$$\rho_0(\eta_0(x))\eta_0'(x) = \bar{\rho}(x). \tag{2.4}$$

In fact, $\eta_0$ can be described by

$$\eta_0(x) = \varphi^{-1}(\phi(x)),$$

where $\phi$ and $\varphi$ are defined by

$$\phi(x) := \int_{-x}^{x} \bar{\rho}(y)\,\mathrm{d}y, \quad -\bar{R} \leq x \leq \bar{R},$$

$$\varphi(\xi) := \int_{-\xi}^{\xi} \rho_0(y)\,\mathrm{d}y, \quad -a_0 \leq \xi \leq a_0.$$

In addition, the third condition in (2.3) also implies the sound speeds of the initial state and the steady state have the same behavior at the vacuum boundary up to their first-order derivatives. Indeed, let us consider the case at $x = \bar{R}$, as the other case at $x = -\bar{R}$ can be discussed in a similar way. By using (2.4), we have

$$\lim_{x \to \bar{R}} \frac{(\rho_0^{\gamma-1}(\eta_0(x)))'}{(\bar{\rho}^{\gamma-1}(x))'} = \lim_{x \to \bar{R}} (\eta_0'(x))^{-\gamma+3} = 1,$$

which indicates the third condition in (1.12).

We define the Lagrangian density and velocity respectively by

$$d(t, x) := \rho(t, \eta(t, x)), \quad v(t, x) := u(t, \eta(t, x)),$$

for $t \geq 0$ and $-\bar{R} \leq x \leq \bar{R}$. Then $(d, v)$ satisfies

$$\partial_t(d\partial_x \eta) = 0, \tag{2.5}$$

$$d\partial_t v + \frac{\partial_x(d^\gamma)}{\partial_x \eta} = \frac{1}{\partial_x \eta}\partial_x\left(\frac{\partial_x v}{\partial_x \eta}\right) - d\partial_x\psi, \tag{2.6}$$

for $t > 0$ and $-\bar{R} < x < \bar{R}$, where $\partial_x\psi$ is defined by

$$\partial_x\psi(t, x) := \int_{-\infty}^{\infty} W'(\eta(t, x) - \eta(t, y))(d(t, y) - d(t, x))\partial_y\tilde{\eta}(t, y)\,\mathrm{d}y,$$

$$\tilde{\eta}(t, x) := \begin{cases} -\eta(t, \bar{R}) + \bar{R} + x, & x < -\bar{R}, \\ -\eta(t, -x), & -\bar{R} \leq x < 0, \\ 0, & x = 0, \\ \eta(t, x), & 0 < x \leq \bar{R}, \\ \eta(t, \bar{R}) - \bar{R} + x, & x > \bar{R}. \end{cases}$$

We note that the Lagrangian density can be solved from (2.5) and (2.4),

$$d(t, x) = \frac{\bar{\rho}(x)}{\partial_x\eta(t, x)}, \tag{2.7}$$

and $v = \partial_t\eta$ due to (2.1), then (2.6) be rewritten as,

$$\bar{\rho}\partial_{tt}\eta + \partial_x\left[\left(\frac{\bar{\rho}}{\partial_x\eta}\right)^\gamma - \bar{\rho}^\gamma\right] + \Phi = \partial_x\left(\frac{\partial_{tx}\eta}{\partial_x\eta}\right), \tag{2.8}$$

where $\Phi$ is given by

$$\Phi(t, x) := \partial_x(\bar{\rho}^\gamma(x))$$
$$+ \frac{\bar{\rho}(x)}{\partial_x\eta(t, x)}\int_{-\infty}^{\infty} W'(\eta(t, x) - \tilde{\eta}(t, y))(\bar{\rho}(y)\partial_x\eta(t, x) - \bar{\rho}(x)\partial_y\tilde{\eta}(t, y))\,\mathrm{d}y. \tag{2.9}$$

Notice that $\Phi(t, x) = 0$ vanishes for $\eta(t, x) = x$, which is the expected long time asymptotics due to the stationary equation (1.4).

Let us derive the corresponding boundary conditions for (2.8). First, by (2.10) and (2.3),

$$\partial_x v(t, \pm\bar{R}) = 0. \tag{2.10}$$

Then since $\eta$ satisfies the differential equation (2.1), we have

$$\eta(t, x) = \eta_0(x) + \int_0^t v(s, x)\,\mathrm{d}s. \tag{2.11}$$

Taking $\partial_x(2.11)$ and then using the boundary conditions (2.10) and (2.3), we find

$$\partial_x\eta(t, \pm\bar{R}) = 1. \tag{2.12}$$

In addition, besides the initial condition (2.2), (1.11) also implies

$$\partial_t\eta(0, x) = u(0, \eta_0(x)) := \eta_1(x). \tag{2.13}$$

2.2. **Main results.** In what follows, we give the definition of strong solutions to (2.8) with (2.2), (2.13) and (2.12).

**Definition 2.1.** $\eta$ is called a strong solution to (2.8) with initial conditions (2.2), (2.13) and boundary condition (2.12) in $[0, T] \times [-\bar{R}, \bar{R}]$ for some constant $T > 0$, if

(1) $c_1 \leq \partial_x \eta(t, x) \leq c_2$ for $0 \leq t \leq T$ and $-\bar{R} \leq x \leq \bar{R}$ with some constants $c_1, c_2 > 0$;

(2) $\partial_t \eta \in C([0, T]; W^{1,\infty}(-\bar{R}, \bar{R}))$;

(3) $\bar{\rho}^{\frac{1}{2}} \partial_t \eta \in C^1([0, T]; L^2(-\bar{R}, \bar{R}))$, $\bar{\rho}^{-\frac{1}{2}} \partial_{xx} \eta \in C([0, T]; L^2(-\bar{R}, \bar{R}))$, and $\bar{\rho}^{-\frac{1}{2}} \partial_x \left( \frac{\partial_x \eta}{\partial_x \eta} \right) \in C([0, T]; L^2(-\bar{R}, \bar{R}))$;

(4) $\partial_x \eta(t, \pm \bar{R}) = 1$ hold in the sense of $H^1$-trace for $0 \leq t \leq T$;

(5) (2.8) holds over $[0, T] \times [-\bar{R}, \bar{R}]$ a.e.

**Remark 2.1.** *The condition (1) states a positive lower bound for $\partial_x \eta$. Then, by the inverse function theorem, the transformation from $x$ to $\xi$ is invertible, which means that the system (1.6)-(1.7) in Eulerian coordinates and the system (2.5)-(2.6) in Lagrangian coordinates are equivalent.*

**Remark 2.2.** *The condition (2) indicates that $u$ is Lipschitz continuous in $\eta$ as $\partial_\eta u(t, \eta(t, x)) = \frac{\partial_x v(t,x)}{\partial_x \eta(t,x)}$. Therefore, by Picard-Lindelöf's theorem, the ordinary differential equation (2.1) admits a unique solution $\eta$ which ensures the well-definiteness of Lagrangian variable.*

**Remark 2.3.** *The condition (3) gives the function spaces such that each term in (2.8) multiplied by $\bar{\rho}^{-1/2}$ is in $C([0, T]; L^2(-\bar{R}, \bar{R}))$. Indeed,*

$$\left\| \bar{\rho}^{-\frac{1}{2}} \partial_x \left[ \left( \frac{\bar{\rho}}{\partial_x \eta} \right)^\gamma \right] \right\|_{L^2(-\bar{R}, \bar{R})}$$
$$\leq C \left( \left\| \frac{\bar{\rho}^{\gamma - \frac{1}{2}} \partial_{xx} \eta}{|\partial_x \eta|^{\gamma+1}} \right\|_{L^2(-\bar{R}, \bar{R})} + \left\| \frac{\bar{\rho}^{-\frac{1}{2}} \partial_x (\bar{\rho}^\gamma)}{(\partial_x \eta)^\gamma} \right\|_{L^2(-\bar{R}, \bar{R})} \right)$$
$$\leq C(\|\bar{\rho}^{\gamma - \frac{1}{2}} \partial_{xx} \eta\|_{L^2(-\bar{R}, \bar{R})} + 1),$$

*which together with $\bar{\rho} \in W^{1,\infty}$ implies that $\bar{\rho}^{-\frac{1}{2}} \partial_x \left[ \left( \frac{\bar{\rho}}{\partial_x \eta} \right)^\gamma - \bar{\rho}^\gamma \right] \in C([0, T]; L^2(-\bar{R}, \bar{R}))$. Moreover, by **Lemma B.4**,*

$$\|\bar{\rho}^{-\frac{1}{2}} \Phi(t)\|_{L^2} \leq C(\|\partial_x \eta(t)\|_{L^\infty}^2 + \|\bar{\rho}^{\frac{2\gamma-1}{2}} \partial_{xx} \eta(t)\|_{L^2}^2),$$

*which implies that $\bar{\rho}^{-\frac{1}{2}} \Phi \in C([0, T]; L^2(-\bar{R}, \bar{R}))$. Therefore, each term in (2.8) belongs to $C([0, T]; L^2(-\bar{R}, \bar{R}))$. As a conclusion, the condition (5) follows.*

**Remark 2.4.** *If conditions (1), (2), and (3) in **Definition 2.1** hold, then*

$$\|\partial_x v(t)\|_{H^1} \leq \|\partial_x \eta(t)\|_{L^\infty} \left\| \frac{\partial_{tx} \eta(t)}{\partial_x \eta(t)} \right\|_{H^1}$$
$$\leq \|\partial_x \eta(t)\| \left( \left\| \frac{1}{\partial_x \eta(t)} \right\|_{L^\infty} \|\partial_t \eta(t)\|_{W^{1,\infty}} + \left\| \partial_x \left( \frac{\partial_{tx} \eta(t)}{\partial_x \eta(t)} \right) \right\|_{L^2} \right),$$

which implies that $\partial_x v \in C([0,T]; H^1(-\bar{R}, \bar{R}))$. Therefore $\partial_x v(t, \pm \bar{R})$ exists in the sense of $H^1$-trace. Since

$$\partial_x \eta(t,x) - \partial_x \eta_0(x) = \int_0^t \partial_x v(s,x)\, \mathrm{d}s,$$

for arbitrary $t > 0$, by condition (4) and $\partial_x \eta_0(\pm \bar{R}) = 1$, we have $\partial_x v(t, \pm \bar{R}) = 0$ in the sense of $H^1$-trace.

**Remark 2.5.** *If $\eta$ exists in the sense of* **Definition 2.1**, *besides the boundary conditions* (2.12), *we have*

$$\eta(t,0) = 0. \tag{2.14}$$

Indeed, let $x = 0$, then (2.1) becomes an ordinary differential equation in $t$. Recall that, at the origin, $u(t,0) = 0$ for $t \geq 0$ and $\eta_0(0) = 0$ due to the oddness. It can be directly verified that $\eta(t,0) = 0$ is a solution to this ordinary differential equation. Since $\partial_t \eta \in W^{1,\infty}$, by Picard-Lindelöf's theorem, $\eta(t,0) = 0$ is also the unique solution.

Having defined the strong solution $\eta$, we can construct a unique global strong solution $(\rho, u, a)$ to the free boundary value problem (1.6)-(1.11). Indeed, by condition (1) in **Definition** 2.1 and the inverse function theorem, $\eta$ defines a diffeomorphism from the fixed reference domain $[-\bar{R}, \bar{R}]$ to the free boundary domain $[-a(t), a(t)]$ with

$$a(t) := \eta(t, \bar{R}), \tag{2.15}$$

for $t \geq 0$. Moreover, we denote the inverse of the map $\eta$ by $\eta^{-1}$. We defined $(\rho, u)$ as

$$\rho(\tau, \xi) := \frac{\bar{\rho}(x)}{\partial_x \eta(\tau, x)}, \quad u(\tau, \xi) := \partial_\tau \eta(\tau, x), \tag{2.16}$$

with $x = \eta^{-1}(\tau, \xi)$ for $\tau \geq 0$ and $-a(\tau) \leq \xi \leq a(\tau)$. Accordingly, we give the definition of strong solutions to (1.6) and (1.7) with (1.9), (1.10) and (1.11).

**Definition 2.2.** *$(\rho, u, a)$ is called a strong solution to the compressible Navier-Stokes-Riesz system* (1.6) *and* (1.7) *with initial condition* (1.11) *and boundary conditions* (1.9), (1.10) *in* $[0,T] \times [0, \bar{R}]$ *for some constant* $T > 0$, *if*

(1) *$\eta$ is a strong solution to* (2.8) *with* (2.2), (2.13) *and* (2.12);
(2) *$(\rho, u, a)$ satisfies* (2.15) *and* (2.16) *a.e.*

Now let us state the first main theorem.

**Theorem 2.1.** *Let $\frac{3}{8} < s < \frac{1}{2}$ and $2(1-s) < \gamma < 1 + 2s/3$. Assume that* (2.3) *is valid and the initial value $\eta_0$ is an odd function satisfying the compatibility condition $\partial_x \eta_0(\pm \bar{R}) = 1$. Moreover, let us assume that that for some $\varepsilon_0 > 0$ we have*

$$\|(\partial_x \eta_0 - 1, \partial_x \eta_1)\|_{L^\infty(-\bar{R}, \bar{R})}^2 + \|\bar{\rho}^{\gamma - \frac{1}{2}} \partial_{xx} \eta_0\|_{L^2(-\bar{R}, \bar{R})}^2 + \|\bar{\rho}^{\frac{1}{2}} \partial_t \eta_1\|_{L^2(-\bar{R}, \bar{R})}^2 \leq \varepsilon_0.$$

*Then the equation* (2.8) *with initial conditions* (2.2), (2.13) *and boundary conditions* (2.12) *admits a unique strong solution such that*

$$\begin{aligned}
\|(\partial_x \eta - 1, &\partial_{tx} \eta)(t)\|_{L^\infty(-\bar{R}, \bar{R})}^2 \\
&+ \|\bar{\rho}^{\gamma - \frac{1}{2}} \partial_{xx} \eta(t)\|_{L^2(-\bar{R}, \bar{R})}^2 + \|\bar{\rho}^{\frac{1}{2}} \partial_{tt} \eta(t)\|_{L^2(-\bar{R}, \bar{R})}^2 \leq C\varepsilon_0, \quad t \geq 0, \quad (2.17)
\end{aligned}$$

*for some constant $C > 0$ independent of time $t$. Furthermore, the solution satisfies the following time decay estimates,*

$$\|\bar{\rho}^{\frac{3}{2}}\partial_t^i(\partial_x\eta - 1)(t)\|_{L^2(-\bar{R},\bar{R})} \le C(1+t)^{-\frac{1}{2}}E^{\frac{1}{2}}(0), \quad t \ge 0, i = 0, 1, \qquad (2.18)$$

$$\|\bar{\rho}^{\frac{1}{2}}\partial_t^i\eta\|_{L^2(-\bar{R},\bar{R})} \le C(1+t)^{-\frac{1}{2}}E^{\frac{1}{2}}(0), \quad t \ge 0, i = 1, 2, \qquad (2.19)$$

*for some constant $C > 0$ independent of time $t$, where*

$$E(t) := \|(\partial_x\eta - 1, \partial_{tx}\eta)(t)\|_{L^\infty(-\bar{R},\bar{R})}^2$$
$$+ \|\bar{\rho}^{-\frac{1}{2}}\partial_{xx}\eta(t)\|_{L^2(-\bar{R},\bar{R})}^2 + \|\bar{\rho}^{\frac{1}{2}}\partial_{tt}\eta(t)\|_{L^2(-\bar{R},\bar{R})}^2. \qquad (2.20)$$

**Remark 2.6.** *Let us emphasize that $\gamma > 2(1-s)$ guarantees the existence of steady solutions and it is also essentially used by the energy estimates in Section 3. In addition, the condition $\gamma < 1+2s/3$ is needed for the definition of strong solution, and more precisely, it is required to have the regularity $\bar{\rho}^{-\frac{1}{2}}\Phi \in C([0,T]; L^2(-\bar{R},\bar{R}))$.*

The second main result concerns the original free boundary problem (1.6)-(1.11) in Eulerian coordinates.

**Theorem 2.2.** *Let $\frac{3}{8} < s < \frac{1}{2}$ and $2(1-s) < \gamma < 1+2s/3$. Assume that the initial condition (1.11) satisfies (1.12) and the compatibility condition $u_0(0) = \partial_\xi u_0(a_0) = 0$. Then $(\rho, u, a)$ defined by (2.15) and (2.16) is the unique global strong solution to the system (1.6)-(1.7) with initial condition (1.11) and boundary conditions (1.9), (1.10).*

**Remark 2.7.** *If the free boundary problem (1.6)-(1.7) with initial condition (1.11) and boundary conditions (1.9), (1.10) admits a classic solution, then by the discussion in Section 2.1, the classic solution is in the form of (2.16).*

**Remark 2.8.** *From the construction of (2.16), the strong solution $\rho$ is non-negative over the free boundary domain and continuously connected to the vacuum near the free boundary.*

**Remark 2.9.** *The sound speed of the flow in our construction is $C^{1/2}$-Hölder continuous near the vacuum. Indeed, note that*

$$|\rho^{\gamma-1}(\tau,\xi)| \le \frac{1}{c_1^{\gamma-1}}\left|\bar{\rho}^{\gamma-1}(x(\tau,\xi))\right|.$$

*This together with the fact that $\bar{\rho}^{\gamma-1} \in W^{1,\infty}$ implies that the sound speed $c = \sqrt{\rho^{\gamma-1}}$ is of $C^{1/2}$-Hölder continuity near the vacuum.*

Let us illustrate the main ideas and the structure of the proof.

The original free boundary problem (1.6)-(1.7) subject to conditions (1.9)-(1.11) is reformulated as an initial boundary value problem (2.8) with constraints (2.12)-(2.13) on a fixed domain through the Lagrangian particle trajectory framework (2.1) and (2.2). Notably, we derive novel estimates for the Riesz potential that address its inherent non-local structure. Unlike the Newtonian potential case analyzed in [24], where localization was feasible, we develop two direct treatments of the Riesz potential's non-local nature. One is the establishment of the estimate (3.10) of the linear operator $\mathcal{L}$ in (3.7) for the basic energy estimate, which also reveals the competition between the fluid pressure and the self-consistent attractie force in the energy (1.5) (see **Remark 3.1**). The other is a strategic decomposition (**Lemma 3.1, Lemma 3.2, Lemma 3.8**) for the higher spatial regularity. Our approach also involves critical aspect of establishing rigorous positive lower and upper bounds for the determinant of the Jacobian matrix associated with the Lagrangian coordinate transformation, specifically $\partial_x \eta$. Building on the methodology in [24], we investigate the representation formula of $\partial_x \eta$ in further details and we advance their analytical framework. Our key innovation lies in establishing a bootstrap pointwise estimate (**Lemma** 3.4), achieved through systematic applications of weighted Sobolev space techniques with carefully constructed test functions.

*Notations.* Throughout this paper, $\varepsilon$, $C_\varepsilon$ (or $C'_\varepsilon$), $C$ (or $C'$) and $C_{s,\gamma}$ (or $C'_{s,\gamma}$) will be used to denote generic positive constants. $L^p(-\bar{R}, \bar{R})$ $(1 \le p < \infty)$ denotes the $p$-th integrable Lebesgue space with norm defined by

$$\|f\|_{L^p} := \left( \int_{-\bar{R}}^{\bar{R}} |f|^p \, \mathrm{d}x \right)^{\frac{1}{p}}.$$

$L^\infty$ denotes the essentially bounded measurable functional space with norm defined by

$$\|f\|_{L^\infty} := \inf\{C \ge 0 : |f| \le C \text{ in } (-\bar{R}, \bar{R}) \text{ a.e.}\}.$$

$H^\alpha(-\bar{R}, \bar{R})$ $(\alpha \in \mathbb{N} \text{ or } \alpha \in (0,1))$ denotes the Sobolev space with norm defined by

$$\|f\|_{H^\alpha} := \|f\|_{L^2} + \sum_{n=1}^{\alpha} \|\partial_x^n f\|_{L^2}.$$

In addition, we use the following definition for the norm of fractional Sobolev spaces $H^\alpha(-\bar{R}, \bar{R})$ $\alpha \in (0,1)$:

$$\|f\|_{H^\alpha} := \|f\|_{L^2} + \int_{-\bar{R}}^{\bar{R}} \int_{-\bar{R}}^{\bar{R}} \frac{|f(x) - f(y)|^2}{|x-y|^{1+2\alpha}} \, \mathrm{d}y \, \mathrm{d}x.$$

We also introduce the concept of the Cauchy principal value: given a function $f$ over $[a, b]$ with a singularity point $a \le c \le b$, we define its Cauchy principal value over $[a, b]$ as

$$\text{p.v.} \int_a^b f(x) dx := \lim_{\varepsilon \to 0} \left( \int_a^{c-\varepsilon} f(x) dx + \int_{c+\varepsilon}^b f(x) dx \right).$$

## 3. A PRIORI ESTIMATES

In the rest of this work, we set $\bar{R}$ to be unity for notation simplicity. Let us denote the perturbation of $\eta$ around $x$ by

$$w(t,x) := \eta(t,x) - x.$$

for $t \ge 0$ and $-1 \le x \le 1$. Then by (2.8), $w$ satisfies

$$\bar{\rho}\partial_{tt}w + \partial_x \left[ \left( \frac{\bar{\rho}}{1 + \partial_x w} \right)^\gamma - \bar{\rho}^\gamma \right] + \Phi = \partial_x \left( \frac{\partial_{tx} w}{1 + \partial_x w} \right), \tag{3.1}$$

for $t > 0$ and $-1 < x < 1$. In addition, (2.12), (2.14), (2.2), and (2.13) imply that

$$\partial_x w(t, \pm 1) = 0, \quad w(t, 0) = 0, \tag{3.2}$$

$$(w, \partial_t w)(0, x) = (w_0, w_1)(x). \tag{3.3}$$

Here, $\Phi$ has the form given in (2.9):

$$\Phi(t,x) = \partial_x(\bar{\rho}^\gamma(x)) + \frac{\bar{\rho}(x)}{\partial_x \eta(t,x)} \int_{-\infty}^{\infty} W'(w(t,x) - \tilde{w}(t,y) + x - y)\Theta(t,x,y) \, \mathrm{d}y,$$

with $\Theta(t,x,y) = \bar{\rho}(y)\partial_x w(t,x) - \bar{\rho}(x)\partial_y \tilde{w}(t,y) + \bar{\rho}(y) - \bar{\rho}(x)$. We further define $\bar{\rho}\partial_t w_1$ by using (3.1) and (3.3),

$$\bar{\rho}(x)\partial_t w_1(x) := -\partial_x\left[\left(\frac{\bar{\rho}(x)}{1+\partial_x w(x)}\right)^\gamma - \bar{\rho}^\gamma(x)\right] - \Phi(0,x) + \partial_x\left(\frac{\partial_x w_1(x)}{1+\partial_x w_0(x)}\right).$$

In addition, we also denote $\tilde{w}$ the extension of $w$ over $\mathbb{R}$.

$$\tilde{w}(t,x) := \tilde{\eta}(t,x) - x,$$

for $t \geq 0$ and $x \in \mathbb{R}$. Moreover, we define the initial perturbation

$$\epsilon_0 := \sup_{-1 \leq x \leq 1}|\partial_x w_0| + \sup_{-1 \leq x \leq 1}|\partial_x w_1|.$$

We make the ***a priori*** assumption that there exists a suitably large constant $A > 0$ such that

$$\sup_{0 \leq t \leq T, -1 \leq x \leq 1}|\partial_x w(t,x)| \leq A\epsilon_0, \qquad \sup_{0 \leq t \leq T, -1 \leq x \leq 1}|\partial_{tx} w(t,x)| \leq A\epsilon_0, \qquad (3.4)$$

for $T > 0$. For the definition of $A$ we refer to (3.43) as introduced in the proof of Lemma 3.4 later on and note that $A$ is independent of $T$.

3.1. **Lower-order estimates.** In this subsection, we will show two results, formulated in **Lemma 3.1** and **Lemma 3.2**, concerning the basic energy structure for (2.8). Moreover a time-decay rate of the energy will be shown in **Lemma 3.3**. In the last part, we will provide a pointwise estimate in **Lemma 3.4**.

3.1.1. *Basic energy estimates.* First, let us give an estimate for the weighted norm of both time and spatial derivative of $w$.

**Lemma 3.1.** *For $\gamma > 2(1-s)$, there exists a $\epsilon_0 > 0$ such that under the **a priori** assumption* (3.4)*, it holds that*

$$\|\bar{\rho}^{\frac{1}{2}}\partial_t w(t)\|_{L^2}^2 + \|\bar{\rho}^{\frac{\gamma}{2}}\partial_x w(t)\|_{L^2}^2 + \int_0^t \|\partial_{tx} w(s)\|_{L^2}^2 \, \mathrm{d}s$$
$$\leq C_{s,\gamma}(\|\bar{\rho}^{\frac{1}{2}}w_1\|_{L^2}^2 + \|\bar{\rho}^{\frac{\gamma}{2}}\partial_x w_0\|_{L^2}^2), \quad 0 \leq t \leq T, \quad (3.5)$$

*for some constant $C_{s,\gamma} > 0$.*

*Proof.* Multiplying (3.1) by $\partial_t w$, then integrating the product with respect to spatial variable, using the integration by parts and (3.2), we have

$$\frac{1}{2}\frac{\mathrm{d}}{\mathrm{d}t}\int_{-1}^1 \bar{\rho}|\partial_t w|^2\,\mathrm{d}x + \underbrace{\frac{\mathrm{d}}{\mathrm{d}t}\int_{-1}^1 \frac{1}{\gamma-1}\bar{\rho}^\gamma\left[\left(\frac{1}{1+\partial_x w}\right)^{\gamma-1} + (\gamma-1)\partial_x w - 1\right]\,\mathrm{d}x}_{(3.6)_1}$$

$$+ \underbrace{\frac{1}{2}\frac{\mathrm{d}}{\mathrm{d}t}\int_{-1}^1 w\cdot\mathcal{L}w\,\mathrm{d}x}_{(3.6)_2} + \underbrace{\frac{1}{2}\frac{\mathrm{d}}{\mathrm{d}t}\int_{-1}^1\int_{-1}^1 R_{0,2}(t,x,y)\bar{\rho}(y)\bar{\rho}(x)\,\mathrm{d}y\,\mathrm{d}x}_{(3.6)_3}$$

$$+ \underbrace{\int_{-1}^1 \frac{|\partial_{tx} w|^2}{1+\partial_x w}\,\mathrm{d}x}_{(3.6)_4} = 0, \quad (3.6)$$

where $\mathcal{L}$ is a linear operator defined by

$$\mathcal{L}w(t,x) := \mathcal{L}_1 w(t,x) + \mathcal{L}_2 w(t,x), \tag{3.7}$$

$$\mathcal{L}_1 w(t,x) := \bar{\rho}^2(x) \int_{-\infty}^{\infty} \frac{x-y}{|x-y|^{3-2s}} (\partial_x \tilde{w}(t,x) - \partial_y \tilde{w}(t,y)) \, \mathrm{d}y,$$

$$\mathcal{L}_2 w(t,x) := -\bar{\rho}(x) \int_{-\infty}^{\infty} \frac{2(s-1)}{|x-y|^{3-2s}} (\bar{\rho}(x) - \bar{\rho}(y))(\tilde{w}(t,x) - \tilde{w}(t,y)) \, \mathrm{d}y,$$

and the nonlinear term $R_{0,2}$ has the following form,

$$R_{0,2}(t,x,y) := \frac{1}{2s-1} \frac{1}{|\eta(t,x) - \eta(t,y)|^{1-2s}} - \frac{1}{2s-1} \frac{1}{|x-y|^{1-2s}}$$
$$- \frac{x-y}{|x-y|^{3-2s}} (w(t,x) - w(t,y)) - \frac{s-1}{|x-y|^{3-2s}} (w(t,x) - w(t,y))^2. \tag{3.8}$$

For $(3.6)_1$, by the estimate $(3.4)$ in **a priori assumption**, we have

$$\left(\frac{\gamma}{2} - A\epsilon_0\right) \int_{-1}^{1} \bar{\rho}^\gamma |\partial_x w|^2 \, \mathrm{d}x \leq (3.6)_1 \leq \gamma \int_{-1}^{1} \bar{\rho}^\gamma |\partial_x w|^2 \, \mathrm{d}x, \tag{3.9}$$

where we used the following elementary algebraic inequality,

$$\left(\frac{\gamma}{2} - \epsilon\right) z^2 \leq \frac{1}{\gamma-1} \left(\frac{1}{1+z}\right)^{\gamma-1} + z - \frac{1}{\gamma-1} \leq \gamma z^2,$$

for $|z| < \epsilon$ with $\epsilon$ sufficiently small, that gives a restriction for $\epsilon_0$ in $(3.4)$.

For $(3.6)_2$ we claim that

$$|(3.6)_2| \leq 2(1-s) \int_{-1}^{1} \bar{\rho}^\gamma |\partial_x w|^2 \, \mathrm{d}x. \tag{3.10}$$

Indeed, recall $(3.7)$, on the one hand, by direct computation, we find

$$\int_{-1}^{1} w \cdot \mathcal{L}w \, \mathrm{d}x = -2(1-s) \iint_{-1 < x < y < 1} \frac{|w(x) - w(y)|^2}{|x-y|^{3-2s}} \bar{\rho}(x)\bar{\rho}(y) \, \mathrm{d}x \, \mathrm{d}y. \tag{3.11}$$

On the other hand, as $\bar{\rho}$ has compactly supports and $\bar{\rho}$ is the solution to $(1.4)$, then

$$\bar{\rho}^\gamma(x) = \int_{x}^{\infty} \int_{-\infty}^{\infty} \frac{y-z}{|y-z|^{3-2s}} \bar{\rho}(y)(\bar{\rho}(z) - \bar{\rho}(y)) \, \mathrm{d}z \, \mathrm{d}y,$$

for any $x \in [0,1]$. Since for arbitrary $\varepsilon > 0$,

$$\left| \int_{y-\varepsilon}^{y+\varepsilon} \frac{y-z}{|y-z|^{3-2s}} \bar{\rho}(y)(\bar{\rho}(z) - \bar{\rho}(y)) \, \mathrm{d}z \right| \leq \bar{\rho}(y) \|\bar{\rho}\|_{W^{1,\infty}} \int_{y-\varepsilon}^{y+\varepsilon} \frac{1}{|y-z|^{1-2s}} \, \mathrm{d}z$$
$$\leq \frac{1}{s} \bar{\rho}(y) \|\bar{\rho}\|_{W^{1,\infty}} \varepsilon^{2s},$$

we have

$$\lim_{\varepsilon \to 0} \int_{y-\varepsilon}^{y+\varepsilon} \frac{y-z}{|y-z|^{3-2s}} \bar{\rho}(y)(\bar{\rho}(z) - \bar{\rho}(y)) \, \mathrm{d}z = 0. \tag{3.12}$$

Noticing $(y-z)/|y-z|^{3-2s}$ is odd symmetric centered at $y$, we deduce that

$$\left( \int_{-\infty}^{y-\varepsilon} + \int_{y+\varepsilon}^{\infty} \right) \frac{y-z}{|y-z|^{3-2s}} \bar{\rho}(y)^2 \, \mathrm{d}z = 0,$$

that together with (3.12) implies

$$\bar{\rho}^\gamma(x) = \int_x^\infty \text{p.v.} \int_{-\infty}^\infty \frac{y-z}{|y-z|^{3-2s}} \bar{\rho}(y)\bar{\rho}(z) \, \mathrm{d}z \, \mathrm{d}y.$$

For notational simplicity, we will omit p.v. in the integrals hereafter. Therefore, it follows that

$$\int_{-1}^1 \bar{\rho}^\gamma |\partial_x w|^2 \, \mathrm{d}x = \int_{-\infty}^\infty \int_x^\infty \int_{-\infty}^\infty \frac{y-z}{|y-z|^{3-2s}} \bar{\rho}(y)\bar{\rho}(z) |\partial_x \tilde{w}(x)|^2 \, \mathrm{d}z \, \mathrm{d}y \, \mathrm{d}x$$
$$= \int_{-\infty}^\infty \int_{-\infty}^y \int_{-\infty}^\infty \frac{y-z}{|y-z|^{3-2s}} \bar{\rho}(y)\bar{\rho}(z) |\partial_x \tilde{w}(x)|^2 \, \mathrm{d}z \, \mathrm{d}x \, \mathrm{d}y,$$

where we have used the Fubini's theorem in the last line. By the symmetry of integration in $y$ and $z$, we also have

$$\int_{-1}^1 \bar{\rho}^\gamma |\partial_x w|^2 \, \mathrm{d}x = \int_{-\infty}^\infty \int_{-\infty}^z \int_{-\infty}^\infty \frac{z-y}{|z-y|^{3-2s}} \bar{\rho}(z)\bar{\rho}(y) |\partial_x \tilde{w}(x)|^2 \, \mathrm{d}y \, \mathrm{d}x \, \mathrm{d}z,$$

then

$$\int_{-1}^1 \bar{\rho}^\gamma |\partial_x w|^2 \, \mathrm{d}x = \frac{1}{2} \int_{-\infty}^\infty \int_{-\infty}^\infty \int_z^y \frac{y-z}{|y-z|^{3-2s}} \bar{\rho}(y)\bar{\rho}(z) |\partial_x \tilde{w}(x)|^2 \, \mathrm{d}x \, \mathrm{d}y \, \mathrm{d}z$$
$$= \iint_{-\infty < z < y < \infty} \frac{y-z}{|y-z|^{3-2s}} \bar{\rho}(y)\bar{\rho}(z) \int_z^y |\partial_x \tilde{w}(x)|^2 \, \mathrm{d}x \, \mathrm{d}y \, \mathrm{d}z. \qquad (3.13)$$

By the fundamental theorem of calculus, for $y > z$, we have

$$|\tilde{w}(y) - \tilde{w}(z)| \le \left( \int_z^y |\partial_x \tilde{w}(x)|^2 \, \mathrm{d}x \right)^{1/2} (y-z)^{1/2},$$

which together with (3.13) implies

$$\int_{-1}^1 \bar{\rho}^\gamma |\partial_x w|^2 \, \mathrm{d}x \ge \iint_{-1 < z < y < 1} \bar{\rho}(y)\bar{\rho}(z) \frac{|w(y) - w(z)|^2}{|y-z|^{3-2s}} \, \mathrm{d}y \, \mathrm{d}z. \qquad (3.14)$$

For $(3.6)_3$, let us firstly prove a pointwise estimate for $R_{0,2}$,

$$|R_{0,2}(t,x,y)| \le 4(1-s)(3-2s)A\epsilon_0 \frac{|w(x) - w(y)|^2}{|x-y|^{3-2s}}, \qquad (3.15)$$

then immediately we have

$$|(3.6)_3| \le 4(1-s)(3-2s)A\epsilon_0 \int_{-1}^1 \int_{-1}^1 \frac{|w(x) - w(y)|^2}{|x-y|^{3-2s}} \bar{\rho}(x)\bar{\rho}(y) \, \mathrm{d}y \, \mathrm{d}x.$$

Next using (3.14) we find that

$$|(3.6)_3| \le 8(1-s)(3-2s)A\epsilon_0 \int_{-1}^1 \bar{\rho}^\gamma |\partial_x w|^2 \, \mathrm{d}x. \qquad (3.16)$$

It remains to prove (3.15). We compute

$$\frac{\mathrm{d}}{\mathrm{d}\alpha} \left( \frac{1}{2s-1} \frac{1}{|\alpha(\eta(t,x) - \eta(t,y)) + (1-\alpha)(x-y)|^{1-2s}} \right)$$
$$= \frac{x-y + \alpha(w(t,x) - w(t,y))}{|x-y + \alpha(w(t,x) - w(t,y))|^{3-2s}} (w(t,x) - w(t,y)), \qquad (3.17)$$

and

$$\frac{\mathrm{d}}{\mathrm{d}\beta}\left(\frac{x-y+\beta(w(t,x)-w(t,y))}{|x-y+\beta(w(t,x)-w(t,y))|^{3-2s}}\right)$$

$$=2(s-1)\frac{w(t,x)-w(t,y)}{|x-y+\beta(w(t,x)-w(t,y))|^{3-2s}}$$

$$+\delta_0(x-y+\beta(w(t,x)-w(t,y)))\frac{w(t,x)-w(t,y)}{|x-y+\beta(w(t,x)-w(t,y))|^{2-2s}}. \quad (3.18)$$

Applying the fundamental theorem of calculus to (3.8), taking (3.17) and (3.18) into account, we have

$$|R_{0,2}(t,x,y)|$$

$$=\left|\int_0^1\frac{x-y+\alpha(w(t,x)-w(t,y))}{|x-y+\alpha(w(t,x)-w(t,y))|^{3-2s}}(w(t,x)-w(t,y))\,\mathrm{d}\alpha\right.$$

$$\left.-\frac{x-y}{|x-y|^{3-2s}}(w(t,x)-w(t,y))-\frac{2s-2}{2}\frac{1}{|x-y|^{3-2s}}(w(t,x)-w(t,y))^2\right|$$

$$=\frac{|w(t,x)-w(t,y)|^2}{|x-y|^{3-2s}}$$

$$\times\left|\int_0^1\int_0^\alpha 2(s-1)\frac{1}{|1+\beta\frac{w(t,x)-w(t,y)}{x-y}|^{3-2s}}\,\mathrm{d}\beta\,\mathrm{d}\alpha-(s-1)\right|$$

$$=2(1-s)\frac{|w(t,x)-w(t,y)|^2}{|x-y|^{3-2s}}\left|\int_0^1\frac{1-\beta}{|1+\beta\frac{w(t,x)-w(t,y)}{x-y}|^{3-2s}}-(1-\beta)\,\mathrm{d}\beta\right|,$$

where we use the **a priori assumption** (3.4) such that the second term on the right hand side of (3.18) vanishes since $|\frac{x-y}{w(t,x)-w(t,y)}|\geq (A\epsilon_0)^{-1}\gg 1$ if $\epsilon_0$ is sufficiently small. Moreover, by **a priori assumption** (3.4) again, then we have

$$\left|\frac{1}{|1+\beta\frac{w(t,x)-w(t,y)}{x-y}|^{3-2s}}-1\right|\leq 2(3-2s)A\epsilon_0,$$

where we used the following elementary algebraic inequality,

$$\left|\frac{1}{|1+z|^{3-2s}}-1\right|\leq 2(3-2s)\epsilon,$$

for $|z|<\epsilon$ with $\epsilon$ sufficiently small. Therefore

$$|R_{0,2}(t,x,y)|\leq 4(1-s)(3-2s)A\epsilon_0\frac{|w(x)-w(y)|^2}{|x-y|^{3-2s}}.$$

For $(3.6)_4$, by the estimate (3.4) in **a priori assumption**, we have

$$(3.6)_4\geq\frac{1}{1+A\epsilon_0}\int_0^t\int_{-1}^1|\partial_{tx}w|^2\,\mathrm{d}x\,\mathrm{d}s, \quad (3.19)$$

Therefore combining the results of (3.9), (3.10), (3.16), (3.19) and noticing that $\gamma>2(1-s)$, then integrating the resultant with respect to time $t$, by choosing $\epsilon_0$ sufficiently small we obtain (3.5) which proves **Lemma 3.1**.          □

**Remark 3.1.** *In fact, the estimate between two terms* $(3.6)_1$ *and* $(3.6)_2$ *reveals the competition between the fluid pressure and the self-consistent attractive force. Let us illustrate it from a variational point of view.*

*We recall from* (1.5) *that in the one-dimensional case,*

$$\mathcal{F}[\rho] = \int_{\mathbb{R}} \frac{1}{\gamma - 1} \rho^\gamma + \frac{1}{2} \rho(W * \rho) \, \mathrm{d}x,$$

*then* $\mathcal{F}$ *can be rewritten in the Lagrangian coordinates and use the mass conservation* (2.7) *to get*

$$\mathcal{F}[\rho] = \frac{1}{\gamma - 1} \int_{-1}^{1} \bar{\rho}^\gamma (\partial_x \eta)^{1-\gamma} \, \mathrm{d}x$$
$$+ \frac{1}{2} \int_{-1}^{1} \int_{-1}^{1} W(\eta(t,x) - \eta(t,y)) \bar{\rho}(x) \bar{\rho}(y) \, \mathrm{d}x \, \mathrm{d}y =: \mathcal{H}[\eta].$$

*We set* $\eta = x + \epsilon w$. *Formally computing,*

$$(\partial_x \eta)^{1-\gamma} = (1 + \epsilon \partial_x w)^{1-\gamma} = 1 + \epsilon(1-\gamma)\partial_x w + \frac{1}{2}(1-\gamma)(-\gamma)(\epsilon \partial_x w)^2 + O(1)\epsilon^3,$$

*and*

$$W(\eta(t,x) - \eta(t,y)) = W(x - y + \epsilon(w(x) - w(y)))$$
$$= W(x - y) + W'(x-y)\epsilon[w(x) - w(y)]$$
$$+ \frac{1}{2}W''(x-y)\epsilon^2[w(x) - w(y)]^2 + O(1)\epsilon^3.$$

*Since the unique global minimizer of the nonlinear energy functional* $\mathcal{F}$ *is attained at* $\bar{\rho}$, *we obtain*

$$\mathcal{H}[x + \epsilon w] = \mathcal{H}(x)$$
$$+ \epsilon \left\{ -\int_{-1}^{1} \bar{\rho}^\gamma \partial_x w \, \mathrm{d}x + \frac{1}{2} \int_{-1}^{1} \int_{-1}^{1} W'(x-y)[w(x) - w(y)] \bar{\rho}(x) \bar{\rho}(y) \, \mathrm{d}y \, \mathrm{d}x \right\}$$
$$+ \frac{1}{2}\epsilon^2 \left\{ \gamma \int_{-1}^{1} \bar{\rho}^\gamma |\partial_x w|^2 \, \mathrm{d}x - 2(1-s) \iint_{x<y} \frac{|w(x) - w(y)|^2}{|x-y|^{3-2s}} \bar{\rho}(x) \bar{\rho}(y) \, \mathrm{d}y \, \mathrm{d}x \right\} + O(1)\epsilon^3.$$

*For the coefficient of* $\epsilon$ *above, using the stationary equation* (1.4), *we have*

$$-\int_{-1}^{1} \bar{\rho}^\gamma \partial_x w \, \mathrm{d}x + \frac{1}{2} \int_{-1}^{1} \int_{-1}^{1} W'(x-y)[w(x) - w(y)] \bar{\rho}(x) \bar{\rho}(y) \, \mathrm{d}y \, \mathrm{d}x$$
$$= \int_{-1}^{1} \partial_x (\bar{\rho}^\gamma) \cdot w \, \mathrm{d}x + \int_{-1}^{1} \int_{-1}^{1} W'(x-y)w(x)\bar{\rho}(x)(\bar{\rho}(y) - \bar{\rho}(x)) \, \mathrm{d}y \, \mathrm{d}x = 0.$$

*Therefore, since* $\mathcal{H}[x + \epsilon w]$ *attains its unique minimum at* $\epsilon = 0$, *it holds at least formally that*

$$\gamma \int_{-1}^{1} \bar{\rho}^\gamma |\partial_x w|^2 \, \mathrm{d}x \geq 2(1-s) \iint_{-1<x<y<1} \frac{|w(x) - w(y)|^2}{|x-y|^{3-2s}} \bar{\rho}(x) \bar{\rho}(y) \, \mathrm{d}y \, \mathrm{d}x. \qquad (3.20)$$

*Since* $\gamma > 2(1-s)$, *both* (3.14) *and* (3.20) *explain that the fluid pressure dominates the self-consistent attractive force.*

In the following, we improved the estimate in **Lemma 3.1** for the spatial derivative of $w$ near the vacuum.

**Lemma 3.2.** *For $2(1-s) < \gamma \leq 2$, there exists an $\epsilon_0 > 0$ such that under the **a prior assumption** (3.4), the following holds*

$$\|\partial_x w(t)\|_{L^2}^2 + \int_0^t \|\bar{\rho}^{\frac{\gamma}{2}} \partial_x w(s)\|_{L^2}^2 \, \mathrm{d}s$$
$$\leq C_{s,\gamma}(\|\partial_x w_0\|_{L^2}^2 + \|\bar{\rho}^{\frac{1}{2}} w_1\|_{L^2}^2 + \|\bar{\rho}^{\frac{\gamma}{2}} \partial_x w_1\|_{L^2}^2), \quad 0 \leq t \leq T, \quad (3.21)$$

*for some constant $C_{s,\gamma} > 0$.*

*Proof.* We rewrite (3.1) as

$$-\partial_x \left( \frac{\partial_{tx} w}{1 + \partial_x w} \right) + \partial_x \left[ \left( \frac{\bar{\rho}}{1 + \partial_x w} \right)^{\gamma} - \bar{\rho}^{\gamma} \right] = -\bar{\rho} \partial_{tt} w - \Phi, \quad (3.22)$$

then multiplying (3.22) by $w$, integrating the product with respect to spatial variable, using the integration by parts and boundary conditions (3.2),

$$\frac{\mathrm{d}}{\mathrm{d}t} \underbrace{\int_{-1}^1 \partial_x w - \ln(1 + \partial_x w) \, \mathrm{d}x}_{(3.23)_1} + \underbrace{\int_{-1}^1 \bar{\rho}^{\gamma} \partial_x w \left[ 1 - \left( \frac{1}{1 + \partial_x w} \right)^{\gamma} \right] \mathrm{d}x}_{(3.23)_2}$$

$$= \frac{\mathrm{d}}{\mathrm{d}t} \underbrace{\int_{-1}^1 -\bar{\rho} w \partial_t w \, \mathrm{d}x}_{(3.23)_3} + \underbrace{\int_{-1}^1 \bar{\rho} |\partial_t w|^2 \, \mathrm{d}x}_{(3.23)_4} + \underbrace{\int_{-1}^1 -w \cdot \mathcal{L}w \, \mathrm{d}x}_{(3.23)_5}$$

$$+ \frac{1}{2} \underbrace{\int_{-1}^1 \int_{-1}^1 -\bar{\rho}(x) \bar{\rho}(y) (w(t,x) - w(t,y)) R_{1,1}(t,x,y) \, \mathrm{d}y \, \mathrm{d}x}_{(3.23)_6}, \quad (3.23)$$

where the nonlinear term $R_{1,1}$ has the form

$$R_{1,1}(t,x,y) := \frac{\eta(t,x) - \eta(t,y)}{|\eta(t,x) - \eta(t,y)|^{3-2s}} - \frac{x-y}{|x-y|^{3-2s}} - 2(s-1) \frac{w(t,x) - w(t,y)}{|x-y|^{3-2s}}. \quad (3.24)$$

For $(3.23)_1$, by **a priori** assumption (3.4), we have

$$\frac{1}{4} \int_{-1}^1 |\partial_x w|^2 \, \mathrm{d}x \leq (3.23)_1 \leq \int_{-1}^1 |\partial_x w|^2 \, \mathrm{d}x, \quad (3.25)$$

where we used the following elementary algebraic inequality

$$\frac{1}{4} z^2 \leq z - \ln(1+z) \leq z^2,$$

for $|z| < \epsilon$ with $\epsilon$ sufficiently small.

For $(3.23)_2$, by **a priori** assumption (3.4), we have

$$(\gamma - 3A\epsilon_0) \int_{-1}^1 \bar{\rho}^{\gamma} |\partial_x w|^2 \, \mathrm{d}x \leq (3.23)_2 \leq \frac{\gamma}{2} \int_{-1}^1 \bar{\rho}^{\gamma} |\partial_x w|^2 \, \mathrm{d}x, \quad (3.26)$$

where we used the following elementary algebraic inequality

$$(\gamma - 3\epsilon) z^2 \leq z \left[ 1 - \left( \frac{1}{1+z} \right)^{\gamma} \right] \leq 2\gamma z^2,$$

for $|z| < \epsilon$ with $\epsilon$ sufficiently small.

For $(3.23)_3$, by Young's inequality and the fundamental theorem of calculus using that $w(t,0) = 0$, we have

$$|(3.23)_3| \leq \varepsilon \int_{-1}^{1} |\partial_x w|^2 \, \mathrm{d}x + C_\varepsilon \int_{-1}^{1} \bar{\rho} |\partial_t w|^2 \, \mathrm{d}x, \tag{3.27}$$

where we have used the fact that $\bar{\rho} \in L^\infty$ for $2(1-s) < \gamma \leq 2$.

For $(3.23)_4$, by the fundamental theorem of calculus using that $\partial_t w(t,0) = 0$, we have

$$|(3.23)_4| \leq C \int_{-1}^{1} |\partial_{tx} w|^2 \, \mathrm{d}x. \tag{3.28}$$

For $(3.23)_5$, we have the same estimate as $(3.14)$,

$$|(3.23)_5| \leq 2(1-s) \int_{-1}^{1} \bar{\rho}^\gamma |\partial_x w|^2 \, \mathrm{d}x. \tag{3.29}$$

For $(3.23)_6$, we claim that

$$|R_{1,1}(t,x,y)| \leq 2(1-s)(3-2s) \frac{A^2 \epsilon_0^2}{|x-y|^{2-2s}}, \tag{3.30}$$

then we have

$$|(3.23)_6| \leq A^2 \epsilon_0^2 (1-s)(3-2s) \int_{-1}^{1} \int_{-1}^{1} \frac{|w(t,x) - w(t,y)|}{|x-y|^{3-2s}} \bar{\rho}(x) \bar{\rho}(y) \, \mathrm{d}y \, \mathrm{d}x.$$

Using the same estimate as $(3.14)$, we can bound $(3.23)_6$ as

$$|(3.23)_6| \leq A^2 \epsilon_0^2 (1-s)(3-2s) \int_{-1}^{1} \bar{\rho}^\gamma |\partial_x w|^2 \, \mathrm{d}x. \tag{3.31}$$

It remains to prove $(3.30)$. We compute

$$\frac{\mathrm{d}}{\mathrm{d}\alpha} \frac{\alpha(\eta(t,x) - \eta(t,y)) + (1-\alpha)(x-y)}{|\alpha(\eta(t,x) - \eta(t,y)) + (1-\alpha)(x-y)|^{3-2s}}$$
$$= 2(s-1) \frac{w(t,x) - w(t,y)}{|x-y + \alpha(w(t,x) - w(t,y))|^{3-2s}}$$
$$+ \delta_0 (x - y + \alpha(w(t,x) - w(t,y))) \frac{w(t,x) - w(t,y)}{|x-y + \alpha(w(t,x) - w(t,y))|^{2-2s}}. \tag{3.32}$$

Applying the fundamental theorem of calculus to $(3.24)$, taking $(3.32)$ into account, we have

$$|R_{1,1}(t,x,y)| = 2(1-s) \frac{|w(t,x) - w(t,y)|}{|x-y|^{3-2s}} \left| \int_{0}^{1} 1 - \frac{1}{|1 + \alpha \frac{w(t,x) - w(t,y)}{x-y}|^{3-2s}} \, \mathrm{d}\alpha \right|,$$

where we use the **a prior** assumption $(3.4)$ such that the second term on the right hand side of $(3.32)$ vanishes since $|\frac{x-y}{w(t,x) - w(t,y)}| \geq (A\epsilon_0)^{-1} \gg 1$ if $\epsilon_0$ is sufficiently small. Moreover, by **a prior** assumption $(3.4)$ again, we have

$$\left| 1 - \frac{1}{|1 + \alpha \frac{w(t,x) - w(t,y)}{x-y}|^{3-2s}} \right| \leq 2(3-2s) A\epsilon_0,$$

where we used the following elementary algebraic inequality,

$$\left| 1 - \frac{1}{|1+z|^{3-2s}} \right| \leq 2(3-2s)\epsilon,$$

for $|z| < \epsilon$ with $\epsilon$ sufficiently small. Therefore

$$|R_{1,1}(t,x,y)| \leq 2(1-s)(3-2s)\frac{A^2 \epsilon_0^2}{|x-y|^{2-2s}}.$$

Therefore combining the results of (3.25), (3.26), (3.27), (3.28), (3.29), (3.31) and noticing that $\gamma > 2(1-s)$, then integrating the resultant with respect to time $t$, by choosing and $\epsilon_0$ sufficiently small, we obtain (3.21) which proves **Lemma 3.2**. $\qquad\square$

Having obtained basic energy estimates, we further deduce a time-weighted estimate for the time and spatial derivatives of $w$.

**Lemma 3.3.** *For $2(1-s) < \gamma < 2$, there exists an $\epsilon_0 > 0$ such that under the **a prior assumption** (3.4) the following holds*

$$(1+t)\|\bar{\rho}^{\frac{1}{2}}\partial_t w(t)\|_{L^2}^2 + (1+t)\|\bar{\rho}^{\frac{\gamma}{2}}\partial_x w(t)\|_{L^2}^2 + \int_0^t (1+s)\|\partial_{tx} w(s)\|_{L^2}^2 \, ds$$

$$\leq C_{s,\gamma}(\|\partial_x w_0\|_{L^2}^2 + \|\bar{\rho}^{\frac{1}{2}} w_1\|_{L^2}^2 + \|\bar{\rho}^{\frac{\gamma}{2}}\partial_x w_0\|_{L^2}^2), \quad 0 \leq t \leq T, \quad (3.33)$$

*for some constant $C_{s,\gamma} > 0$.*

*Proof.* Recall that from **Lemma 3.1**, we have

$$\frac{1}{2}\frac{d}{dt}\int_{-1}^1 \bar{\rho}|\partial_t w|^2 \, dx + \frac{d}{dt}\int_{-1}^1 \frac{1}{\gamma-1}\bar{\rho}^\gamma \left[ \left(\frac{1}{1+\partial_x w}\right)^{\gamma-1} + (\gamma-1)\partial_x w - 1 \right] dx$$

$$+ \frac{1}{2}\frac{d}{dt}\int_{-1}^1 w \cdot \mathcal{L}w \, dx + \frac{d}{dt}\int_{-1}^1 \int_{-1}^1 \bar{\rho}(x)\bar{\rho}(y)R_{0,2}(t,x,y) \, dx \, dy$$

$$+ \int_{-1}^1 \frac{|\partial_{tx} w|^2}{1+\partial_x w} \, dx = 0,$$

Multiplying the above equation by $1+t$ we find that

$$\frac{1}{2}\frac{d}{dt}\left[ (1+t)\int_{-1}^1 \bar{\rho}|\partial_t w|^2 \, dx \right] - \frac{1}{2}\underbrace{\int_{-1}^1 \bar{\rho}|\partial_t w|^2 \, dx}_{(3.34)_1}$$

$$+ \frac{d}{dt}\left\{ \underbrace{(1+t)\int_{-1}^1 \frac{1}{\gamma-1}\bar{\rho}^\gamma \left[ \left(\frac{1}{1+\partial_x w}\right)^{\gamma-1} + (\gamma-1)\partial_x w - 1 \right] dx}_{(3.34)_2} \right\}$$

$$- \underbrace{\int_{-1}^1 \frac{1}{\gamma-1}\bar{\rho}^\gamma \left[ \left(\frac{1}{1+\partial_x w}\right)^{\gamma-1} + (\gamma-1)\partial_x w - 1 \right] dx}_{(3.34)_3}$$

$$+ \frac{1}{2} \frac{\mathrm{d}}{\mathrm{d}t} \Bigg[ \underbrace{(1+t) \int_{-1}^{1} w \cdot \mathcal{L}w \, \mathrm{d}x}_{(3.34)_4} \Bigg] - \frac{1}{2} \underbrace{\int_{-1}^{1} w \cdot \mathcal{L}w \, \mathrm{d}x}_{(3.34)_5}$$

$$+ \frac{1}{2} \frac{\mathrm{d}}{\mathrm{d}t} \Bigg( \underbrace{(1+t) \int_{-1}^{1} \int_{-1}^{1} \bar{\rho}(x)\bar{\rho}(y) R_{0,2}(t,x,y) \, \mathrm{d}y \, \mathrm{d}x}_{(3.34)_6} \Bigg)$$

$$- \frac{1}{2} \underbrace{\int_{-1}^{1} \int_{-1}^{1} \bar{\rho}(x)\bar{\rho}(y) R_{0,2}(t,x,y) \, \mathrm{d}y \, \mathrm{d}x}_{(3.34)_7} + \underbrace{(1+t) \int_{-1}^{1} \frac{|\partial_{tx} w|^2}{1 + \partial_x w} \, \mathrm{d}x}_{(3.34)_8} = 0. \quad (3.34)$$

For $(3.34)_1$, by the fundamental theorem of calculus using that $\partial_t w(t,0) = 0$, we have

$$|(3.34)_1| \leq C \int_{-1}^{1} |\partial_{tx} w|^2 \, \mathrm{d}x. \quad (3.35)$$

For $(3.34)_2$ and $(3.34)_3$, by the estimate $(3.4)$ in **a priori** assumption,

$$\left( \frac{\gamma}{2} - A\epsilon_0 \right)(1+t) \int_{-1}^{1} \bar{\rho}^{\gamma} |\partial_x w|^2 \, \mathrm{d}x \leq (3.34)_2 \leq \gamma(1+t) \int_{-1}^{1} \bar{\rho}^{\gamma} |\partial_x w|^2 \, \mathrm{d}x, \quad (3.36)$$

where we used the following elementary algebraic inequality,

$$\left( \frac{\gamma}{2} - \epsilon \right) z^2 \leq \frac{1}{\gamma - 1} \left( \frac{1}{1+z} \right)^{\gamma - 1} + z - \frac{1}{\gamma - 1} \leq \gamma z^2,$$

for $|z| < \epsilon$ with $\epsilon$ sufficiently small. We find also that

$$|(3.34)_3| \leq C \int_{-1}^{1} \bar{\rho}^{\gamma} |\partial_x w|^2 \, \mathrm{d}x. \quad (3.37)$$

For $(3.34)_4$ and $(3.34)_5$, using $(3.11)$ and $(3.14)$ we obtain

$$|(3.34)_4| \leq 2(1-s)(1+t) \int_{-1}^{1} \bar{\rho}^{\gamma} |\partial_x w|^2 \, \mathrm{d}x, \quad (3.38)$$

$$|(3.34)_5| \leq 2(1-s) \int_{-1}^{1} \bar{\rho}^{\gamma} |\partial_x w|^2 \, \mathrm{d}x. \quad (3.39)$$

For $(3.34)_6$ and $(3.34)_7$, similar as the estimate for $(3.6)_3$ in **Lemma 3.1**,

$$|(3.34)_6| \leq 8(1-s)(3-2s)A\epsilon_0(1+t) \int_{-1}^{1} \bar{\rho}^{\gamma} |\partial_x w|^2 \, \mathrm{d}x, \quad (3.40)$$

and

$$|(3.34)_7| \leq 8(1-s)(3-2s)A\epsilon_0 \int_{-1}^{1} \bar{\rho}^{\gamma} |\partial_x w|^2 \, \mathrm{d}x. \quad (3.41)$$

For $(3.34)_8$, by the estimate $(3.4)$ in **a priori** assumption,

$$\frac{1}{1 + A\epsilon_0}(1+t) \int_{-1}^{1} |\partial_{tx} w|^2 \, \mathrm{d}x \leq (3.34)_8 \leq \frac{1}{1 - A\epsilon_0}(1+t) \int_{-1}^{1} |\partial_{tx} w|^2 \, \mathrm{d}x. \quad (3.42)$$

Therefore combining the results of (3.35), (3.36), (3.37), (3.38), (3.39), (3.40), (3.41), (3.42) and noticing that $\gamma > 2(1 - s)$, then integrating with respect to $t$, by choosing $\epsilon_0$ sufficiently small, we obtain (3.33) which prove **Lemma 3.3**.    $\square$

3.1.2. *Pointwise estimates on Lagrangian variable.* Let us give a pointwise estimate on $\eta_x(t, x)$ for all $0 \le t \le T$. The key ingredient in the following lemma is that this estimate holds for arbitrary range of time which extends the estimate (3.4) in **a prior assumption** in time.

**Lemma 3.4.** *For $\gamma > 2(1 - s)$, let $N \in \mathbb{N}$ and $A > 0$ be large enough such that*

$$\begin{cases} (1 - 10^{-N})\gamma > 2(1 - s), \\ A \ge \dfrac{C_{s,\gamma}}{1 - 10^{-N} - 2(1 - s)/\gamma}, \end{cases} \tag{3.43}$$

*for a fixed constant $C_{s,\gamma} > 0$ determined below, then under the **a prior assumption** (3.4), it holds that*

$$\|\partial_x w(t)\|_{L^\infty} \le (1 - 10^{-N-1})A\epsilon_0, \quad 0 \le t \le T. \tag{3.44}$$

*Proof.* Rewriting the equation (2.8) as

$$\partial_{tx}(\ln \partial_x \eta) = \bar{\rho}\partial_{tt}\eta + \partial_x \left[ \left( \frac{\bar{\rho}}{\partial_x \eta} \right)^\gamma - \bar{\rho}^\gamma \right] + \Phi. \tag{3.45}$$

From the initial condition (2.2) and boundary conditions (2.12), by fundamental theorem of calculus, we have

$$\partial_x \eta(t, x) = \partial_x \eta_0(x) \exp \left[ -\int_0^t \int_x^1 \partial_{tx}(\ln \partial_x \eta) \, \mathrm{d}y \, \mathrm{d}\tau \right],$$

therefore by (3.45),

$$\partial_x \eta(t, x) = \partial_x \eta_0(x) \exp \left( F(t, x) - \int_0^t \int_x^1 G(\tau, y) \, \mathrm{d}y \, \mathrm{d}\tau \right), \tag{3.46}$$

where we set

$$F := \int_0^t \left( \frac{\bar{\rho}}{\partial_x \eta} \right)^\gamma \mathrm{d}\tau,$$
$$G := \bar{\rho}\partial_{tt}\eta - \partial_x(\bar{\rho}^\gamma) + \Phi.$$

Using again (3.46), we note that

$$\partial_t F = \left( \frac{\bar{\rho}}{\partial_x \eta} \right)^\gamma = \left( \frac{\bar{\rho}}{\partial_x \eta_0} \right)^\gamma \exp(-\gamma F) \exp \left( \gamma \int_0^t \int_x^1 G(\tau, y) \, \mathrm{d}y \, \mathrm{d}\tau \right).$$

Multiplying the above resultant by $\gamma \exp(\gamma F)$ gives

$$\partial_t \exp(\gamma F) = \gamma \left( \frac{\bar{\rho}}{\partial_x \eta_0} \right)^\gamma \exp \left( \gamma \int_0^t \int_x^1 G(\tau, y) \, \mathrm{d}y \, \mathrm{d}\tau \right).$$

Taking integration in time over $[0, t]$ and noting $F(0, x) = 0$ it follows that

$$\exp(F) = \left[ 1 + \gamma \left( \frac{\bar{\rho}}{\partial_x \eta_0} \right)^\gamma \int_0^t \exp \left( \gamma \int_0^s \int_x^1 G(\tau, y) \, \mathrm{d}y \, \mathrm{d}\tau \right) \mathrm{d}s \right]^{1/\gamma}.$$

Recall (3.46), we have

$$\partial_x \eta(t,x) = \partial_x \eta_0(x) \exp\left(-\int_0^t \int_x^1 G(\tau,y)\,\mathrm{d}y\,\mathrm{d}\tau\right)$$

$$\times \left[1 + \gamma\left(\frac{\bar\rho(x)}{\partial_x \eta_0(x)}\right)^\gamma \int_0^t \exp\left(\gamma \int_0^s \int_x^1 G(\tau,y)\,\mathrm{d}y\,\mathrm{d}\tau\right)\mathrm{d}s\right]^{1/\gamma}$$

$$= \partial_x \eta_0(x) \left[\exp\left(-\gamma \bar\rho^\gamma(x)t - \gamma H_1(t,x) - \gamma H_1(0,x) - \gamma \int_0^t H_2(s,x)\,\mathrm{d}s\right)\right.$$

$$+ \gamma \left(\frac{\bar\rho(x)}{\partial_x \eta_0(x)}\right)^\gamma \int_0^t \exp\left(-\gamma \bar\rho^\gamma(x)(t-s) - \gamma H_1(t,x) + \gamma H_1(s,x)\right.$$

$$\left.\left. - \gamma \int_s^t H_2(\tau,x)\,\mathrm{d}\tau\right)\mathrm{d}s\right]^{1/\gamma}, \quad (3.47)$$

where we denote

$$H_1(t,x) := \int_x^1 \bar\rho(y)\partial_t \eta(t,y)\,\mathrm{d}y,$$

$$H_2(t,x) := \int_x^1 \Phi(t,x)\,\mathrm{d}x.$$

In what follows, we give estimates for $H_1$ and $H_2$.

For $H_1$, applying Hölder's inequality and **Lemma 3.1**, we have

$$|H_1(t,x)| \le \left(\int_x^1 \bar\rho(y)\,\mathrm{d}y\right)^{1/2}\left(\int_x^1 \bar\rho(y)|\partial_t \eta(t,y)|^2\,\mathrm{d}y\right)^{1/2} \le C_{s,\gamma}\epsilon_0, \quad (3.48)$$

For $H_2$, recall the (2.9) and (1.4), we have

$$H_2(t,x) = -\int_x^1 \int_{-\infty}^\infty \bar\rho(y)\frac{y-z}{|y-z|^{3-2s}}(\bar\rho(z) - \bar\rho(y))\,\mathrm{d}z\,\mathrm{d}y$$

$$+ \int_x^1 \int_{-\infty}^\infty \frac{\bar\rho(y)}{\partial_y \eta(t,y)}\frac{\eta(t,y) - \eta(t,z)}{|\eta(t,y) - \eta(t,z)|^{3-2s}}(\bar\rho(z)\partial_y \eta(t,y) - \bar\rho(y)\partial_z \tilde\eta(t,z))\,\mathrm{d}z\,\mathrm{d}y$$

$$= \int_x^1 \int_{-1}^x \bar\rho(y)\bar\rho(z) R_{1,0}(t,y,z)\,\mathrm{d}z\,\mathrm{d}y,$$

where $R_{1,0}$ has the form

$$R_{1,0}(t,x,y) = \frac{\eta(t,x) - \eta(t,y)}{|\eta(t,x) - \eta(t,y)|^{3-2s}} - \frac{x-y}{|x-y|^{3-2s}}. \quad (3.49)$$

Now let us estimate $R_{1,0}$. We compute

$$\frac{\mathrm{d}}{\mathrm{d}\alpha}\frac{\alpha(\eta(t,x) - \eta(t,y)) + (1-\alpha)(x-y)}{|\alpha(\eta(t,x) - \eta(t,y)) + (1-\alpha)(x-y)|^{3-2s}}$$

$$= 2(s-1)\frac{w(t,x) - w(t,y)}{|x-y + \alpha(w(t,x) - w(t,y))|^{3-2s}}$$

$$+ \delta_0(x - y + \alpha(w(t,x) - w(t,y)))\frac{w(t,x) - w(t,y)}{|x-y + \alpha(w(t,x) - w(t,y))|^{2-2s}}. \quad (3.50)$$

Applying the fundamental theorem of calculus to (3.49), taking (3.50) into account, we obtain

$$|R_{1,0}(t,x,y)| = 2(1-s)\frac{|w(t,x)-w(t,y)|}{|x-y|^{3-2s}} \int_0^1 \frac{1}{|1+\alpha\frac{w(t,x)-w(t,y)}{x-y}|^{3-2s}}\,\mathrm{d}\alpha,$$

where we use the **a prior assumption** (3.4) such that $|\frac{x-y}{w(t,x)-w(t,y)}| \geq A\epsilon_0 \gg 1$ if $\epsilon_0$ is sufficiently small. Moreover, by the **a prior assumption** (3.4) again, we have

$$\frac{1}{|1+\lambda\frac{w(t,x)-w(t,y)}{x-y}|^{3-2s}} \leq \frac{1}{(1-A\epsilon_0)^{3-2s}},$$

therefore

$$|R_{1,0}(t,x,y)| \leq \frac{2(1-s)A\epsilon_0}{(1-A\epsilon_0)^{3-2s}}\frac{1}{|x-y|^{2-2s}}. \tag{3.51}$$

Applying the above estimate, by noting (1.4), we have

$$|H_2(t,x)| \leq 2(1-s)A\epsilon_0 \int_x^1 \int_0^x \bar{\rho}(y)\bar{\rho}(z)\frac{1}{|y-z|^{2-2s}}\,\mathrm{d}z\,\mathrm{d}y$$
$$\leq \frac{2(1-s)A\epsilon_0}{(1-A\epsilon_0)^{3-2s}}\bar{\rho}^\gamma(x), \tag{3.52}$$

where we have used the fact that,

$$W'(y-z) = \frac{1}{|y-z|^{2-2s}}$$

for $y > z$.

Combining (3.48) and (3.52), we obtain the lower and upper bounds for $\partial_x\eta$ by (3.47),

$$L(t,x,\epsilon_0) \leq \partial_x\eta(t,x) \leq U(t,x,\epsilon_0),$$

where $U$ and $L$ respectively have the form

$$U(t,x,\epsilon) := (1+\epsilon)(1+C_{s,\gamma}\epsilon)\left[g_1(t,x,\epsilon) + \frac{1+C_{s,\gamma}\epsilon}{(1-\epsilon)^\gamma}\frac{1}{1-\frac{2(1-s)A\epsilon}{(1-A\epsilon)^{3-2s}}}\left(1-g_1(t,x,\epsilon)\right)\right]^{1/\gamma},$$

$$L(t,x,\epsilon) := (1-\epsilon)(1-C_{s,\gamma}\epsilon)\left[g_2(t,x,\epsilon) + \frac{1-C_{s,\gamma}\epsilon}{(1+\epsilon)^\gamma}\frac{1}{1+\frac{2(1-s)A\epsilon}{(1-A\epsilon)^{3-2s}}}\left(1-g_2(t,x,\epsilon)\right)\right]^{1/\gamma},$$

with

$$g_1(t,x,\epsilon) = e^{-\gamma\bar{\rho}^\gamma(x)[1-\frac{2(1-s)A\epsilon}{(1-A\epsilon)^{3-2s}}]t} \quad \text{and} \quad g_2(t,x,\epsilon) = e^{-\gamma\bar{\rho}^\gamma(x)[1+\frac{2(1-s)A\epsilon}{(1-A\epsilon)^{3-2s}}]t}.$$

In order to estimate $U$, let us compute

$$\partial_\epsilon U(t,x,\epsilon) = (1+C_{s,\gamma}\epsilon)\left[g_1(t,x,\epsilon) + \frac{1+C_{s,\gamma}\epsilon}{(1-\epsilon)^\gamma[1-\frac{2(1-s)A\epsilon}{(1-A\epsilon)^{3-2s}}]}\left(1-g_1(t,x,\epsilon)\right)\right]^{1/\gamma}$$

$$+ C_{s,\gamma}(1+\epsilon)\left[g_1(t,x,\epsilon) + \frac{1+C_{s,\gamma}\epsilon}{(1-\epsilon)^\gamma[1-\frac{2(1-s)A\epsilon}{(1-A\epsilon)^{3-2s}}]}\left(1-g_1(t,x,\epsilon)\right)\right]^{1/\gamma}$$

$$+ \frac{(1+\epsilon)(1+C_{s,\gamma}\epsilon)}{\gamma}\left[g_1(t,x,\epsilon) + \frac{1+C_{s,\gamma}\epsilon}{(1-\epsilon)^\gamma[1-\frac{2(1-s)A\epsilon}{(1-A\epsilon)^{3-2s}}]}\left(1-g_1(t,x,\epsilon)\right)\right]^{(1-\gamma)/\gamma}$$

$$\times\left[\left(2(1-s)A\frac{1+2(1-s)A\epsilon}{(1-A\epsilon)^{4-2s}} - \frac{2(1-s)A\frac{1+2(1-s)A\epsilon}{(1-A\epsilon)^{4-2s}}(1+C_{s,\gamma}\epsilon)}{(1-\epsilon)^\gamma[1-\frac{2(1-s)A\epsilon}{(1-A\epsilon)^{3-2s}}]}\right)\gamma\bar{\rho}^\gamma(x)tg_1(t,x,\epsilon)\right.$$

$$+ \left(\frac{C_{s,\gamma}}{(1-\epsilon)^\gamma[1-\frac{2(1-s)A\epsilon}{(1-A\epsilon)^{3-2s}}]} + \frac{\gamma(1+C_{s,\gamma}\epsilon)}{(1-\epsilon)^{\gamma+1}[1-\frac{2(1-s)A\epsilon}{(1-A\epsilon)^{3-2s}}]}\right.$$

$$\left.\left.+ \frac{2(1-s)A\frac{1+2(1-s)A\epsilon}{(1-A\epsilon)^{4-2s}}(1+C_{s,\gamma}\epsilon)}{(1-\epsilon)^\gamma[1-\frac{2(1-s)A\epsilon}{(1-A\epsilon)^{3-2s}}]^2}\right)\left(1-g_1(t,x,\epsilon)\right)\right], \quad (3.53)$$

and

$$\partial_{\epsilon\epsilon}U(t,x,\epsilon) = 2C_{s,\gamma}\left[g_1(t,x,\epsilon) + \frac{1+C_{s,\gamma}\epsilon}{(1-\epsilon)^\gamma[1-\frac{2(1-s)A\epsilon}{(1-A\epsilon)^{3-2s}}]}\left(1-g_1(t,x,\epsilon)\right)\right]^{1/\gamma}$$

$$+ 2(1+C_{s,\gamma}+2C_{s,\gamma}\epsilon)\frac{1}{\gamma}\left[g_1(t,x,\epsilon) + \frac{1+C_{s,\gamma}\epsilon}{(1-\epsilon)^\gamma[1-\frac{2(1-s)A\epsilon}{(1-A\epsilon)^{3-2s}}]}\left(1-g_1(t,x,\epsilon)\right)\right]^{(1-\gamma)/\gamma}$$

$$\times\left[\left(2(1-s)A\frac{1+2(1-s)A\epsilon}{(1-A\epsilon)^{4-2s}} - \frac{2(1-s)A\frac{1+2(1-s)A\epsilon}{(1-A\epsilon)^{4-2s}}(1+C_{s,\gamma}\epsilon)}{(1-\epsilon)^\gamma[1-\frac{2(1-s)A\epsilon}{(1-A\epsilon)^{3-2s}}]}\right)\gamma\bar{\rho}^\gamma(x)tg_1(t,x,\epsilon)\right.$$

$$+ \left(\frac{C_{s,\gamma}}{(1-\epsilon)^\gamma[1-\frac{2(1-s)A\epsilon}{(1-A\epsilon)^{3-2s}}]} + \frac{\gamma(1+C_{s,\gamma}\epsilon)}{(1-\epsilon)^{\gamma+1}[1-\frac{2(1-s)A\epsilon}{(1-A\epsilon)^{3-2s}}]}\right.$$

$$\left.\left.+ \frac{2(1-s)A\frac{1+2(1-s)A\epsilon}{(1-A\epsilon)^{4-2s}}(1+C_{s,\gamma}\epsilon)}{(1-\epsilon)^\gamma[1-\frac{2(1-s)A\epsilon}{(1-A\epsilon)^{3-2s}}]^2}\right)\left(1-g_1(t,x,\epsilon)\right)\right]$$

$$+ (1+\epsilon)(1+C_{s,\gamma}\epsilon)\frac{1-\gamma}{\gamma^2}\left[g_1(t,x,\epsilon) + \frac{1+C_{s,\gamma}\epsilon}{(1-\epsilon)^\gamma[1-\frac{2(1-s)A\epsilon}{(1-A\epsilon)^{3-2s}}]}\left(1-g_1(t,x,\epsilon)\right)\right]^{(1-2\gamma)/\gamma}$$

$$\times\left[\left(2(1-s)A\frac{1+2(1-s)A\epsilon}{(1-A\epsilon)^{4-2s}} - \frac{2(1-s)A\frac{1+2(1-s)A\epsilon}{(1-A\epsilon)^{4-2s}}(1+C_{s,\gamma}\epsilon)}{(1-\epsilon)^\gamma[1-\frac{2(1-s)A\epsilon}{(1-A\epsilon)^{3-2s}}]}\right)\gamma\bar{\rho}^\gamma(x)tg_1(t,x,\epsilon)\right.$$

$$+ \left(\frac{C_{s,\gamma}}{(1-\epsilon)^\gamma[1-\frac{2(1-s)A\epsilon}{(1-A\epsilon)^{3-2s}}]} + \frac{\gamma(1+C_{s,\gamma}\epsilon)}{(1-\epsilon)^{\gamma+1}[1-\frac{2(1-s)A\epsilon}{(1-A\epsilon)^{3-2s}}]}\right.$$

$$+ \frac{2(1-s)A\frac{1+2(1-s)A\epsilon}{(1-A\epsilon)^{4-2s}}(1+C_{s,\gamma}\epsilon)}{(1-\epsilon)^{\gamma}[1-\frac{2(1-s)A\epsilon}{(1-A\epsilon)^{3-2s}}]^2}\Bigg)(1-g_1(t,x,\epsilon))\Bigg]^2$$

$$+(1+\epsilon)(1+C_{s,\gamma}\epsilon)\frac{1}{\gamma}\Bigg[g_1(t,x,\epsilon)+\frac{1+C_{s,\gamma}\epsilon}{(1-\epsilon)^{\gamma}[1-\frac{2(1-s)A\epsilon}{(1-A\epsilon)^{3-2s}}]}(1-g_1(t,x,\epsilon))\Bigg]^{(1-\gamma)/\gamma}$$

$$\times\Bigg[\Bigg(\frac{4(1-s)^2A^2(1+2(1-s)A\epsilon)^2}{(1-A\epsilon)^{8-4s}}-\frac{4(1-s)^2A^2\frac{(1+2(1-s)A\epsilon)^2}{(1-A\epsilon)^{8-4s}}(1+C_{s,\gamma}\epsilon)}{(1-\epsilon)^{\gamma}[1-\frac{2(1-s)A\epsilon}{(1-A\epsilon)^{3-2s}}]}\Bigg)\gamma^2\bar{\rho}^{2\gamma}(x)t^2g_1(t,x,\epsilon)$$

$$+\Bigg(-\frac{4(1-s)A\frac{1+2(1-s)A\epsilon}{(1-A\epsilon)^{4-2s}}C_{s,\gamma}}{(1-\epsilon)^{\gamma}[1-\frac{2(1-s)A\epsilon}{(1-A\epsilon)^{3-2s}}]}-\frac{4(1-s)A\frac{1+2(1-s)A\epsilon}{(1-A\epsilon)^{4-2s}}(1+C_{s,\gamma}\epsilon)}{(1-\epsilon)^{\gamma+1}[1-\frac{2(1-s)A\epsilon}{(1-A\epsilon)^{3-2s}}]}$$

$$-\frac{8(1-s)^2A^2\frac{(1+2(1-s)A\epsilon)^2}{(1-A\epsilon)^{8-4s}}(1+C_{s,\gamma}\epsilon)}{(1-\epsilon)^{\gamma}[1-\frac{2(1-s)A\epsilon}{(1-A\epsilon)^{3-2s}}]^2}+4(1-s)(3-2s)A^2\frac{1+(1-s)A\epsilon}{(1-A\epsilon)^{5-2s}}$$

$$+\frac{4(1-s)(3-2s)A^2\frac{1+(1-s)A\epsilon}{(1-A\epsilon)^{5-2s}}(1+C_{s,\gamma}\epsilon)}{(1-\epsilon)^{\gamma}[1-\frac{2(1-s)A\epsilon}{(1-A\epsilon)^{3-2s}}]}\Bigg)\gamma\bar{\rho}^{\gamma}(x)tg_1(t,x,\epsilon)$$

$$+\Bigg(\frac{2\gamma C_{s,\gamma}}{(1-\epsilon)^{\gamma+1}[1-\frac{2(1-s)A\epsilon}{(1-A\epsilon)^{3-2s}}]}+\frac{\gamma(\gamma+1)(1+C_{s,\gamma}\epsilon)}{(1-\epsilon)^{\gamma+2}[1-\frac{2(1-s)A\epsilon}{(1-A\epsilon)^{3-2s}}]}+\frac{4(1-s)A\frac{1+2(1-s)A\epsilon}{(1-A\epsilon)^{4-2s}}C_{s,\gamma}}{(1-\epsilon)^{\gamma}[1-\frac{2(1-s)A\epsilon}{(1-A\epsilon)^{3-2s}}]^2}$$

$$+\frac{4(1-s)A\frac{1+2(1-s)A\epsilon}{(1-A\epsilon)^{4-2s}}\gamma(1+C_{s,\gamma}\epsilon)}{(1-\epsilon)^{\gamma+1}[1-\frac{2(1-s)A\epsilon}{(1-A\epsilon)^{3-2s}}]^2}+\frac{8(1-s)^2A^2\frac{(1+2(1-s)A\epsilon)^2}{(1-A\epsilon)^{8-4s}}(1+C_{s,\gamma}\epsilon)}{(1-\epsilon)^{\gamma}[1-\frac{2(1-s)A\epsilon}{(1-A\epsilon)^{3-2s}}]^3}$$

$$+\frac{4(1-s)(3-2s)A^2\frac{1+(1-s)A\epsilon}{(1-A\epsilon)^{5-2s}}(1+C_{s,\gamma}\epsilon)}{(1-\epsilon)^{\gamma}[1-\frac{2(1-s)A\epsilon}{(1-A\epsilon)^{3-2s}}]^2}\Bigg)(1-g_1(t,x,\epsilon))\Bigg]. \quad (3.54)$$

Applying Taylor's theorem to $U$ with respect to $\epsilon$ near $\epsilon=0$ for fixed $t\in[0,T]$ and $x\in[-1,1]$, we have

$$U(t,x,\epsilon_0)=1+\partial_\epsilon U(t,x,0)\epsilon_0+\frac{\partial_{\epsilon\epsilon}U(t,x,\epsilon')}{2}\epsilon_0^2,$$

for some $\epsilon'\in(0,\epsilon_0)$. Using (3.53) and (3.54), we find the following uniform estimates independent of time,

$$\partial_\epsilon U(t,x,0)=1+C_{s,\gamma}+\frac{C_{s,\gamma}+\gamma+2(1-s)A}{\gamma}(1-e^{-\gamma\bar{\rho}^{\gamma}(x)t}$$
$$\leq\frac{2(1-s)}{\gamma}A+\frac{\gamma+1}{\gamma}(1+C_{s,\gamma}),$$

and

$$\partial_{\epsilon\epsilon}U(t,x,\epsilon')\leq 2C_{s,\gamma}\left(1+\frac{1+C_{s,\gamma}\epsilon'}{(1-\epsilon')^{\gamma}[1-\frac{2(1-s)A\epsilon'}{(1-A\epsilon')^{3-2s}}]}\right)^{1/\gamma}$$

$$+ 2(1 + C_{s,\gamma} + 2C_{s,\gamma}\epsilon')\frac{1}{\gamma}\left(1 + \frac{1 + C_{s,\gamma}\epsilon'}{(1-\epsilon')^\gamma[1 - \frac{2(1-s)A\epsilon'}{(1-A\epsilon')^{3-2s}}]}\right)^{(1-\gamma)/\gamma}$$

$$\times\left(2(1-s)A\frac{1 + 2(1-s)A\epsilon'}{(1-A\epsilon')^{4-2s}e} + \frac{C_{s,\gamma}}{(1-\epsilon')^\gamma[1 - \frac{2(1-s)A\epsilon'}{(1-A\epsilon')^{3-2s}}]} + \frac{\gamma(1 + C_{s,\gamma}\epsilon')}{(1-\epsilon')^{\gamma+1}[1 - \frac{2(1-s)A\epsilon'}{(1-A\epsilon')^{3-2s}}]}\right.$$

$$+ \frac{2(1-s)A\frac{1+2(1-s)A\epsilon'}{(1-A\epsilon')^{4-2s}}(1 + C_{s,\gamma}\epsilon')}{(1-\epsilon')^\gamma[1 - \frac{2(1-s)A\epsilon'}{(1-A\epsilon')^{3-2s}}]^2} + \left.\frac{2(1-s)A\frac{1+2(1-s)A\epsilon'}{(1-A\epsilon')^{4-2s}}(1 + C_{s,\gamma}\epsilon')}{(1-\epsilon')^\gamma[1 - \frac{2(1-s)A\epsilon'}{(1-A\epsilon')^{3-2s}}]^2 e}\right)$$

$$+ (1+\epsilon')(1 + C_{s,\gamma}\epsilon')\frac{1-\gamma}{\gamma^2}\left(1 + \frac{1 + C_{s,\gamma}\epsilon'}{(1-\epsilon')^\gamma[1 - \frac{2(1-s)A\epsilon'}{(1-A\epsilon')^{3-2s}}]}\right)^{(1-2\gamma)/\gamma}$$

$$\times\left(2(1-s)A\frac{1 + 2(1-s)A\epsilon'}{(1-A\epsilon')^{3-2s}e} + \frac{C_{s,\gamma}}{(1-\epsilon')^\gamma[1 - \frac{2(1-s)A\epsilon'}{(1-A\epsilon')^{3-2s}}]} + \frac{\gamma(1 + C_{s,\gamma}\epsilon')}{(1-\epsilon')^{\gamma+1}[1 - \frac{2(1-s)A\epsilon'}{(1-A\epsilon')^{3-2s}}]}\right.$$

$$+ \left.\frac{2(1-s)A\frac{1+2(1-s)A\epsilon'}{(1-A\epsilon')^{3-2s}}(1 + C_{s,\gamma}\epsilon')}{(1-\epsilon')^\gamma[1 - \frac{2(1-s)A\epsilon'}{(1-A\epsilon')^{3-2s}}]^2} + \frac{2(1-s)A\frac{1+2(1-s)A\epsilon'}{(1-A\epsilon')^{3-2s}}(1 + C_{s,\gamma}\epsilon')}{e(1-\epsilon')^\gamma[1 - \frac{2(1-s)A\epsilon'}{(1-A\epsilon')^{3-2s}}]^2}\right)^2$$

$$+ (1+\epsilon')(1 + C_{s,\gamma}\epsilon')\frac{1}{\gamma}\left(1 + \frac{1 + C_{s,\gamma}\epsilon'}{(1-\epsilon')^\gamma[1 - \frac{2(1-s)A\epsilon'}{(1-A\epsilon')^{3-2s}}]}\right)^{(1-\gamma)/\gamma}$$

$$\times\left(4(1-s)^2A^2\frac{(1 + 2(1-s)A\epsilon)^2}{(1-A\epsilon')^{8-4s}e^2} + \frac{2\gamma C_{s,\gamma}}{(1-\epsilon')^{\gamma+1}[1 - \frac{2(1-s)A\epsilon'}{(1-A\epsilon')^{3-2s}}]} + \frac{\gamma(\gamma+1)(1 + C_{s,\gamma}\epsilon')}{(1-\epsilon')^{\gamma+2}[1 - \frac{2(1-s)A\epsilon'}{(1-A\epsilon')^{3-2s}}]}\right.$$

$$+ \frac{4(1-s)A\frac{1+2(1-s)A\epsilon}{(1-A\epsilon)^{4-2s}}C_{s,\gamma}}{(1-\epsilon')^\gamma[1 - \frac{2(1-s)A\epsilon'}{(1-A\epsilon')^{3-2s}}]^2} + \frac{8(1-s)^2A^2\frac{(1+2(1-s)A\epsilon')^2}{(1-A\epsilon')^{8-4s}}(1 + C_{s,\gamma}\epsilon')}{(1-\epsilon')^\gamma[1 - \frac{2(1-s)A\epsilon'}{(1-A\epsilon')^{3-2s}}]^3}$$

$$+ 4(1-s)(3-2s)A^2\frac{1 + (1-s)A\epsilon'}{(1-A\epsilon')^{5-2s}e} + \frac{4(1-s)(3-2s)A^2\frac{1+(1-s)A\epsilon'}{(1-A\epsilon')^{5-2s}}(1 + C_{s,\gamma}\epsilon')}{(1-\epsilon')^\gamma[1 - \frac{2(1-s)A\epsilon'}{(1-A\epsilon')^{3-2s}}]e}$$

$$+ \left.\frac{4(1-s)(3-2s)A^2\frac{1+(1-s)A\epsilon}{(1-A\epsilon')^{5-2s}}(1 + C_{s,\gamma}\epsilon')}{(1-\epsilon')^\gamma[1 - \frac{2(1-s)A\epsilon'}{(1-A\epsilon')^{3-2s}}]^2}\right)$$

$$\leq 4C_{s,\gamma} + (1 + C_{s,\gamma})\frac{1}{\gamma}2^{(1+\gamma)/\gamma}(2(1-s)A(1 + 2e^{-1}) + C_{s,\gamma} + \gamma)$$

$$+ \frac{1-\gamma}{\gamma^2}2^{(1-\gamma)/\gamma}(2(1-s)A(1 + 2e^{-1}) + C_{s,\gamma} + \gamma)^2$$

$$+ \frac{2^{1/\gamma}}{\gamma}(4(1-s^2)A^2e^{-2} + 2\gamma C_{s,\gamma} + \gamma(\gamma+1) + 4(1-s)AC_{s,\gamma} + 8(1-s)^2A^2$$

$$+ 4(1-s)(3-2s)A^2e^{-1} + 4(1-s)(3-2s)A^2e^{-1} + 4(1-s)(3-2s)A^2)$$

$$:= C'_{s,\gamma}.$$

Therefore we deduce that

$$U(t, x, \epsilon_0) \leq 1 + \left(\frac{2(1-s)}{\gamma}A + \frac{\gamma+1}{\gamma}(1 + C_{s,\gamma}) + \frac{1}{2}C'_{s,\gamma}\epsilon_0\right)\epsilon_0.$$

Since $\gamma > 2(1-s)$, let $A \geq \frac{C_{s,\gamma}}{1-10^{-N}-2(1-s)/\gamma}$ for some $N \in \mathbb{N}$ such that $(1-10^{-N})\gamma > 2(1-s)$, then

$$U(t,x,\epsilon_0) \leq 1 + (1-10^{-N-1})A\epsilon_0,$$

for $\epsilon_0$ sufficiently small. Let us notice that the upper bound is independent of $t$ and $x$, therefore we have

$$\partial_x \eta(t,x) - 1 \leq (1-10^{-N-1})A\epsilon_0, \tag{3.55}$$

for all $t \in [0,T]$ and $x \in [-1,1]$.

Similarly, we can estimate $L$ following the similar procedure as estimating $U$ to obtain

$$\partial_x \eta(t,x) - 1 \geq -(1-10^{-N-1})A\epsilon_0, \tag{3.56}$$

for all $t \in [0,T]$ and $x \in [-1,1]$.

Therefore (3.55) and (3.56) prove (3.44) in **Lemma 3.4**. $\qquad\square$

3.2. **Higher-order estimates.** In this subsection, we will give two similar results as **Lemma 3.1**, **Lemma 3.2** and **Lemma 3.3** involving higher regularity in time variable, see **Lemma 3.5**, **Lemma 3.6** and **Lemma 3.7**. An energy estimate for higher regularity in spatial variable will also be proved in **Lemma 3.8**.

First, let us give an estimate for the weighted norm of both second time derivative and mixed time-space derivative of $w$.

**Lemma 3.5.** *For $\gamma > 2(1-s)$, there exists an $\epsilon_0 > 0$ such that under the **a prior assumption** (3.4), it holds that*

$$\|\bar{\rho}^{\frac{1}{2}}\partial_{tt}w(t)\|_{L^2}^2 + \|\bar{\rho}^{\frac{\gamma}{2}}\partial_{tx}w(t)\|_{L^2}^2 + \int_0^t \|\partial_{ttx}w(s)\|_{L^2}^2\,\mathrm{d}s$$

$$\leq C_{s,\gamma}(\|\bar{\rho}^{\frac{1}{2}}\partial_t w_1\|_{L^2}^2 + \|\bar{\rho}^{\frac{\gamma}{2}}\partial_x w_1\|_{L^2}), \quad 0 \leq t \leq T, \tag{3.57}$$

*for some constant $C_{s,\gamma} > 0$.*

*Proof.* Differentiating (3.1) with respect to $t$ and multiplying the resultant by $\partial_{tt}w$, then integrating the product with respect to spatial variable, using the integration by parts and boundary conditions (3.2), we have

$$\frac{1}{2}\frac{\mathrm{d}}{\mathrm{d}t}\int_{-1}^1 \bar{\rho}|\partial_{tt}w|^2\,\mathrm{d}x + \frac{\gamma}{2}\frac{\mathrm{d}}{\mathrm{d}t}\underbrace{\int_{-1}^1 \bar{\rho}^\gamma \left(\frac{1}{1+\partial_x w}\right)^{\gamma+1}|\partial_{tx}w|^2\,\mathrm{d}x}_{(3.58)_1}$$

$$+ \frac{\gamma(\gamma+1)}{2}\underbrace{\int_{-1}^1 \bar{\rho}^\gamma \left(\frac{1}{1+\partial_x w}\right)^{\gamma+2}(\partial_{tx}w)^3\,\mathrm{d}x}_{(3.58)_2}$$

$$+ \frac{1}{2}\frac{\mathrm{d}}{\mathrm{d}t}\underbrace{\int_{-1}^1 \partial_t w \cdot \mathcal{L}\partial_t w\,\mathrm{d}x}_{(3.58)_3}$$

$$+ \frac{1-s}{2}\frac{\mathrm{d}}{\mathrm{d}t}\underbrace{\int_{-1}^1\int_{-1}^1 \bar{\rho}(x)\bar{\rho}(y)|\partial_t w(t,x) - \partial_t w(t,y)|^2 R_3(t,x,y)\,\mathrm{d}y\,\mathrm{d}x}_{(3.58)_4}$$

$$+ \frac{(1-s)(3-2s)}{2} \underbrace{\int_{-1}^1 \int_{-1}^1 \bar{\rho}(x)\bar{\rho}(y)(\partial_t w(t,x) - \partial_t w(t,y))^3 R_4(t,x,y)\,\mathrm{d}y\,\mathrm{d}x}_{(3.58)_5}$$

$$+ \frac{(1-s)(3-2s)}{2} \underbrace{\int_{-1}^1 \int_{-1}^1 -\bar{\rho}(x)\bar{\rho}(y)\frac{(\partial_t w(t,x) - w(t,y))^3}{|x-y|^{4-2s}}\,\mathrm{d}y\,\mathrm{d}x}_{(3.58)_6}$$

$$+ \underbrace{\int_{-1}^1 \frac{|\partial_{ttx} w|^2}{1+\partial_x w}\,\mathrm{d}x}_{(3.58)_7} = \underbrace{\int_{-1}^1 \left(\frac{\partial_{tx} w}{1+\partial_x w}\right)^2 \partial_{ttx} w\,\mathrm{d}x}_{(3.58)_8}, \quad (3.58)$$

where $R_3$ and $R_4$ denote

$$R_3(t,x,y) := -\frac{1}{|\eta(t,x) - \eta(t,y)|^{3-2s}} + \frac{1}{|x-y|^{3-2s}}, \qquad (3.59)$$

$$R_4(t,x,y) := -\frac{1}{|\eta(t,x) - \eta(t,y)|^{4-2s}} + \frac{1}{|x-y|^{4-2s}}. \qquad (3.60)$$

For $(3.58)_1$ and $(3.58)_7$, by the estimate $(3.4)$ in **_a priori_ assumption**, we have

$$\left(\frac{1}{1+A\epsilon_0}\right)^{\gamma+1} \int_{-1}^1 \bar{\rho}^\gamma |\partial_{tx} w|^2 \,\mathrm{d}x \leq (3.58)_1$$

$$\leq \left(\frac{1}{1-A\epsilon_0}\right)^{\gamma+1} \int_{-1}^1 \bar{\rho}^\gamma |\partial_{tx} w|^2 \,\mathrm{d}x \quad (3.61)$$

and

$$\frac{1}{1+A\epsilon_0} \int_{-1}^1 |\partial_{ttx} w|^2 \,\mathrm{d}x \leq (3.58)_7 \leq \frac{1}{1-A\epsilon_0} \int_{-1}^1 |\partial_{ttx} w|^2 \,\mathrm{d}x. \qquad (3.62)$$

For $(3.58)_2$, by the estimate $(3.4)$ in **_a priori_ assumption**, we have

$$|(3.58)_2| \leq \frac{\gamma(\gamma+1)A\epsilon_0}{2} \left(\frac{1}{1-A\epsilon_0}\right)^{\gamma+2} \int_{-1}^1 \bar{\rho}^\gamma |\partial_{tx} w|^2 \,\mathrm{d}x. \qquad (3.63)$$

For $(3.58)_3$, using a similar approach as proving the estimate $(3.14)$, we have

$$|(3.58)_3| \leq 2(1-s) \int_{-1}^1 \bar{\rho}^\gamma |\partial_{tx} w|^2 \,\mathrm{d}x. \qquad (3.64)$$

For $(3.58)_4$, we claim that

$$|R_3(t,x,y)| \leq \frac{(3-2s)A\epsilon_0(1+A\epsilon_0)}{(1-A\epsilon_0)^{5-2s}} \frac{1}{|x-y|^{3-2s}}. \qquad (3.65)$$

Then we have

$$|(3.58)_4| \leq \frac{(3-2s)A\epsilon_0(1+A\epsilon_0)}{(1-A\epsilon_0)^{5-2s}} \int_{-1}^1 \int_{-1}^1 \frac{|\partial_t w(t,x) - \partial_t w(t,y)|^2}{|x-y|^{3-2s}} \bar{\rho}(x)\bar{\rho}(y)\,\mathrm{d}y\,\mathrm{d}x.$$

Using the same estimate as dealing with $(3.58)_4$, we can bound $(3.58)_4$ as

$$|(3.58)_4| \leq 2(3-2s)A\epsilon_0 \int_{-1}^1 \bar{\rho}^\gamma |\partial_{tx} w|^2 \,\mathrm{d}x. \qquad (3.66)$$

It remains to prove (3.65). We compute

$$\frac{\mathrm{d}}{\mathrm{d}\alpha} \frac{1}{|\alpha(\eta(t,x) - \eta(t,y)) + (1-\alpha)(x-y)|^{3-2s}}$$
$$= (2s-3) \frac{x-y + \alpha(w(t,x) - w(t,y))}{|x-y+\alpha(w(t,x)-w(t,y))|^{5-2s}}(w(t,x) - w(t,y)). \quad (3.67)$$

Applying the fundamental theorem of calculus to (3.59) and taking (3.67) into account, we have

$$|R_3(t,x,y)| = (3-2s)\frac{|w(t,x) - w(t,y)|}{|x-y|^{4-2s}} \left| \int_0^1 \frac{1 + \alpha \frac{\eta(t,x)-\eta(t,y)}{x-y}}{|1 + \alpha \frac{\eta(t,x)-\eta(t,y)}{x-y}|^{5-2s}} \, \mathrm{d}\alpha \right|.$$

Let $\epsilon_0 > 0$ in the (3.4) in **a prior assumption** be sufficiently small, then we have

$$\frac{1 + \lambda \frac{\eta(t,x)-\eta(t,y)}{x-y}}{|1 + \lambda \frac{\eta(t,x)-\eta(t,y)}{x-y}|^{5-2s}} \leq \frac{1 + A\epsilon_0}{(1-A\epsilon_0)^{5-2s}},$$

therefore

$$|R_3(t,x,y)| \leq \frac{(3-2s)A\epsilon_0(1+A\epsilon_0)}{(1-A\epsilon_0)^{5-2s}}\frac{1}{|x-y|^{3-2s}}.$$

For $(3.58)_5$, we claim a pointwise estimate,

$$|R_4(t,x,y)| \leq \frac{(4-2s)A\epsilon_0(1+A\epsilon_0)}{(1-A\epsilon_0)^{6-2s}}\frac{1}{|x-y|^{4-2s}}. \quad (3.68)$$

Then we have

$$|(3.58)_5| \leq \frac{(4-2s)A\epsilon_0(1+A\epsilon_0)}{(1-A\epsilon_0)^{6-2s}} \int_{-1}^1 \int_{-1}^1 \frac{|\partial_t w(t,x) - \partial_t w(t,y)|^3}{|x-y|^{4-2s}} \bar{\rho}(x)\bar{\rho}(y) \, \mathrm{d}y \, \mathrm{d}x$$
$$\leq \frac{(4-2s)A^2\epsilon_0^2(1+A\epsilon_0)}{(1-A\epsilon_0)^{6-2s}} \int_{-1}^1 \int_{-1}^1 \frac{|\partial_t w(t,x) - \partial_t w(t,y)|^2}{|x-y|^{3-2s}} \bar{\rho}(x)\bar{\rho}(y) \, \mathrm{d}y \, \mathrm{d}x.$$

Using the same estimate as dealing with $(3.58)_3$, we can bound $(3.58)_5$ as

$$|(3.58)_5| \leq \frac{(4-2s)A^2\epsilon_0^2(1+A\epsilon_0)}{(1-A\epsilon_0)^{6-2s}} \int_{-1}^1 \int_{-1}^1 \bar{\rho}^\gamma |\partial_{tx}w|^2 \, \mathrm{d}y \, \mathrm{d}x. \quad (3.69)$$

It remains to prove (3.68). We compute

$$\frac{\mathrm{d}}{\mathrm{d}\alpha} \frac{1}{|\alpha(\eta(t,x) - \eta(t,y)) + (1-\alpha)(x-y)|^{4-2s}}$$
$$= (2s-4) \frac{x-y + \alpha(w(t,x) - w(t,y))}{|x-y+\alpha(w(t,x)-w(t,y))|^{6-2s}}(w(t,x) - w(t,y)). \quad (3.70)$$

Applying the fundamental theorem of calculus to (3.60) and taking (3.70) into account, we have

$$|R_4(t,x,y)| = (4-2s)\frac{|w(t,x) - w(t,y)|}{|x-y|^{5-2s}} \left| \int_0^1 \frac{1 + \alpha \frac{\eta(t,x)-\eta(t,y)}{x-y}}{|1 + \alpha \frac{\eta(t,x)-\eta(t,y)}{x-y}|^{6-2s}} \, \mathrm{d}\alpha \right|.$$

Let $\epsilon_0 > 0$ in the (3.4) in **a prior assumption** be sufficiently small, then we have

$$\frac{1 + \lambda \frac{\eta(t,x)-\eta(t,y)}{x-y}}{|1 + \lambda \frac{\eta(t,x)-\eta(t,y)}{x-y}|^{6-2s}} \leq \frac{1 + A\epsilon_0}{(1-A\epsilon_0)^{6-2s}}$$

and consequently

$$|R_4(t,x,y)| \leq \frac{(4-2s)A\epsilon_0(1+A\epsilon_0)}{(1-A\epsilon_0)^{6-2s}} \frac{1}{|x-y|^{4-2s}}.$$

For $(3.58)_6$, by the estimate $(3.4)$ in **a priori assumption**, we have

$$|(3.58)_6| \leq A\epsilon_0 \int_{-1}^{1}\int_{-1}^{1} \frac{|\partial_t w(t,x) - \partial_t w(t,y)|^2}{|x-y|^{3-2s}} \bar\rho(x)\bar\rho(y)\,\mathrm{d}y\,\mathrm{d}x.$$

Using the same estimate as dealing with $(3.58)_3$, we can bound $(3.58)_6$ as

$$|(3.58)_6| \leq A\epsilon_0 \int_{-1}^{1} \bar\rho^\gamma |\partial_{tx} w|^2 \,\mathrm{d}x. \tag{3.71}$$

For $(3.58)_8$, by the estimate $(3.4)$ in **a priori assumption** and Hölder's inequality, we have

$$|(3.58)_8| \leq \varepsilon \int_{-1}^{1} |\partial_{ttx} w|^2\,\mathrm{d}x + C_\varepsilon A^2 \epsilon_0^2 \int_{-1}^{1} |\partial_{tx} w|^2\,\mathrm{d}x. \tag{3.72}$$

Combining the results of $(3.61)$, $(3.63)$, $(3.64)$, $(3.66)$, $(3.69)$, $(3.71)$, $(3.62)$, $(3.72)$ and noticing that $\gamma > 2(1-s)$, by **Lemma 3.1** and choosing $\varepsilon$ and $\epsilon_0$ sufficiently small, we obtain $(3.57)$ which proves the **Lemma 3.5**. $\square$

In the following, we improve the estimate in **Lemma 3.5** for the mixed time-space derivative of $w$ near the vacuum.

**Lemma 3.6.** *For $2(1-s) < \gamma < 2$, there exists an $\epsilon_0 > 0$ such that under the **a prior assumption** $(3.4)$, it holds that*

$$\|\partial_{tx} w(t)\|_{L^2}^2 + \int_0^t \|\bar\rho^{\frac{\gamma}{2}} \partial_{tx} w(s)\|_{L^2}^2\,\mathrm{d}s$$
$$\leq C_{s,\gamma}(\|\partial_x w_1\|_{L^2}^2 + \|\bar\rho^{\frac{1}{2}}\partial_t w_1\|_{L^2}^2 + \|\bar\rho^{\frac{\gamma}{2}}\partial_x w_1\|_{L^2}^2), \quad 0 \leq t \leq T, \tag{3.73}$$

*for some constant $C_{s,\gamma} > 0$.*

*Proof.* Differentiating $(3.1)$ with respect to the time variable, then multiplying the resultant by $\partial_t w$ and integrating the product with respect to the spatial variable, using the integration by parts and boundary conditions $(3.2)$, we have

$$\frac{1}{2}\frac{\mathrm{d}}{\mathrm{d}t}\underbrace{\int_{-1}^{1}\frac{|\partial_{tx} w|^2}{1+\partial_x w}\,\mathrm{d}x}_{(3.74)_1} + \gamma\underbrace{\int_{-1}^{1}\frac{\bar\rho^\gamma |\partial_{tx} w|^2}{(1+\partial_x w)^\gamma}\,\mathrm{d}x}_{(3.74)_2}$$

$$+ \underbrace{\int_{-1}^{1}\partial_t w \cdot \mathcal{L}\partial_t w\,\mathrm{d}x}_{(3.74)_3} + \frac{1}{2}\underbrace{\int_{-1}^{1}\int_{-1}^{1}\bar\rho(x)\bar\rho(y)(\partial_t w(t,x)-\partial_t w(x,y))^2 R_3(t,x,y)\,\mathrm{d}y\,\mathrm{d}x}_{(3.74)_4}$$

$$= \frac{\mathrm{d}}{\mathrm{d}t}\underbrace{\int_{-1}^{1}-\bar\rho\partial_t w\cdot\partial_{tt} w\,\mathrm{d}x}_{(3.74)_5} + \int_{-1}^{1}\bar\rho|\partial_{tt} w|^2\,\mathrm{d}x + \frac{1}{2}\underbrace{\int_{-1}^{1}\frac{(\partial_{tx} w)^3}{|1+\partial_x w|^2}\,\mathrm{d}x}_{(3.74)_6}. \tag{3.74}$$

For $(3.74)_1$, $(3.74)_2$ and $(3.74)_6$, by **a priori assumption** $(3.4)$, we have

$$\frac{1}{1+A\epsilon_0}\int_{-1}^{1}|\partial_{tx} w|^2\,\mathrm{d}x \leq (3.74)_1 \leq \frac{1}{1-A\epsilon_0}\int_{-1}^{1}|\partial_{tx} w|^2\,\mathrm{d}x, \tag{3.75}$$

$$\frac{1}{(1+A\epsilon_0)^\gamma} \int_{-1}^{1} \bar\rho^\gamma |\partial_{tx} w|^2 \, \mathrm{d}x \leq (3.74)_2 \leq \frac{1}{(1-A\epsilon_0)^\gamma} \int_{-1}^{1} \bar\rho^\gamma |\partial_{tx} w|^2 \, \mathrm{d}x \quad (3.76)$$

and

$$|(3.74)_6| \leq \frac{A\epsilon_0}{|1+A\epsilon_0|^2} \int_{-1}^{1} |\partial_{tx} w|^2 \, \mathrm{d}x. \quad (3.77)$$

For $(3.74)_3$, using the same estimate as dealing with $(3.58)_3$ in **Lemma 3.5**, we have

$$|(3.74)_3| \leq 2(1-s) \int_{-1}^{1} \bar\rho^\gamma |\partial_{tx} w|^2 \, \mathrm{d}x. \quad (3.78)$$

For $(3.74)_4$, recall the estimate $(3.65)$, we have

$$|(3.74)_4| \leq \frac{(3-2s)A\epsilon_0(1+A\epsilon_0)}{(1-A\epsilon_0)^{5-2s}} \int_{-1}^{1} \int_{-1}^{1} \frac{|\partial_t w(t,x) - \partial_t w(t,y)|^2}{|x-y|^{3-2s}} \bar\rho(x) \bar\rho(y) \, \mathrm{d}y \, \mathrm{d}x.$$

Using the same estimate as dealing with $(3.58)_3$ in **Lemma 3.5**, we can bound $(3.74)_3$ as follows

$$|(3.74)_4| \leq \frac{(3-2s)A\epsilon_0(1+A\epsilon_0)}{(1-A\epsilon_0)^{5-2s}} \int_{-1}^{1} \bar\rho^\gamma |\partial_{tx} w|^2 \, \mathrm{d}x. \quad (3.79)$$

For $(3.74)_5$, by Hölder's inequality and Poincaré's inequality, we have

$$|(3.74)_5| \leq \varepsilon \int_{-1}^{1} |\partial_{tx} w|^2 \, \mathrm{d}x + C_\varepsilon \int_{-1}^{1} \bar\rho |\partial_{tt} w|^2 \, \mathrm{d}x, \quad (3.80)$$

where we have used the fact that $\partial_x \bar\rho \in L^\infty$ for $2(1-s) < \gamma < 2$.

Combining the results of $(3.75)$, $(3.76)$, $(3.78)$, $(3.79)$, $(3.80)$, $(3.77)$ and noticing that $\gamma > 2(1-s)$, by **Lemma 3.1** and choosing $\varepsilon$ and $\epsilon_0$ sufficiently small, we obtain $(3.73)$ which proves **Lemma 3.6**.                                      $\square$

Having obtained basic energy estimates, we further deduce a time-weighted estimate for the second-order time derivative and mixed time-space derivative of $w$.

**Lemma 3.7.** *For $2(1-s) < \gamma < 2$, there exists an $\epsilon_0 > 0$ such that under the **a prior assumption** $(3.4)$, it holds that*

$$(1+t)\|\bar\rho^{\frac{1}{2}} \partial_{tt} w(t)\|_{L^2}^2 + (1+t)\|\bar\rho^{\frac{\gamma}{2}} \partial_{tx} w(t)\|_{L^2}^2 + \int_0^t (1+s)\|\partial_{ttx} w(s)\|_{L^2}^2 \, \mathrm{d}s$$

$$\leq C_{s,\gamma} \sum_{i=0}^{1} (\|\partial_t^i \partial_x w_0\|_{L^2}^2 + \|\bar\rho^{\frac{1}{2}} \partial_t^i w_1\|_{L^2}^2 + \|\bar\rho^{\frac{\gamma}{2}} \partial_t^i \partial_x w_0\|_{L^2}^2), \quad 0 \leq t \leq T, \quad (3.81)$$

*for some constant $C_{s,\gamma} > 0$.*

*Proof.* Recall that from Lemma 3.5, we have

$$\frac{1}{2}\frac{\mathrm{d}}{\mathrm{d}t}\int_{-1}^{1}\bar{\rho}|\partial_{tt}w|^2\,\mathrm{d}x + \frac{\gamma}{2}\frac{\mathrm{d}}{\mathrm{d}t}\int_{-1}^{1}\bar{\rho}^{\gamma}\left(\frac{1}{1+\partial_x w}\right)^{\gamma+1}|\partial_{tx}w|^2\,\mathrm{d}x$$

$$+ \frac{\gamma(\gamma+1)}{2}\int_{-1}^{1}\bar{\rho}^{\gamma}\left(\frac{1}{1+\partial_x w}\right)^{\gamma+2}|\partial_{tx}w|^2\,\mathrm{d}x$$

$$+ \frac{1}{2}\frac{\mathrm{d}}{\mathrm{d}t}\int_{-1}^{1}\partial_t w \cdot \mathcal{L}\partial_t w\,\mathrm{d}x$$

$$+ \frac{1-s}{2}\frac{\mathrm{d}}{\mathrm{d}t}\int_{-1}^{1}\int_{-1}^{1}\bar{\rho}(x)\bar{\rho}(y)|\partial_t w(t,x) - \partial_t w(t,y)|^2 R_3(t,x,y)\,\mathrm{d}y\,\mathrm{d}x$$

$$+ \frac{(1-s)(3-2s)}{2}\int_{-1}^{1}\int_{-1}^{1}(\partial_t w(t,x) - \partial_t w(t,y))^3 R_4(t,x,y)\,\mathrm{d}y\,\mathrm{d}x$$

$$+ \frac{(1-s)(3-2s)}{2}\int_{-1}^{1}\int_{-1}^{1}\frac{(\partial_t w(t,x) - \partial_t w(t,y))^3}{|x-y|^{4-2s}}\,\mathrm{d}y\,\mathrm{d}x$$

$$+ \int_{-1}^{1}\frac{|\partial_{ttx}w|^2}{1+\partial_x w}\,\mathrm{d}x = \int_{-1}^{1}\left(\frac{\partial_{tx}w}{1+\partial_x w}\right)^2\partial_{ttx}w\,\mathrm{d}x.$$

Multiplying the above equation by $1+t$ we have

$$\frac{1}{2}\frac{\mathrm{d}}{\mathrm{d}t}\left[(1+t)\int_{-1}^{1}\bar{\rho}|\partial_{tt}w|^2\,\mathrm{d}x\right] - \underbrace{\frac{1}{2}\int_{-1}^{1}\bar{\rho}|\partial_{tt}w|^2\,\mathrm{d}x}_{(3.82)_1}$$

$$+ \frac{\gamma}{2}\frac{\mathrm{d}}{\mathrm{d}t}\left[\underbrace{(1+t)\int_{-1}^{1}\bar{\rho}^{\gamma}\left(\frac{1}{1+\partial_x w}\right)^{\gamma+1}|\partial_{tx}w|^2\,\mathrm{d}x}_{(3.82)_2}\right]$$

$$- \frac{\gamma}{2}\underbrace{\int_{-1}^{1}\bar{\rho}^{\gamma}\left(\frac{1}{1+\partial_x w}\right)^{\gamma+1}|\partial_{tx}w|^2\,\mathrm{d}x}_{(3.82)_3}$$

$$+ \frac{\gamma(\gamma+1)}{2}\underbrace{(1+t)\int_{-1}^{1}\bar{\rho}^{\gamma}\left(\frac{1}{1+\partial_x w}\right)^{\gamma+2}|\partial_{tx}w|^2\,\mathrm{d}x}_{(3.82)_4}$$

$$+ \frac{1}{2}\frac{\mathrm{d}}{\mathrm{d}t}\left[\underbrace{(1+t)\int_{-1}^{1}\partial_t w \cdot \mathcal{L}\partial_t w\,\mathrm{d}x}_{(3.82)_5}\right] - \frac{1}{2}\underbrace{\int_{-1}^{1}\partial_t w \cdot \mathcal{L}\partial_t w\,\mathrm{d}x}_{(3.82)_6}$$

$$+ \frac{1-s}{2}\frac{\mathrm{d}}{\mathrm{d}t}\left[\underbrace{(1+t)\int_{-1}^{1}\int_{-1}^{1}\bar{\rho}(x)\bar{\rho}(y)|\partial_t w(t,x) - \partial_t w(t,y)|^2 R_3(t,x,y)\,\mathrm{d}y\,\mathrm{d}x}_{(3.82)_7}\right]$$

$$- \frac{1-s}{2} \underbrace{\int_{-1}^{1} \int_{-1}^{1} \bar{\rho}(x)\bar{\rho}(y)|\partial_t w(t,x) - \partial_t w(t,y)|^2 R_3(t,x,y) \,\mathrm{d}y \,\mathrm{d}x}_{(3.82)_8}$$

$$+ \frac{(1-s)(3-2s)}{2} \underbrace{(1+t)\int_{-1}^{1}\int_{-1}^{1} (\partial_t w(t,x) - \partial_t w(t,y))^3 R_4(t,x,y) \,\mathrm{d}y \,\mathrm{d}x}_{(3.82)_9}$$

$$+ \frac{(1-s)(3-2s)}{2} \underbrace{(1+t)\int_{-1}^{1}\int_{-1}^{1} \frac{(\partial_t w(t,x) - \partial_t w(t,y))^3}{|x-y|^{4-2s}} \,\mathrm{d}y \,\mathrm{d}x}_{(3.82)_{10}}$$

$$+ \underbrace{(1+t)\int_{-1}^{1} \frac{|\partial_{tx} w|^2}{1+\partial_x w} \,\mathrm{d}x}_{(3.82)_{11}} = \underbrace{(1+t)\int_{-1}^{1} \left(\frac{\partial_{tx} w}{1+\partial_x w}\right)^2 \partial_{ttx} w \,\mathrm{d}x}_{(3.82)_{12}}. \quad (3.82)$$

For $(3.82)_1$, by Poincaré's inequality, we find

$$|(3.82)_1| \leq C \int_{-1}^{1} |\partial_{ttx} w|^2 \,\mathrm{d}x. \quad (3.83)$$

For $(3.82)_2$, $(3.82)_3$, $(3.82)_4$, and $(3.82)_{11}$, by the estimate $(3.4)$ in **a priori assumption**,

$$\frac{1}{(1+A\epsilon_0)^{\gamma+2}}(1+t)\int_{-1}^{1} \bar{\rho}^\gamma |\partial_{tx} w|^2 \,\mathrm{d}x \leq (3.82)_2, (3.82)_3$$

$$\leq \frac{1}{(1-A\epsilon_0)^{\gamma+2}}(1+t)\int_{-1}^{1} \bar{\rho}^\gamma |\partial_{tx} w|^2 \,\mathrm{d}x, \quad (3.84)$$

$$|(3.82)_4| \leq \frac{1}{(1-A\epsilon_0)^{\gamma+1}} \int_{-1}^{1} \bar{\rho}^\gamma |\partial_{tx} w|^2 \,\mathrm{d}x, \quad (3.85)$$

and

$$\frac{1}{1+A\epsilon_0}(1+t)\int_{-1}^{1} |\partial_{ttx} w|^2 \,\mathrm{d}x \leq (3.82)_{11} \leq \frac{1}{1-A\epsilon_0}(1+t)\int_{-1}^{1} |\partial_{ttx} w|^2 \,\mathrm{d}x. \quad (3.86)$$

For $(3.82)_5$ and $(3.82)_6$, using the same estimate as dealing with $(3.58)_3$ in **Lemma 3.5**, we have

$$|(3.82)_5| \leq (1+t)\int_{-1}^{1} \bar{\rho}^\gamma |\partial_{tx} w|^2 \,\mathrm{d}x \quad (3.87)$$

and

$$|(3.82)_6| \leq \int_{-1}^{1} \bar{\rho}^\gamma |\partial_{tx} w|^2 \,\mathrm{d}x. \quad (3.88)$$

For $(3.82)_7$, using $(3.65)$ and applying the same estimate as dealing with $(3.58)_3$ in **Lemma 3.5**, we have

$$|(3.82)_7| \leq \frac{(3-2s)A\epsilon_0(1+A\epsilon_0)}{(1-A\epsilon_0)^{5-2s}}(1+t)\int_{-1}^{1} \bar{\rho}^\gamma |\partial_{tx} w|^2 \,\mathrm{d}x. \quad (3.89)$$

For $(3.82)_8$, using $(3.68)$ applying the same estimate as dealing with $(3.58)_3$ in **Lemma 3.5**, we find that

$$|(3.82)_8| \leq \frac{(4-2s)A\epsilon_0(1+A\epsilon_0)}{(1-A\epsilon_0)^{6-2s}} \int_{-1}^{1} \bar{\rho}^\gamma |\partial_{tx}w|^2 \, \mathrm{d}x. \qquad (3.90)$$

For $(3.82)_9$, using $(3.68)$, we have

$$|(3.82)_9| \leq A^2\epsilon_0^2(1+t) \int_{-1}^{1}\int_{-1}^{1} \frac{|\partial_t w(t,x) - \partial_t w(t,y)|^2}{|x-y|^{3-2s}} \bar{\rho}(x)\bar{\rho}(y) \, \mathrm{d}y \, \mathrm{d}x. \qquad (3.91)$$

For $(3.82)_{10}$, by the estimate $(3.4)$ in **a priori** assumption, we get

$$|(3.82)_{10}| \leq A\epsilon_0(1+t) \int_{-1}^{1}\int_{-1}^{1} \frac{|\partial_t w(t,x) - \partial_t w(t,y)|^2}{|x-y|^{3-2s}} \, \mathrm{d}y \, \mathrm{d}x.$$

Using the same estimate as $(3.58)_3$ in **Lemma 3.5**, we obtain

$$|(3.82)_{10}| \leq A\epsilon_0(1+t) \int_{-1}^{1} \bar{\rho}^\gamma |\partial_{tx}w|^2 \, \mathrm{d}x. \qquad (3.92)$$

For $(3.82)_{12}$, by the estimate $(3.4)$ in **a priori** assumption and Hölder's inequality, we have

$$|(3.82)_{12}| \leq \frac{1}{2(1-A\epsilon_0)^2}(1+t) \int_{-1}^{1} |\partial_{tx}w|^2 \, \mathrm{d}x + \frac{A^2\epsilon_0^2}{2(1-A\epsilon_0)^2}(1+t) \int_{-1}^{1} |\partial_{ttx}w|^2 \, \mathrm{d}x. \qquad (3.93)$$

Combining the results of $(3.83)$, $(3.84)$, $(3.85)$, $(3.87)$, $(3.88)$, $(3.89)$, $(3.90)$, $(3.91)$, $(3.92)$, $(3.86)$, $(3.93)$, and noticing that $\gamma > 2(1-s)$, then integrating the resultant with respect to time $t$, by choosing $\epsilon_0$ sufficiently small, we obtain $(3.81)$ which proves **Lemma 3.7**. $\qquad \square$

In the following lemma, we prove a higher-order estimate in spatial variable.

**Lemma 3.8.** *For $2(1-s) < \gamma < 4/3$ there exists a constant $C > 0$ such that*

$$\|\bar{\rho}^{\frac{2\gamma-1}{2}}\partial_{xx}w(t)\|_{L^2}^2 + \int_0^t \|\bar{\rho}^{\frac{3\gamma-1}{2}}\partial_{xx}w(s)\|_{L^2}^2 \, \mathrm{d}s$$

$$\leq C_{s,\gamma}(\|\bar{\rho}^{\frac{1}{2}}w_1\|_{L^2}^2 + \|\partial_x w_0\|_{L^2}^2 + \|\bar{\rho}^{\frac{\gamma}{2}}\partial_x w_1\|_{L^2}^2 + \|\bar{\rho}^{\frac{2\gamma-2}{2}}\partial_{xx}w_0\|_{L^2}^2), \quad 0 \leq t \leq T, \qquad (3.94)$$

*for some constant $C_{s,\gamma} > 0$.*

*Proof.* Let us set

$$\mathcal{G} := \ln \partial_x \eta, \quad \tilde{\mathcal{G}} := \ln \partial_x \tilde{\eta}. \qquad (3.95)$$

Then, we have that for arbitrary $\alpha > 0$, under the **a priori** assumption $(3.4)$,

$$\frac{1}{2}\bar{\rho}^\alpha \partial_x w \leq \bar{\rho}^\alpha \mathcal{G} \leq \bar{\rho}^\alpha \partial_x w,$$

where we used the following elementary algebraic inequality

$$\frac{1}{2}z \leq \ln(1+z) \leq z$$

for $|z| \ll \epsilon$ with $\epsilon$ sufficiently small. Moreover, we have

$$\frac{1}{1+A\epsilon_0}\bar{\rho}^\alpha \partial_{xx}w \leq \bar{\rho}^\alpha \partial_x \mathcal{G} \leq \frac{1}{1-A\epsilon_0}\bar{\rho}^\alpha \partial_{xx}w.$$

Now let us consider the equation for $\mathcal{G}$. From the equation (2.8) we have

$$\partial_{tx}\mathcal{G} + \gamma\left(\frac{\bar{\rho}}{\partial_x\eta}\right)^{\gamma}\partial_x\mathcal{G} = \bar{\rho}\partial_{tt}\eta + \bar{\rho}(W'*\bar{\rho})\left[\left(\frac{1}{\partial_x\eta}\right)^{\gamma}-1\right] + \Phi. \qquad (3.96)$$

Multiplying (3.96) by $\bar{\rho}^{2\gamma-1}\partial_x\mathcal{G}$, then integrating the product with respect to spatial variable, we have

$$\frac{1}{2}\frac{\mathrm{d}}{\mathrm{d}t}\int_{-1}^{1}\bar{\rho}^{2\gamma-1}|\partial_x\mathcal{G}|^2\,\mathrm{d}x + \underbrace{\int_{-1}^{1}\gamma\frac{\bar{\rho}^{3\gamma-1}}{(\partial_x\eta)^{\gamma}}|\partial_x\mathcal{G}|^2\,\mathrm{d}x}_{(3.97)_1}$$

$$= \underbrace{\int_{-1}^{1}\bar{\rho}^{2\gamma}\partial_x\mathcal{G}\cdot\partial_{tt}\eta\,\mathrm{d}x}_{(3.97)_2} + \underbrace{\int_{-1}^{1}\bar{\rho}^{2\gamma}(W'*\bar{\rho})\partial_x\mathcal{G}\cdot\left[\left(\frac{1}{\partial_x\eta}\right)^{\gamma}-1\right]\,\mathrm{d}x}_{(3.97)_3} + \underbrace{\int_{-1}^{1}\bar{\rho}^{2\gamma-1}\partial_x\mathcal{G}\cdot\Phi\,\mathrm{d}x}_{(3.97)_4}. \qquad (3.97)$$

In what follows, we estimate $(3.97)_1 - (3.97)_4$ term by term where $\varepsilon > 0$ is chosen to be sufficiently small.

For $(3.97)_1$, using (3.4), we find that

$$\frac{\gamma}{(1+A\epsilon)^{\gamma}}\int_{-1}^{1}\bar{\rho}^{3\gamma-1}|\partial_x\mathcal{G}|^2\,\mathrm{d}x \leq (3.97)_1 \leq \frac{\gamma}{(1-A\epsilon)^{\gamma}}\int_{-1}^{1}\bar{\rho}^{3\gamma-1}|\partial_x\mathcal{G}|^2\,\mathrm{d}x. \qquad (3.98)$$

For $(3.97)_2$, by Young's inequality, we have

$$(3.97)_2 \leq \varepsilon\int_{-1}^{1}\bar{\rho}^{3\gamma-1}|\partial_x\mathcal{G}|^2\,\mathrm{d}x + C_{\varepsilon}\int_{-1}^{1}\bar{\rho}^{\gamma+1}|\partial_{tt}w|^2\,\mathrm{d}x. \qquad (3.99)$$

For $(3.97)_3$ it follows from (3.4) and Young's inequality that

$$(3.97)_3 \leq \varepsilon\int_{-1}^{1}\bar{\rho}^{3\gamma-1}|\partial_x\mathcal{G}|^2\,\mathrm{d}x + C_{\varepsilon}\int_{-1}^{1}\bar{\rho}^{\gamma+1}|\partial_x w|^2\,\mathrm{d}x, \qquad (3.100)$$

where we used the following elementary algebraic inequality

$$\left|1-\left(\frac{1}{1+z}\right)^{\gamma}\right| \leq 2\gamma|z|.$$

For $(3.97)_4$, we introduce the expression of $\Phi$ in $(2.9)$ and by splitting it into four terms, we deduce

$$
\begin{aligned}
(3.97)_4 &= \underbrace{\int_{-1}^{1} \bar{\rho}^{2\gamma-1} \partial_x \mathcal{G} \cdot \mathcal{L}w \, dx}_{(3.101)_1} \\
&\quad + \underbrace{\int_{-1}^{1} \bar{\rho}^{2\gamma}(x) \partial_x \mathcal{G}(x) \int_{-1}^{1} (\bar{\rho}(y) - \bar{\rho}(x)) R_{1,1}(x,y) \, dy \, dx}_{(3.101)_2} \\
&\quad + \underbrace{\int_{-1}^{1} \bar{\rho}^{2\gamma+1}(x) \partial_x \mathcal{G}(x) \int_{-1}^{1} (\partial_x w(x) - \partial_y w(y)) R_{1,0}(x,y) \frac{1}{1+\partial_x w} \, dy \, dx}_{(3.101)_3} \\
&\quad + \underbrace{\int_{-1}^{1} \bar{\rho}^{2\gamma+1}(x) \partial_x \mathcal{G}(x) \int_{-1}^{1} (\partial_x w(x) - \partial_y w(y)) \frac{x-y}{|x-y|^{3-2s}} \frac{-\partial_x w(x)}{1+\partial_x w(x)} \, dy \, dx}_{(3.101)_4}.
\end{aligned}
\tag{3.101}
$$

For $(3.101)_1$, we need to estimate two terms using the expression of $\mathcal{L}$ in $(3.7)$,

$$
\begin{aligned}
(3.101)_1 &= \underbrace{\int_{-\infty}^{\infty} \bar{\rho}^{2\gamma+1}(x) \partial_x \tilde{\mathcal{G}}(x) \int_{-\infty}^{\infty} \frac{x-y}{|x-y|^{3-2s}} (\partial_x \tilde{w}(x) - \partial_y \tilde{w}(y)) \, dy \, dx}_{(3.102)_1} \\
&\quad - \underbrace{\int_{-\infty}^{\infty} \bar{\rho}^{2\gamma}(x) \partial_x \tilde{\mathcal{G}}(x) \int_{-\infty}^{\infty} \frac{-2(1-s)}{|x-y|^{3-2s}} (\bar{\rho}(x) - \bar{\rho}(y))(\tilde{w}(x) - \tilde{w}(y)) \, dy \, dx}_{(3.102)_2}.
\end{aligned}
\tag{3.102}
$$

For $(3.102)_1$, we further split it into three parts,

$$
\begin{aligned}
(3.102)_1 &= \underbrace{\int_{-\infty}^{\infty} \bar{\rho}^{\frac{\gamma+3}{2}}(x) \partial_x \tilde{\mathcal{G}}(x) \int_{-\infty}^{\infty} \frac{x-y}{|x-y|^{3-2s}} (\bar{\rho}^{\frac{3\gamma-1}{2}}(x) \partial_x \tilde{w}(x) - \bar{\rho}^{\frac{3\gamma-1}{2}}(y) \partial_y \tilde{w}(y)) \, dy \, dx}_{(3.103)_1} \\
&\quad - \underbrace{\int_{-\infty}^{\infty} \bar{\rho}^{\frac{\gamma+1}{2}}(x) \partial_x \tilde{\mathcal{G}}(x) \int_{-\infty}^{\infty} \frac{x-y}{|x-y|^{3-2s}} (\bar{\rho}^{\frac{3\gamma-1}{2}}(x) - \bar{\rho}^{\frac{3\gamma-1}{2}}(y))(\bar{\rho}(x) - \bar{\rho}(y)) \partial_y \tilde{w}(y) \, dy \, dx}_{(3.103)_2} \\
&\quad - \underbrace{\int_{-\infty}^{\infty} \bar{\rho}^{\frac{\gamma+1}{2}}(x) \partial_x \tilde{\mathcal{G}}(x) \int_{-\infty}^{\infty} \frac{x-y}{|x-y|^{3-2s}} (\bar{\rho}^{\frac{3\gamma-1}{2}}(x) - \bar{\rho}^{\frac{3\gamma-1}{2}}(y)) \bar{\rho}(y) \partial_y \tilde{w}(y) \, dy \, dx}_{(3.103)_3}.
\end{aligned}
\tag{3.103}
$$

For $(3.103)_1$, by Hölder's inequality, Lemma B.3, and the definition of $H^{1-s}$, we obtain that

$$|(3.103)_1| \leq \left( \int_{-1}^{1} \int_{-1}^{1} \frac{1}{|x-y|^{1-2s}} \bar{\rho}^{\gamma+3}(x) |\partial_{xx} w|^2 \, \mathrm{d}x \, \mathrm{d}y \right)^{\frac{1}{2}}$$

$$\times \left( \int_{-1}^{1} \int_{-1}^{1} \frac{|\bar{\rho}^{\frac{3\gamma-1}{2}}(x) \partial_x w(x) - \bar{\rho}^{\frac{3\gamma-1}{2}}(y) \partial_y w(y)|^2}{|x-y|^{3-2s}} \, \mathrm{d}x \, \mathrm{d}y \right)^{\frac{1}{2}}$$

$$\leq \|\bar{\rho}^{\frac{\gamma+3}{2}} \partial_x \mathcal{G}\|_{L^2} \|\bar{\rho}^{\frac{3\gamma-1}{2}} \partial_x w\|_{H^{1-s}}.$$

Then by Gagliardo-Nirenberg's inequality we have

$$\|\bar{\rho}^{\frac{3\gamma-1}{2}} \partial_x w\|_{H^{1-s}} \leq \|\bar{\rho}^{\frac{3\gamma-1}{2}} \partial_x w\|_{L^2}^s \|\bar{\rho}^{\frac{3\gamma-1}{2}} \partial_x w\|_{H^1}^{1-s}$$

$$\leq C \|\bar{\rho}^{\frac{3\gamma-1}{2}} \partial_x w\|_{L^2} + C \|\bar{\rho}^{\frac{3\gamma-1}{2}} \partial_x w\|_{L^2}^s (\|\bar{\rho}^{\frac{3\gamma-3}{2}} \partial_x w\|_{L^2} + \|\bar{\rho}^{\frac{3\gamma-1}{2}} \partial_x \mathcal{G}\|_{L^2})^{1-s},$$

where we have used the fact that $\bar{\rho} \in W^{1,\infty}$. By Hardy-type inequality **Corollary B.1**, we find

$$\|\bar{\rho}^{\frac{3\gamma-3}{2}} \partial_x w\|_{L^2} \leq \|\bar{\rho}^{\frac{3\gamma-3}{2}+1} \partial_x w\|_{L^2} + \|\bar{\rho}^{\frac{3\gamma-3}{2}+1} \partial_x \mathcal{G}\|_{L^2}.$$

Therefore, by Young's inequality, we infer

$$|(3.103)_1| \leq \varepsilon \|\bar{\rho}^{\frac{3\gamma-1}{2}} \partial_x \mathcal{G}\|_{L^2}^2 + C_\epsilon \|\bar{\rho}^{\frac{\gamma}{2}} \partial_x w\|_{L^2}^2. \tag{3.104}$$

For $(3.103)_2$, applying integration by parts, we have

$$(3.103)_2 =$$

$$\underbrace{\int_{-\infty}^{\infty} \partial_x (\bar{\rho}^{\frac{\gamma+1}{2}}(x)) \cdot \tilde{\mathcal{G}}(x) \int_{-\infty}^{\infty} \frac{x-y}{|x-y|^{3-2s}} (\bar{\rho}^{\frac{3\gamma-1}{2}}(x) - \bar{\rho}^{\frac{3\gamma-1}{2}}(y))(\bar{\rho}(x) - \bar{\rho}(y)) \partial_y \tilde{w}(y) \, \mathrm{d}y \, \mathrm{d}x}_{(3.105)_1}$$

$$+ \underbrace{\int_{-\infty}^{\infty} \bar{\rho}^{\frac{\gamma+1}{2}}(x) \tilde{\mathcal{G}}(x) \int_{-\infty}^{\infty} \frac{-2(1-s)}{|x-y|^{3-2s}} (\bar{\rho}^{\frac{3\gamma-1}{2}}(x) - \bar{\rho}^{\frac{3\gamma-1}{2}}(y))(\bar{\rho}(x) - \bar{\rho}(y)) \partial_y \tilde{w}(y) \, \mathrm{d}y \, \mathrm{d}x}_{(3.105)_2}$$

$$+ \underbrace{\int_{-\infty}^{\infty} \bar{\rho}^{\frac{\gamma+1}{2}}(x) \tilde{\mathcal{G}}(x) \int_{-\infty}^{\infty} \frac{x-y}{|x-y|^{3-2s}} \partial_x (\bar{\rho}^{\frac{3\gamma-1}{2}}(x)) \cdot (\bar{\rho}(x) - \bar{\rho}(y)) \partial_y \tilde{w}(y) \, \mathrm{d}y \, \mathrm{d}x}_{(3.105)_3}$$

$$+ \underbrace{\int_{-\infty}^{\infty} \bar{\rho}^{\frac{\gamma+1}{2}}(x) \tilde{\mathcal{G}}(x) \int_{-\infty}^{\infty} \frac{x-y}{|x-y|^{3-2s}} (\bar{\rho}^{\frac{3\gamma-1}{2}}(x) - \bar{\rho}^{\frac{3\gamma-1}{2}}(y)) \partial_x (\bar{\rho}(x)) \cdot \partial_y \tilde{w}(y) \, \mathrm{d}y \, \mathrm{d}x}_{(3.105)_4}.$$

$$\tag{3.105}$$

For $(3.105)_1$, we further split it into two terms,

$$(3.105)_1 =$$

$$\underbrace{\int_{-\infty}^{\infty} \frac{\gamma+1}{2} \bar{\rho}^{\frac{\gamma+1}{2}}(x) \partial_x \bar{\rho}(x) \cdot \mathcal{G}(x) \int_{-\infty}^{\infty} \frac{x-y}{|x-y|^{3-2s}} (\bar{\rho}^{\frac{3\gamma-1}{2}}(x) - \bar{\rho}^{\frac{3\gamma-1}{2}}(y)) \partial_y \tilde{w}(y) \, \mathrm{d}y \, \mathrm{d}x}_{(3.106)_1}$$

$$- \underbrace{\int_{-\infty}^{\infty} \frac{\gamma+1}{2} \bar{\rho}^{\frac{\gamma-1}{2}}(x) \partial_x \bar{\rho}(x) \cdot \mathcal{G}(x) \int_{-\infty}^{\infty} \frac{x-y}{|x-y|^{3-2s}} (\bar{\rho}^{\frac{3\gamma-1}{2}}(x) - \bar{\rho}^{\frac{3\gamma-1}{2}}(y)) \bar{\rho}(y) \partial_y \tilde{w}(y) \, \mathrm{d}y \, \mathrm{d}x}_{(3.106)_2}.$$

$$\tag{3.106}$$

For $(3.106)_1$, by Hölder's inequality and Lemma B.3, then choosing a fixed constant $\delta \in (-2s, 4 - 3\gamma]$, we have

$$|(3.106)_1| \leq \left( \int_{-1}^{1} \int_{-1}^{1} \frac{\bar{\rho}^{-1+\delta}(y)}{|x-y|^{1-2s}} \bar{\rho}^{\gamma+1}(x) |\mathcal{G}(x)|^2 \, \mathrm{d}y \, \mathrm{d}x \right)^{\frac{1}{2}}$$
$$\times \left( \int_{-1}^{1} \int_{-1}^{1} \frac{1}{|x-y|^{1-2s}} \bar{\rho}^{1-\delta}(y) |\partial_y \tilde{w}(y)|^2 \, \mathrm{d}y \, \mathrm{d}x \right)^{\frac{1}{2}}$$
$$\leq \|\bar{\rho}^{\frac{\gamma+2s+\delta}{2}} \mathcal{G}\|_{L^2} \|\bar{\rho}^{\frac{1-\delta}{2}} \partial_x w\|_{L^2}.$$

By the choice of $\delta$ such that $\delta > -2s$, we find

$$\|\bar{\rho}^{\frac{\gamma+2s+\delta}{2}} \mathcal{G}\|_{L^2} \leq C \|\bar{\rho}^{\frac{\gamma}{2}} \partial_x w\|_{L^2}.$$

In addition, by **Corollary B.1** and the choice of $\delta$ such that $\delta \leq 4 - 3\gamma$, we have

$$\|\bar{\rho}^{\frac{1-\delta}{2}} \partial_y w\|_{L^2} \leq C(\|\bar{\rho}^{\frac{3-\delta}{2}} \partial_x w\|_{L^2} + \|\bar{\rho}^{\frac{3-\delta}{2}} \partial_{xx} w\|_{L^2}) \leq C(\|\bar{\rho}^{\frac{\gamma}{2}} \partial_x w\|_{L^2} + \|\bar{\rho}^{\frac{3\gamma-1}{2}} \partial_x \mathcal{G}\|_{L^2}).$$

Therefore Young's inequality gives

$$|(3.106)_1| \leq \varepsilon \|\bar{\rho}^{\frac{3\gamma-1}{2}} \partial_x \mathcal{G}\|_{L^2}^2 + C_\varepsilon \|\bar{\rho}^{\frac{\gamma}{2}} \partial_x w\|_{L^2}^2,$$

for $\varepsilon > 0$ being sufficiently small.

Similarly, we estimate $(3.106)_2$ and $(3.105)_2 - (3.105)_4$, by Hölder's inequality and choosing $\delta \in (-2s, 4 - 3\gamma]$, we obtain

$$|(3.106)_2| \leq \left( \int_{-1}^{1} \int_{-1}^{1} \frac{1}{|x-y|^{1-2s}} \bar{\rho}^{\gamma-\delta}(x) |\mathcal{G}(x)|^2 \, \mathrm{d}y \, \mathrm{d}x \right)^{\frac{1}{2}}$$
$$\times \left( \int_{-1}^{1} \int_{-1}^{1} \frac{\bar{\rho}^{-1+\delta}(x)}{|x-y|^{1-2s}} \bar{\rho}^2(y) |\partial_y w(y)|^2 \, \mathrm{d}y \, \mathrm{d}x \right)^{\frac{1}{2}},$$

$$|(3.105)_2| \leq \left( \int_{-1}^{1} \int_{-1}^{1} \frac{\bar{\rho}^{-1+\delta}(y)}{|x-y|^{1-2s}} \bar{\rho}^{\gamma+1}(x) |\mathcal{G}(x)|^2 \, \mathrm{d}y \, \mathrm{d}x \right)^{\frac{1}{2}}$$
$$\times \left( \int_{-1}^{1} \frac{1}{|x-y|^{1-2s}} \bar{\rho}^{1-\delta}(y) |\partial_y w(y)|^2 \, \mathrm{d}x \right)^{\frac{1}{2}},$$

$$|(3.105)_3| \leq \left( \int_{-1}^{1} \int_{-1}^{1} \frac{\bar{\rho}^{-1+\delta}(y)}{|x-y|^{1-2s}} \bar{\rho}^{2\gamma-1}(x) |\mathcal{G}(x)|^2 \, \mathrm{d}y \, \mathrm{d}x \right)^{\frac{1}{2}}$$
$$\times \left( \int_{-1}^{1} \int_{-1}^{1} \frac{1}{|x-y|^{1-2s}} \bar{\rho}^{1-\delta}(y) |\partial_y w(y)|^2 \, \mathrm{d}y \, \mathrm{d}x \right)^{\frac{1}{2}},$$

and

$$|(3.105)_4| \leq \left( \int_{-1}^{1} \int_{-1}^{1} \frac{\bar{\rho}^{-1+\delta}(y)}{|x-y|^{1-2s}} \bar{\rho}^{\frac{\gamma+1}{2}}(x) |\mathcal{G}(x)|^2 \, \mathrm{d}y \, \mathrm{d}x \right)^{\frac{1}{2}}$$
$$\times \left( \int_{-1}^{1} \int_{-1}^{1} \frac{1}{|x-y|^{1-2s}} \bar{\rho}^{1-\delta}(y) |\partial_y w(y)|^2 \, \mathrm{d}y \, \mathrm{d}x \right)^{\frac{1}{2}},$$

Then following the same approach as estimating $(3.106)_1$ by using Lemma B.3, **Corollary B.1**, and Young's inequality, we have

$$|(3.106)_2|, |(3.105)_2|, |(3.105)_3|, |(3.105)_4| \le \varepsilon \|\bar{\rho}^{\frac{3\gamma-1}{2}} \partial_x \mathcal{G}\|_{L^2}^2 + C_\varepsilon \|\bar{\rho}^{\frac{\gamma}{2}} \partial_x w\|_{L^2}.$$

Therefore the following holds

$$|(3.103)_2| \le \varepsilon \|\bar{\rho}^{\frac{3\gamma-1}{2}} \partial_x \mathcal{G}\|_{L^2}^2 + C_\varepsilon \|\bar{\rho}^{\frac{\gamma}{2}} \partial_x w\|_{L^2}^2. \tag{3.107}$$

For $(3.103)_3$, by Hölder's inequality and Young's inequality, we obtain

$$|(3.103)_3| \le \left(\int_{-1}^1 \int_{-1}^1 \frac{1}{|x-y|^{1-2s}} \bar{\rho}^{\gamma+2-\delta}(x) |\partial_x \mathcal{G}(x)|^2 \, \mathrm{d}y \, \mathrm{d}x\right)^{\frac12}$$
$$\times \left(\int_{-1}^1 \int_{-1}^1 \frac{\bar{\rho}^{-1+\delta}(x)}{|x-y|^{1-2s}} \bar{\rho}^2(y) |\partial_x \tilde{w}(x)|^2 \, \mathrm{d}y \, \mathrm{d}x\right)^{\frac12},$$

where $\delta \in (-2s, 4-3\gamma)$. Then following the same approach as estimating $(3.106)_1$ by using Lemma B.3, **Corollary B.1**, and Young's inequality, we find

$$|(3.103)_3| \le \varepsilon \|\bar{\rho}^{\frac{3\gamma-1}{2}} \partial_x \mathcal{G}\|_{L^2}^2 + C_\varepsilon \|\bar{\rho}^{\frac{\gamma}{2}} \partial_x w\|_{L^2}^2. \tag{3.108}$$

for $\varepsilon > 0$ being sufficiently small.

Combining the results of (3.104), (3.107), and (3.108), we have

$$|(3.102)_1| \le \varepsilon \|\bar{\rho}^{\frac{3\gamma-1}{2}} \partial_x \mathcal{G}\|_{L^2}^2 + C_\varepsilon \|\bar{\rho}^{\frac{\gamma}{2}} \partial_x w\|_{L^2}^2. \tag{3.109}$$

For $(3.102)_2$, we further split it into two terms,

$$(3.102)_2 \le \underbrace{\int_{-\infty}^\infty \bar{\rho}^{2\gamma-\frac12}(x) |\partial_x \tilde{\mathcal{G}}|(x) \int_{-\infty}^\infty \frac{2(1-s)}{|x-y|^{3-2s}} |\bar{\rho}(x) - \bar{\rho}(y)|^{\frac32} |\tilde{w}(x) - \tilde{w}(y)| \, \mathrm{d}y \, \mathrm{d}x}_{(3.110)_1}$$
$$+ \underbrace{\int_{-\infty}^\infty \bar{\rho}^{2\gamma-\frac12}(x) |\partial_x \tilde{\mathcal{G}}|(x) \int_{-\infty}^\infty \frac{2(1-s)}{|x-y|^{3-2s}} |\bar{\rho}(x) - \bar{\rho}(y)| \bar{\rho}^{\frac12}(y) |\tilde{w}(x) - \tilde{w}(y)| \, \mathrm{d}y \, \mathrm{d}x}_{(3.110)_2}. \tag{3.110}$$

For $(3.110)_1$, by Hölder's inequality and (3.14), we deduce that

$$|(3.110)_1| \le \left(\int_{-\infty}^\infty \int_{-\infty}^\infty \bar{\rho}^{-\frac{3-4s}{3-2s}}(y) \bar{\rho}^{4\gamma-1-\frac{3-4s}{3-2s}}(x) |\partial_x \tilde{\mathcal{G}}(x)|^2 \, \mathrm{d}y \, \mathrm{d}x\right)^{\frac12}$$
$$\times \left(\int_{-\infty}^\infty \int_{-\infty}^\infty \frac{|\tilde{w}(x) - \tilde{w}(y)|^2}{|x-y|^{3-4s}} \bar{\rho}^{\frac{3-4s}{3-2s}}(x) \bar{\rho}^{\frac{3-4s}{3-2s}}(y) \, \mathrm{d}y \, \mathrm{d}x\right)^{\frac12}$$
$$\le C \left(\int_{-1}^1 \bar{\rho}^{4\gamma-1-\frac{3-4s}{3-2s}}(x) |\partial_x \mathcal{G}(x)|^2 dx\right)^{\frac12}$$
$$\times \left(\int_{-1}^1 \int_{-1}^1 \frac{|w(x) - w(y)|^2}{|x-y|^{3-2s}} \bar{\rho}(x) \bar{\rho}(y) \, \mathrm{d}y \, \mathrm{d}x\right)^{\frac{3-4s}{2(3-2s)}} \|w\|_{L^\infty}^{\frac{2}{2(3-2s)}}.$$

By Newton-Leibnitz's rule and Hölder's inequality, then choosing a fixed constant $\delta' \in (0, 4-3\gamma)$, we have

$$\|w\|_{L^\infty} \le \|\partial_x w\|_{L^1} \le \|\bar{\rho}^{-\frac{1-\delta'}{2}}\|_{L^2} \|\bar{\rho}^{\frac{1-\delta'}{2}} \partial_x w\|_{L^2} = C \|\bar{\rho}^{\frac{1-\delta'}{2}} \partial_x w\|_{L^2}.$$

Then by **Corollary B.1**, we have

$$\|w\|_{L^\infty} \le C(\|\bar{\rho}^{\frac{3-\delta'}{2}} \partial_x w\|_{L^2} + \|\bar{\rho}^{\frac{3-\delta'}{2}} \partial_x \mathcal{G}\|_{L^2}).$$

Therefore, by (3.14), Young's inequality and the choice of $\delta'$ such that $\delta' < 4 - 3\gamma$, we have

$$(3.110)_1 \le \varepsilon \|\bar{\rho}^{\frac{3\gamma-1}{2}} \partial_x \mathcal{G}\|_{L^2}^2 + C_\varepsilon \|\bar{\rho}^{\frac{\gamma}{2}} \partial_x w\|_{L^2}^2.$$

For $(3.110)_2$, by Hölder's inequality, (3.14), and Young's inequality, we obtain

$$(3.110)_2 \le \left( \int_{-1}^1 \int_{-1}^1 \frac{1}{|x-y|^{1-2s}} \bar{\rho}^{4\gamma-2}(x) |\partial_x \tilde{\mathcal{G}}(x)|^2 \, dy \, dx \right)^{\frac{1}{2}}$$

$$\times \left( \int_{-1}^1 \int_{-1}^1 \frac{|\tilde{w}(x) - \tilde{w}(y)|^2}{|x-y|^{3-2s}} \bar{\rho}(x)\bar{\rho}(y) \, dy \, dx \right)^{\frac{1}{2}}$$

$$\le \varepsilon \|\bar{\rho}^{\frac{3\gamma-1}{2}} \partial_x \mathcal{G}\|_{L^2}^2 + C_\varepsilon \|\bar{\rho}^{\frac{\gamma}{2}} \partial_x w\|.$$

Then we can infer that

$$|(3.102)_2| \le \varepsilon \|\bar{\rho}^{\frac{3\gamma-1}{2}} \partial_x \mathcal{G}\|_{L^2}^2 + C_\varepsilon \|\bar{\rho}^{\frac{\gamma}{2}} \partial_x w\|_{L^2}^2. \tag{3.111}$$

Combining the results of (3.109) and (3.111), we find

$$|(3.101)_1| \le \varepsilon \|\bar{\rho}^{\frac{3\gamma-1}{2}} \partial_x \mathcal{G}\|_{L^2}^2 + C_\varepsilon \|\bar{\rho}^{\frac{\gamma}{2}} \partial_x w\|_{L^2}^2. \tag{3.112}$$

For $(3.101)_2$, by (3.30), we have

$$|(3.101)_2| \le 2(1-s)(3-2s)A^2\epsilon_0^2$$

$$\times \int_{-1}^1 \bar{\rho}^{2\gamma}(x) |\partial_x \mathcal{G}(x)| \int_{-1}^1 |\bar{\rho}(y) - \bar{\rho}(x)| \frac{|w(x) - w(y)|}{|x-y|^{2-2s}} \, dy \, dx.$$

Then by a similar approach as the estimate for part $(3.102)_2$, we can bound $(3.101)_2$ as

$$|(3.101)|_2 \le A\epsilon_0(\|\bar{\rho}^{\frac{3\gamma-1}{2}} \partial_x \mathcal{G}\|_{L^2}^2 + \|\bar{\rho}^{\frac{\gamma}{2}} \partial_x w\|_{L^2}^2). \tag{3.113}$$

For $(3.101)_3$ and $(3.101)_4$, by (3.51), we can show that

$$|(3.101)_3|, |(3.101)_4|$$

$$\le \frac{2(1-s)A\epsilon_0}{(1-A\epsilon_0)^{3-2s}} \underbrace{\int_{-1}^1 \bar{\rho}^{\frac{\gamma+3}{2}}(x) |\partial_x \mathcal{G}(x)| \int_{-1}^1 \frac{|\bar{\rho}^{\frac{3\gamma-1}{2}}(y) \partial_y w(y) - \bar{\rho}^{\frac{3\gamma-1}{2}}(x) \partial_x w(x)|}{|x-y|^{2-2s}} \, dy \, dx}_{(3.114)_1}$$

$$+ \frac{2(1-s)A\epsilon_0}{(1-A\epsilon_0)^{3-2s}} \underbrace{\int_{-1}^1 \bar{\rho}^{\frac{\gamma+3}{2}}(x) |\partial_x \mathcal{G}(x)| \int_{-1}^1 \frac{|\bar{\rho}^{\frac{3\gamma-1}{2}}(x) - \bar{\rho}^{\frac{3\gamma-1}{2}}(y)| \cdot |\partial_y w(y)|}{|x-y|^{2-2s}} \, dy \, dx}_{(3.114)_2}. \tag{3.114}$$

For $(3.114)_1$ and $(3.114)_2$, following the same approach as estimating $(3.106)_1$ by using Hölder's inequality, **Lemma B.3**, **Corollary B.1**, and Young's inequality, we have

$$|(3.114)_1|, |(3.114)_2| \le \varepsilon \|\bar{\rho}^{\frac{3\gamma-1}{2}} \partial_x \mathcal{G}\|_{L^2}^2 + C_\varepsilon \|\bar{\rho}^{\frac{\gamma}{2}} \partial_x w\|_{L^2}^2.$$

Therefore, we obtain

$$|(3.101)_3|, |(3.101)_4| \le A\epsilon_0(\|\bar{\rho}^{\frac{\gamma+3}{2}} \partial_x \mathcal{G}\|_{L^2}^2 + \|\bar{\rho}^{\frac{\gamma}{2}} \partial_x w\|_{L^2}^2). \tag{3.115}$$

Combining the results of (3.112), (3.113), and (3.115), we find

$$|(3.97)_4| \leq \varepsilon \|\bar{\rho}^{\frac{3\gamma-1}{2}} \partial_x \mathcal{G}\|_{L^2}^2 + C_\varepsilon \|\bar{\rho}^{\frac{\gamma}{2}} \partial_x w\|_{L^2}^2. \tag{3.116}$$

Collecting the results of (3.98), (3.99), (3.100), and (3.116), by choosing $\varepsilon$ and $\epsilon_0$ sufficiently small, we obtain (3.94) which proves **Lemma 3.8**.

<div style="text-align: right">□</div>

## 4. Proof of main theorems

### 4.1. **Proof of Theorem 2.1.**
The proof consists of two parts. First, we prove the local existence and uniqueness of strong solutions on some finite time interval. For this purpose, recall that, from (2.8), we have

$$\partial_t \eta = v, \tag{4.1}$$

$$\bar{\rho} \partial_t v + \partial_x \left[ \left( \frac{\bar{\rho}}{\partial_x \eta} \right)^\gamma - \bar{\rho}^\gamma \right] + \Phi = \partial_x \left( \frac{\partial_x v}{\partial_x \eta} \right), \tag{4.2}$$

where $\Phi$ is given by

$$\Phi(t,x) := \partial_x(\bar{\rho}^\gamma(x))$$
$$+ \frac{\bar{\rho}(x)}{\partial_x \eta(t,x)} \int_{-\infty}^\infty W'(\eta(t,x) - \tilde{\eta}(t,y))(\bar{\rho}(y)\partial_x \eta(t,x) - \bar{\rho}(x)\partial_y \tilde{\eta}(t,y)) \, \mathrm{d}y.$$

The finite difference scheme is defined as follows. Let $N$ be a positive integer and $h := 1/N$. We decompose the domain $[-1,1]$ by using a mesh $x_k := kh$, $k = -N, \cdots, 0, \cdots, N$. Then we approximate (4.1) and (4.2) by the following ordinary differential equations,

$$\frac{\mathrm{d}\eta_k}{\mathrm{d}t} = v_k, \tag{4.3}$$

$$\bar{\rho}_k \frac{\mathrm{d}v_k}{\mathrm{d}t} + \frac{1}{h} \left\{ \bar{\theta}_{k+1} \left[ \left( \frac{h}{\eta_{k+1} - \eta_k} \right)^\gamma - 1 \right] - \bar{\theta}_k \left[ \left( \frac{h}{\eta_k - \eta_{k-1}} \right)^\gamma - 1 \right] \right\} + \Phi_k$$
$$= \frac{1}{h} \left( \frac{v_{k+1} - v_k}{\eta_{k+1} - \eta_k} - \frac{v_k - v_{k-1}}{\eta_k - \eta_{k-1}} \right), \tag{4.4}$$

where

$$\bar{\rho}_k := \bar{\rho}(x_k),$$
$$\bar{\theta}_k := \sum_{i \geq k+1} \sum_{j \neq i} W'(x_i - x_j)\bar{\rho}_i(\bar{\rho}_j - \bar{\rho}_i)h^2,$$
$$\Phi_k := - \bar{\rho}_k \sum_{j \neq k} hW'(x_k - x_j)(\bar{\rho}_j - \bar{\rho}_k)$$
$$+ \frac{h\bar{\rho}_k}{\eta_{k+1} - \eta_k} \sum_{j \neq k} hW'(\eta_k - \tilde{\eta}_j) \left( \bar{\rho}_j \frac{\eta_{k+1} - \eta_k}{h} - \bar{\rho}_k \frac{\tilde{\eta}_{j+1} - \tilde{\eta}_j}{h} \right),$$
$$\eta_{-k} := - \eta_k,$$
$$\tilde{\eta}_k := \begin{cases} \eta_{-N} + 1 + kh, & k \leq -N-1, \\ \eta_k, & -N \leq k \leq N, \\ \eta_N - 1 + kh, & k \geq N+1. \end{cases}$$

Notice that the formulas for the approximation of $\Phi$ and $\bar{\rho}^\gamma(x)$ in (2.9) and (4.4) are chosen such that they are consistent with the steady state equation in (1.4). Moreover, this will allow us to obtain a discrete version of the energy estimate. From (2.2) and (2.13), this system is supplemented by the following initial conditions,

$$\eta_k(0) = \eta_0(x_k), \quad v_k(0) = \eta_1(x_k), \tag{4.5}$$

and to match the boundary conditions (2.12) and (2.14), we further impose

$$\eta_0(t) = 0, \quad \frac{\eta_N(t) - \eta_{N-1}(t)}{h} = 1. \tag{4.6}$$

Then, we use (4.6) to determine $\eta_N$ and $v_N$ in terms of $\eta_{N-1}$ and $v_{N-1}$. It follows that

$$\eta_N = \eta_{N-1} + h, \quad v_N = v_{N-1}. \tag{4.7}$$

In addition, we define the approximation of $E$ in (2.20),

$$E_N(t) := \max_{-N \le k \le N} \left| \left( \frac{\eta_k(t) - \eta_{k-1}(t)}{h} - 1, \frac{v_k(t) - v_{k-1}(t)}{h} \right) \right|^2$$
$$+ \sum_{k=-N+1}^{N-1} h \bar{\rho}_k^{2\gamma-1} \left| \frac{1}{h} \left( \frac{\eta_{k+1}(t) - \eta_k(t)}{h} - \frac{\eta_k(t) - \eta_{k-1}(t)}{h} \right) \right|^2 + \sum_{k=-N}^{N} h \bar{\rho}_k \left| \frac{dv_k(t)}{dt} \right|^2.$$

Then following the argument in [24], we can prove the existence of $\eta_1, \cdots, \eta_{n-1}$ and obtain the following lemma.

**Lemma 4.1.** *Let* $0 < s < 1/2$ *and* $\gamma > 2(1-s)$. *Suppose* $E(0) < \infty$ *and* $\partial_x \eta_0(\pm 1) = 1$. *Then there exist* $N_0 \in \mathbb{N}$ *and* $\epsilon_0 > 0$ *such that if* $E(0) < \epsilon_0$ *and* $N > N_0$, *then the problem admits a unique solution* $(\eta_1, \cdots, \eta_{N-1})$ *on* $[0, T^*]$ *for some constant* $T^* > 0$ *independent of* $N$ *satisfying*

$$E_N(t) \le C E_N(0) \le 2 C E(0), \quad 0 \le t \le T^*,$$

*for some constant* $C > 0$ *independent of* $N$.

*Proof.* It follows from the ODE theory that the ODE system (4.3) and (4.4) with the initial data (4.5) has a solution on a time interval. Let $T_N > 0$ be the maximum existence time. Then for $k = -N, \cdots, 0, \cdots, N$, and $t \in [0, T_N)$,

$$\left| \frac{\eta_k(t) - \eta_{k-1}(t)}{h} - 1 \right| \le \left| \frac{\eta_0(x_k) - \eta_0(x_{k-1})}{h} - 1 \right| + \int_0^t \left| \frac{v_k(s) - v_{k-1}(s)}{h} \right| \, \mathrm{d}s$$
$$\le \|\eta_{0x} - 1\|_{L^\infty} + t \sup_{s \in [0,t]} \sqrt{E_N(s)}.$$

This implies, due to $\eta_0(t) = 0$, that for $k = -N, \cdots, 0, \cdots, N$ and $t \in [0, T_N)$,

$$\left| \frac{\eta_k(t)}{x_k} - 1 \right| \le \max_{-N \le k \le N} \left| \frac{\eta_k(t) - \eta_{k-1}(t)}{h} - 1 \right| \le \|\eta_{0x} - 1\|_{L^\infty} + t \sup_{s \in [0,t]} \sqrt{E_N(s)}.$$

Therefore, for $k = -N, \cdots, 0, \cdots, N$ and $t \in [0, T]$,

$$\left| \frac{\eta_k(t) - \eta_{k-1}(t)}{h} - 1 \right| \le 2\epsilon, \quad \left| \frac{\eta_k(t)}{x_k} - 1 \right| \le 2\epsilon, \tag{4.8}$$

provided that

$$\|\eta_{0x} - 1\|_{L^\infty} \le \epsilon, \quad T \sup_{s \in [0,T]} \sqrt{E_N(s)} \le \epsilon, \tag{4.9}$$

for some constant $\epsilon > 0$.

Similar as the proof in [24], we prove that for sufficiently large $N \in \mathbb{N}$, there exists a constant $\bar{C}$ independent of $N$ such that if $T$ satisfies (4.9), then

$$E_N(t) \leq \bar{C} E_N(0) + \bar{C} \int_0^t E_N(s) + E_N^2(s)\, \mathrm{d}s, \qquad (4.10)$$

for $0 \leq t \leq T$.

First, taking the time derivative of both sides of (4.4), we obtain

$$\bar{\rho}_k \frac{\mathrm{d}^2 v_k}{\mathrm{d}t^2} + \frac{\gamma}{h}\left[ -\bar{\theta}_{k+1}\left(\frac{h}{\eta_{k+1} - \eta_k}\right)^{\gamma+1} \frac{v_{k+1} - v_k}{h} \right.$$
$$\left. + \bar{\theta}_k\left(\frac{h}{\eta_k - \eta_{k-1}}\right)^{\gamma+1} \frac{v_k - v_{k-1}}{h} \right] + \frac{\mathrm{d}\Phi_k}{\mathrm{d}t}$$
$$= \frac{1}{h}\left[ -\left(\frac{v_{k+1} - v_k}{\eta_{k+1} - \eta_k}\right)^2 + \frac{\frac{\mathrm{d}v_{k+1}}{\mathrm{d}t} - \frac{\mathrm{d}v_k}{\mathrm{d}t}}{\eta_{k+1} - \eta_k} - \left(\frac{v_k - v_{k-1}}{\eta_k - \eta_{k-1}}\right)^2 + \frac{\frac{\mathrm{d}v_k}{\mathrm{d}t} - \frac{\mathrm{d}v_{k-1}}{\mathrm{d}t}}{\eta_k - \eta_{k-1}} \right]. \quad (4.11)$$

Multiplying (4.11) by $h\frac{\mathrm{d}v_k}{\mathrm{d}t}$ and taking the summation of the product from $k = -N + 1$ to $k = N - 1$, since most of the terms are similar to the case in [24], we only discuss the third term on the left-hand side of (4.11).

For the third term, we have

$$\sum_{k=-N+1}^{N-1} h\frac{\mathrm{d}v_k}{\mathrm{d}t}\frac{\mathrm{d}\Phi_k}{\mathrm{d}t} = (s-1)\frac{\mathrm{d}}{\mathrm{d}t}\underbrace{\sum_{-N \leq j < k \leq N} h^2 \bar{\rho}_j \bar{\rho}_k \frac{1}{|x_j - x_k|^{1-2s}}\left(\frac{v_j - v_k}{x_j - x_k}\right)^2}_{(4.12)_1}$$

$$+ \frac{1-s}{2}\frac{\mathrm{d}}{\mathrm{d}t}\underbrace{\sum_{k=-N}^{N}\sum_{j \neq k} h^2 \bar{\rho}_j \bar{\rho}_k \left(-\frac{1}{|\eta_j - \eta_k|^{3-2s}} + \frac{1}{|x_j - x_k|^{3-2s}}\right)(v_j - v_k)^2}_{(4.12)_2}$$

$$+ \frac{(1-s)(3-2s)}{2}\underbrace{\sum_{k=-N}^{N}\sum_{j \neq k} h^2 \bar{\rho}_j \bar{\rho}_k \left(-\frac{\eta_j - \eta_k}{|\eta_j - \eta_k|^{5-2s}} + \frac{x_j - x_k}{|x_j - x_k|^{5-2s}}\right)(v_j - v_k)^3}_{(4.12)_3}.$$
$$(4.12)$$

In addition, we derive from (4.8) that

$$(1 - 2\epsilon)(x_k - x_{k-1}) \leq \eta_k - \eta_{k-1} \leq (1 + 2\epsilon)(x_k - x_{k-1}),$$

which further implies that

$$(1 - 2\epsilon)(x_k - x_j) \leq \eta_k - \eta_j \leq (1 + 2\epsilon)(x_k - x_j), \quad k > j.$$

Therefore, for $k \neq j$, we have

$$\left| -\frac{1}{|\eta_j - \eta_k|^{3-2s}} + \frac{1}{|x_j - x_k|^{3-2s}} \right| \leq 2(3 - 2s)\epsilon \frac{1}{|x_j - x_k|^{3-2s}},$$

where we have used the following elementary algebraic inequalities,

$$\left| 1 - \frac{1}{|1 + z|^{3-2s}} \right| \leq 2(3 - 2s)\epsilon.$$

Then, we estimate $(4.12)_2$ as follows,

$$|(4.12)_2| \leq 4(3-2s)\epsilon \sum_{-N \leq j < k \leq N} h^2 \bar{\rho}_j \bar{\rho}_k \frac{1}{|x_j - x_k|^{1-2s}} \left( \frac{v_j - v_k}{x_j - x_k} \right)^2. \quad (4.13)$$

For $(4.12)_3$, similarly, we obtain

$$\left| -\frac{1}{|\eta_j - \eta_k|^{4-2s}} + \frac{1}{|x_j - x_k|^{4-2s}} \right| \leq 2(4-2s)\epsilon \frac{1}{|x_j - x_k|^{4-2s}},$$

where we have used the algebraic inequality,

$$\left| 1 - \frac{1}{|1+z|^{4-2s}} \right| \leq 2(4-2s)\epsilon.$$

Moreover, we deduce from the definition of $E_N$ that,

$$-E_N(t)(x_k - x_{k-1}) \leq v_k - v_{k-1} \leq E_N(t)(x_k - x_{k-1}),$$

which further implies that

$$-E_N(t)(x_k - x_j) \leq v_k - v_j \leq E_N(t)(x_k - x_j).$$

Therefore, we estimate $(4.12)_3$ as follows,

$$|(4.12)_3| \leq 4(4-2s)\epsilon E_N^{\frac{1}{2}}(t) \sum_{-N \leq j < k \leq N} h^2 \bar{\rho}_j \bar{\rho}_k \frac{1}{|x_j - x_k|^{1-2s}} \left( \frac{v_j - v_k}{x_j - x_k} \right)^2. \quad (4.14)$$

It remains to consider $(4.12)_1$. By the definition of $E_N(t)$, we have

$$|v_k(t) - v_{k-1}(t)|^2 \leq E_N(t),$$

which implies that

$$|v_j(t) - v_k(t)|^2 \leq (k-j)E_N(t).$$

Therefore, we can estimate $(4.12)_1$ as

$$\sum_{-N \leq j < k \leq N} h^2 \bar{\rho}_j \bar{\rho}_k \frac{1}{|x_j - x_k|^{1-2s}} \left( \frac{v_j - v_k}{x_j - x_k} \right)^2$$
$$\leq \sum_{-N \leq j < k \leq N} h \bar{\rho}_i \bar{\rho}_j \frac{E_N(t)}{|x_j - x_k|^{1-2s}} \leq C E_N(t). \quad (4.15)$$

Indeed, by Cauchy's inequality,

$$|v_i - v_j| = \left| \sum_{k=j}^{i-1} \frac{v_{k+1} - v_k}{h} h \right| \leq \left[ \sum_{k=j}^{i-1} \left( \frac{v_{k+1} - v_k}{h} \right)^2 h \right]^{\frac{1}{2}} (x_i - x_j)^{\frac{1}{2}} \leq E_N^{\frac{1}{2}}(t)(x_i - x_j),$$

which implies $(4.15)$. Therefore, by applying $(4.15)$ to $(4.13)$ and $(4.14)$, we can obtain

$$\sum_{k=-N}^{N} h \bar{\rho}_k \left| \frac{\mathrm{d}v_k}{\mathrm{d}t} \right|^2 \leq \bar{C} E_N(0) + \bar{C} \int_0^t E_N(s) + E_N^2(s) \, \mathrm{d}s. \quad (4.16)$$

Second, we denote $\mathcal{G}_k = \ln \frac{\eta_{k+1} - \eta_k}{h}$. Now, we rewrite (4.4) as follows,

$$\frac{\frac{\mathrm{d}\mathcal{G}_{k+1}}{\mathrm{d}t} - \frac{\mathrm{d}\mathcal{G}_k}{\mathrm{d}t}}{h} - \bar{\theta}_{k+1} \frac{1}{h} \left[ \left( \frac{h}{\eta_{k+1} - \eta_k} \right)^\gamma - \left( \frac{h}{\eta_k - \eta_{k-1}} \right)^\gamma \right]$$
$$= \bar{\rho}_k \frac{\mathrm{d}v_k}{\mathrm{d}t} + \frac{\bar{\theta}_{k+1} - \bar{\theta}_k}{h} \left[ \left( \frac{h}{\eta_{k+1} - \eta_k} \right)^\gamma - 1 \right] + \Phi_k. \quad (4.17)$$

Multiplying (4.17) by $\bar{\rho}_k^{2\gamma-1} \frac{1}{h} (\mathcal{G}_{k+1} - \mathcal{G}_k)$ and taking the summation of the product from $k = -N + 1$ to $k = N - 1$, since most of the terms are similar to the case in [24], we only discuss the third term on the right-hand side.

For the third on the right-hand side, by direct computation, we have

$$\sum_{k=-N}^{N} h \bar{\rho}_k^{2\gamma-1} \frac{\mathcal{G}_{k+1} - \mathcal{G}_k}{h} \Phi_k$$

$$= \underbrace{\sum_{k=-N}^{N} h^2 \bar{\rho}_k^{2\gamma+1} \frac{\mathcal{G}_{k+1} - \mathcal{G}_k}{h} \sum_{j \neq k} \frac{x_k - x_j}{|x_k - x_j|^{3-2s}} \left( \frac{\eta_{k+1} - \eta_k}{h} - \frac{\eta_{j+1} - \eta_j}{h} \right)}_{(4.18)_1}$$

$$+ \underbrace{\sum_{k=-N}^{N} h^2 \bar{\rho}_k^{2\gamma} \frac{\mathcal{G}_{k+1} - \mathcal{G}_k}{h} \sum_{j \neq k} \frac{x_k - x_j}{|x_k - x_j|^{3-2s}} (\bar{\rho}_k - \bar{\rho}_j)(\eta_k - \eta_j)}_{(4.18)_2}$$

$$+ \sum_{k=-N}^{N} h^2 \bar{\rho}_k^{2\gamma+1} \frac{\mathcal{G}_{k+1} - \mathcal{G}_k}{h}$$
$$\underbrace{\times \sum_{j \neq k} \left( \frac{\eta_k - \eta_j}{|\eta_k - \eta_j|^{3-2s}} - \frac{x_k - x_j}{|x_k - x_j|^{3-2s}} \right) \left( \frac{\eta_{k+1} - \eta_k}{h} - \frac{\eta_{j+1} - \eta_j}{h} \right)}_{(4.18)_3}$$

$$+ \sum_{k=-N}^{N} h^2 \bar{\rho}_k^{2\gamma+1} \frac{\mathcal{G}_{k+1} - \mathcal{G}_k}{h}$$
$$\underbrace{\times \sum_{j \neq k} \frac{1}{|\eta_k - \eta_j|^{3-2s}} \left( \frac{h}{\eta_{k+1} - \eta_k} - 1 \right) \left( \frac{\eta_{k+1} - \eta_k}{h} - \frac{\eta_{j+1} - \eta_j}{h} \right)}_{(4.18)_4}$$

$$+ \sum_{k=-N}^{N} h^2 \bar{\rho}_k^{2\gamma} \frac{\mathcal{G}_{k+1} - \mathcal{G}_k}{h}$$
$$\underbrace{\times \sum_{j \neq k} \left( \frac{\eta_k - \eta_j}{|\eta_k - \eta_j|^{3-2s}} - \frac{x_k - x_j}{|x_k - x_j|^{3-2s}} - \frac{2(s-1)}{|x_k - x_j|^{3-2s}} (\eta_k - \eta_j) \right) (\bar{\rho}_j - \bar{\rho}_k)}_{(4.18)_5}.$$

$$(4.18)$$

Then, for $(4.18)_1$, we derive from $(4.8)$ that,

$$(1 - 2\epsilon)(x_k - x_{k-1}) \leq \eta_k - \eta_{k-1} \leq (1 + 2\epsilon)(x_k - x_{k-1}),$$

which further implies that

$$(1 - 2\epsilon)(x_k - x_j) \leq \eta_k - \eta_j, \eta_{k+1} - \eta_{j+1} \leq (1 + 2\epsilon)(x_k - x_j), \quad k > j.$$

Then we can obtain

$$|(4.18)_1| \leq \sum_{k=-N}^{N} h^2 \bar{\rho}_k^{2\gamma+1} \left| \frac{\mathcal{G}_{k+1} - \mathcal{G}_k}{h} \right| \sum_{j \neq k} \frac{4\epsilon}{|x_k - x_j|^{1-2s}} \leq \bar{C}(1 + E_N(t)), \quad (4.19)$$

where we have used the fact that

$$\left| \frac{\mathcal{G}_{k+1} - \mathcal{G}_k}{h} \right| \leq \frac{1}{1 - 2\epsilon} \left| \frac{1}{h} \left( \frac{\eta_{k+1} - \eta_k}{h} - \frac{\eta_k - \eta_{k-1}}{h} \right) \right|,$$

since $(4.8)$ holds and there exists $\theta \in (0, 1)$ such that

$$\ln x - \ln y = \frac{1}{\theta x + (1 - \theta)y}(x - y),$$

by the Lagrange's mean value theorem. Similarly, we can show that

$$|(4.18)_2, (4.18)_3, (4.18)_4, (4.18)_5|$$
$$\leq \sum_{k=-N}^{N} h^2 \bar{\rho}_k^{2\gamma} \frac{\mathcal{G}_{k+1} - \mathcal{G}_k}{h} \sum_{j \neq k} |x_k - x_j|^{2s} \leq \bar{C}(1 + E_N(t)). \quad (4.20)$$

Combining $(4.19)$ and $(4.20)$, we have

$$\sum_{k=-N}^{N} h \bar{\rho}_k^{2\gamma-1} \frac{\mathcal{G}_{k+1} - \mathcal{G}_k}{h} \Phi_k \leq \bar{C} E_N(t). \quad (4.21)$$

Therefore, we can obtain

$$\sum_{k=-N}^{N} h \bar{\rho}_k^{2\gamma-1} \left| \frac{1}{h} \left( \frac{\eta_{k+1} - \eta_k}{h} - \frac{\eta_{k+1} - \eta_k}{h} \right) \right|^2 \leq \bar{C} E_N(0) + \bar{C} \int_0^t E_N(s) + E_N^2(s) \, \mathrm{d}s. \quad (4.22)$$

Then $(4.16)$ and $(4.22)$ imply $(4.10)$. Set $T^* := \sup\{t : E_N(t) \leq 2\bar{C} E_N(0)\}$, where $\bar{C} > 1$ is the constant in $(4.10)$. It suffices to consider the case $T^* < \epsilon/\sqrt{2\bar{C} E_N(0)}$; otherwise the result is already obtained since $E_N(0) \to E(0)$ as $N \to \infty$. By the definition of $T^*$, we have

$$T^* \sup \sqrt{E_N(s)} \leq T^* \sqrt{2\bar{C} E_N(0)} \leq \epsilon,$$

implying that $T^*$ satisfies $(4.9)$. Therefore, it follows from $(4.10)$ that

$$2\bar{C} E_N(0) = E_N(T^*) \leq \bar{C} E_N(0) + \bar{C} \int_0^{T^*} E_N(s) + E_N^2(s) \, \mathrm{d}s$$
$$\leq \bar{C} E_N(0) + 2\bar{C}^2 E_N(0)(1 + 2\bar{C} E_N(0)) T^*,$$

which implies

$$T^* \geq \frac{1}{2\bar{C}(1 + 2\bar{C} E_N(0))} \geq \frac{1}{2\bar{C}(1 + 4\bar{C} E(0))}.$$

The lemma is proved by taking $C = 2\bar{C}$.     □

To obtain the approximation sequence of $\eta$, we define $\eta_h$ by

$$\eta_h(t,x) := \eta_{k-1}(t) + \frac{\eta_k(t) - \eta_{k-1}(t)}{h}(x - x_{k-1}),$$

for $x_{k-1} < x < x_k$. Then by using Lemma 4.1, we can obtain all the needed estimates as in [24, Lemma A.5], then following the exact compactness argument in [24, Proposition A.6], we can prove the convergence of $\{\eta_h\}_h$ as $h$ goes to zero and obtain the unique strong solution in the sense of **Definition 2.1**.

**Lemma 4.2.** *Let $0 < s < 1/2$ and $\gamma > 2(1-s)$. Suppose $E(0) < \infty$ and $\partial_x \eta_0(\pm 1) = 1$. Then there exists an $\epsilon_0 > 0$ such that if $E(0) < \epsilon_0$, then the problem admits a unique solution $\eta$ on $[0, T^*]$ for some constant $T^* > 0$ satisfying*

$$E(t) \leq C E(0), \quad 0 \leq t \leq T^*, \tag{4.23}$$

*for some constant $C > 0$.*

Second, to obtain the global well-posedness of strong solutions from local existence and uniqueness theorem through the standard continuation argument, we need to derive the uniform-in-time boundedness of the following energy functional in (2.20),

$$E(t) = \|(\partial_x \eta - 1, \partial_{tx} \eta)(t)\|_{L^\infty}^2 + \|\bar{\rho}^{\frac{2\gamma-1}{2}} \partial_{xx} \eta(t)\|_{L^2}^2 + \|\bar{\rho}^{\frac{1}{2}} \partial_{tt} \eta(t)\|_{L^2}^2.$$

The uniform bounds of $\|\partial_x \eta(t) - 1\|_{L^\infty}$ can be proved by using **Lemma 3.4** and the standard continuity argument. The estimates of $\|\bar{\rho}^{\frac{2\gamma-1}{2}} \partial_{xx} \eta(t)\|_{L^2}$ and $\|\bar{\rho}^{\frac{1}{2}} \partial_{tt} \eta(t)\|_{L^2}$ are given by **Lemma 3.5** and **Lemma 3.8**. It only remains to prove the uniform boundedness of $\|\partial_{tx} \eta\|_{L^\infty}$. By Gagliardo-Nirenberg inequality, we have

$$\begin{aligned}
\|\partial_{tx} \eta(t)\|_{L^\infty} &\leq C \|\partial_x \eta(t)\|_{L^\infty} \left\| \frac{\partial_{tx} \eta(t)}{\partial_x \eta(t)} \right\|_{L^2}^{\frac{1}{2}} \left\| \frac{\partial_{tx} \eta(t)}{\partial_x \eta(t)} \right\|_{H^1}^{\frac{1}{2}} \\
&\leq C \|\partial_x \eta(t)\|_{L^\infty} \left( \left\| \frac{\partial_{tx} \eta(t)}{\partial_x \eta(t)} \right\|_{L^2} + \left\| \left( \frac{\partial_{tx} \eta(t)}{\partial_x \eta(t)} \right)_x \right\|_{L^2} \right).
\end{aligned}$$

Recall the definition (3.95), we find

$$\begin{aligned}
\|\partial_{tx} \eta(t)\|_{L^\infty} &\leq C \|\partial_x \eta(t)\|_{L^\infty} \left( \left\| \frac{\partial_{tx} \eta(t)}{\partial_x \eta(t)} \right\|_{L^2} + \|\partial_{tx} \mathcal{G}(t)\|_{L^2} \right) \\
&\leq C \|\partial_x \eta(t)\|_{L^\infty} \left( \left\| \frac{1}{\partial_x \eta(t)} \right\|_{L^\infty} \|\partial_{tx} \eta(t)\|_{L^2} + \|\partial_{tx} \mathcal{G}(t)\|_{L^2} \right).
\end{aligned}$$

Since the uniform bound of $\|\partial_x \eta(t) - 1\|_{L^\infty}$ indicates that $\|\partial_x \eta(t)\|_{L^\infty}$ and $\|1/\partial_x \eta(t)\|_{L^\infty}$ are uniformly bounded, while the uniform bound of $\|\partial_{tx} \eta(t)\|_{L^2}$ is given by **Lemma 3.6**. Let us consider $\|\partial_{tx} \mathcal{G}(t)\|_{L^2}$. Recall (3.96), we can obtain

$$\begin{aligned}
\|\partial_{tx} \mathcal{G}(t)\|_{L^2} &\leq C_\gamma (1 + \|\bar{\rho}^\gamma \partial_x \mathcal{G}(t)\|_{L^2} + \|\bar{\rho} \partial_{tt} \eta(t)\|_{L^2} + \|\Phi(t)\|_{L^2}) \\
&\leq C_\gamma (1 + \|\bar{\rho}^{\frac{2\gamma-1}{2}} \partial_{xx} \eta(t)\|_{L^2} + \|\bar{\rho}^{\frac{1}{2}} \partial_{tt} \eta(t)\|_{L^2} + \|\Phi(t)\|_{L^2}),
\end{aligned}$$

for some constant $C_\gamma > 0$. As $\|\bar{\rho}^{\frac{2\gamma-1}{2}} \partial_{xx} \eta(t)\|_{L^2}$ and $\|\bar{\rho}^{\frac{1}{2}} \partial_{tt} \eta(t)\|_{L^2}$ are uniformly bounded, it suffices to estimate $\|\Phi(t)\|_{L^2}$. Since

$$\|\Phi(t)\|_{L^2} \leq \|\bar{\rho}\|_{L^\infty}^{\frac{1}{2}} \|\bar{\rho}^{-\frac{1}{2}} \Phi(t)\|_{L^2},$$

and

$$\|\bar{\rho}^{-\frac{1}{2}} \Phi(t)\|_{L^2} \leq C(\|\partial_x \eta(t)\|_{L^\infty}^2 + \|\bar{\rho}^{\frac{2\gamma-1}{2}} \partial_{xx} \eta(t)\|_{L^2}^2),$$

we have that $\|\Phi(t)\|_{L^2}$ is uniformly bounded by **Lemma B.6**. Therefore $\|\partial_{tx}\mathcal{G}(t)\|_{L^2}$ is uniformly bounded. Therefore we have the uniform boundedness of $\|\partial_{tx}\eta(t)\|_{L^\infty}$. This finishes the proof of (2.17) in **Theorem 2.1**.

The time-decay estimates (2.18) and (2.19) in **Theorem 2.1** are direct results of **Lemma 3.3** and **Lemma 3.7**.

4.2. **Proof of Theorem 2.2.** The result follows by the construction (2.16) of $(\rho, u)$.

## Appendix A. Stationary states

In this subsection, let us recall some important theorems concerning the existence and uniqueness modulo translations corresponding to (1.4). This appendix is a summary of the main results in [8, 4, 9, 2]. Consider the equation

$$\nabla(\bar{\rho}^\gamma) + \bar{\rho}\nabla\Psi_k[\bar{\rho}] = 0, \quad \mathbf{x} \in \mathbb{R}^n. \tag{A.1}$$

**Remark A.1.** *We recall the critical diffusion exponent $\gamma_c := 1 - k/n$. Then $\gamma = \gamma_c$ is called the fair-competition regime; $\gamma > \gamma_c$ is called the diffusion-dominated regime; $0 < \gamma < \gamma_c$ is called the attraction-dominated regime.*

Let us give the definition of solutions to (A.1) according to [8, 4, 9, 2].

**Definition A.1.** *Given $\bar{\rho} \in L^1_+(\mathbb{R}^n) \cap L^\infty(\mathbb{R}^n)$ with $\|\bar{\rho}\|_{L^1} = 1$. We call $\bar{\rho}$ a solution to (A.1), if $\bar{\rho}^\gamma \in W^{1,2}_{loc}(\mathbb{R}^n)$, $\nabla\Psi_k[\bar{\rho}] \in L^1_{loc}(\mathbb{R}^n)$ and it satisfies (A.1) in the sense of distributions in $\mathbb{R}^n$. If $-n < k \leq 1 - n$, we further require $\bar{\rho} \in C^{0,\alpha}(\mathbb{R}^n)$ for some $\alpha \in (1 - k - n, 1)$.*

For our purpose, let us list the results for the diffusion-dominated regime. We recall that

$$\gamma^* := \begin{cases} \frac{k-1}{k}, & \text{for } n \geq 1 \text{ and } 0 < s < \frac{1}{2}, \\ +\infty, & \text{for } n \geq 2 \text{ and } \frac{1}{2} \leq s < \frac{n}{2}. \end{cases}$$

**Theorem A.1.** *Let $n \geq 1$ and $k \in (-n, 0)$. All solutions of (A.1) are radially symmetric non-increasing. Let $\mathcal{F}$ be defined by (1.5). If $\gamma > \gamma_c$, then there exists a global minimiser $\bar{\rho}$ of $\mathcal{F}$ on $\mathcal{Y}$. Further, all global minimisers $\bar{\rho} \in \mathcal{Y}$ are radially symmetric non-increasing, compactly supported, uniformly bounded and $C^\infty$ inside their support. Moreover, all global minimisers of $\mathcal{F}$ are solutions of (A.1), according to **Definition A.1**, whenever $\gamma_c < \gamma < \gamma^*$. Finally, if $\gamma_c < \gamma \leq 2$, we have $\bar{\rho} \in W^{1,\infty}(\mathbb{R}^n)$.*

**Theorem A.2.** *Let $n = 1$, $k \in (-1, 0)$, and $\gamma > \gamma_c$. All solutions of (A.1) are global minimisers of the energy $\mathcal{F}$ on $\mathcal{Y}$. Further, solutions of (A.1) in $\mathcal{Y}$ are unique.*

## Appendix B. Some important inequalities

In this subsection, we list some inequalities extensively used in our proof.

The first one is a Hardy-type inequality.

**Lemma B.1** (Hardy-type inequality I)**.** *Let $k > 1$ be a given real number and $f$ be a function satisfying*

$$\int_0^{\frac{1}{2}} x^k(|f|^2 + |\partial_x f|^2)\,\mathrm{d}x < \infty,$$

*then it holds that*

$$\int_0^{\frac{1}{2}} x^{k-2}|f|^2 \, \mathrm{d}x \leq C \int_0^{\frac{1}{2}} x^k (|f|^2 + |\partial_x f|^2) \, \mathrm{d}x.$$

As a consequence of **Lemma B.1**, we also note the following result.

**Lemma B.2** (Hardy-type inequality II). *Let $k > 1$ be a given real number and $f$ be a function satisfying*

$$\int_{\frac{1}{2}}^1 (1-x)^k (|f|^2 + |\partial_x f|^2) \, \mathrm{d}x < \infty,$$

*then it holds that*

$$\int_{\frac{1}{2}}^1 (1-x)^{k-2}|f|^2 \, \mathrm{d}x \leq C \int_{\frac{1}{2}}^1 (1-x)^k (|f|^2 + |\partial_x f|^2) \, \mathrm{d}x.$$

The proof of **Lemma B.1** and **Lemma B.2** can be found in [21].

By **Theorem A.1**, for $m_c < m < 2$, $\bar{\rho} \in W^{1,\infty}(\mathbb{R})$, therefore for $\varepsilon > 0$ sufficiently small, we have

$$c_1(x + \bar{R}) \leq \bar{\rho}(x) \leq c_2(x + \bar{R}), \quad -\bar{R} \leq x \leq -\bar{R} + \varepsilon,$$
$$c_1(\bar{R} - x) \leq \bar{\rho}(x) \leq c_2(\bar{R} - x), \quad \bar{R} - \varepsilon \leq x \leq \bar{R},$$

for some constant $c_1, c_2 > 0$. Therefore

$$\bar{\rho}(x) \sim 1 - |x| \text{ as } x \text{ is close to } \pm 1. \tag{B.1}$$

Then the following result follows directly by applying **Lemma B.1** and **Lemma B.2**.

**Corollary B.1.** *Let $k > 1$ be a given real number, $\bar{\rho}$ be the steady solution in* **Theorem A.1** *for $N = 1$ and $m_c < m < 2$ compactly supported in $[-\bar{R}, \bar{R}]$ for some constant $\bar{R} > 0$. If $f$ satisfies*

$$\int_{-\bar{R}}^{\bar{R}} \bar{\rho}^k (|f|^2 + |\partial_x f|^2) \, \mathrm{d}x < \infty,$$

*then it holds that*

$$\int_{-\bar{R}}^{\bar{R}} \bar{\rho}^{k-2}|f|^2 \, \mathrm{d}x \leq C \int_{-\bar{R}}^{\bar{R}} \bar{\rho}^k (|f|^2 + |\partial_x f|^2) \, \mathrm{d}x.$$

In addition, we need the following results in our proof.

**Lemma B.3.** *For $0 < \alpha < 1$, there exists a constant $C > 0$ such that*

$$\int_{-1}^1 \frac{1}{\bar{\rho}^\alpha(y)|x-y|^{1-2s}} dy \leq C \bar{\rho}^{2s-\alpha}(x). \tag{B.2}$$

*Furthermore, if $1 + 2s > \alpha + \beta$ and $\beta > 0$, then*

$$\int_{-1}^1 \int_{-1}^1 \frac{1}{\bar{\rho}^\alpha(x)\bar{\rho}^\beta(y)|x-y|^{1-2s}} \, \mathrm{d}x \, \mathrm{d}y \leq C. \tag{B.3}$$

*Proof.* To prove (B.2), by (B.1), we have

$$\int_{-1}^{1} \frac{1}{\bar{\rho}^{\alpha}(y)|x-y|^{1-2s}} \, \mathrm{d}y$$
$$\leq C \left( \int_{-1}^{x} \frac{1}{(1+y)^{\alpha}(x-y)^{1-2s}} \, \mathrm{d}y + \int_{x}^{1} \frac{1}{(1-y)^{\alpha}(y-x)^{1-2s}} \, \mathrm{d}y \right).$$

Let us compute

$$I_1(x) := \int_{-1}^{x} \frac{1}{(1+y)^{\alpha}(x-y)^{1-2s}} \, \mathrm{d}y,$$

$$I_2(x) := \int_{x}^{1} \frac{1}{(1-y)^{\alpha}(y-x)^{1-2s}} \, \mathrm{d}y.$$

For $I_1$, by using the change of variable $t = \frac{1+y}{1+x}$, then

$$I_1 = (1+x)^{2s-\alpha} \int_{0}^{1} t^{-\alpha}(1-t)^{-1+2s} \, \mathrm{d}t.$$

We recall the definition of Beta function,

$$B(a,b) := \int_{0}^{1} t^{a-1}(1-t)^{b-1} \, \mathrm{d}t,$$

then

$$I_1(x) = (1+x)^{2s-\alpha} B(1-\alpha, 2s). \tag{B.4}$$

Similarly, we have

$$I_2(x) = (1-x)^{2s-\alpha} B(1-\alpha, 2s). \tag{B.5}$$

Therefore we conclude (B.2) by (B.4), (B.5) and (B.1).

Now, we use (B.2) to prove (B.3) and we find that

$$\int_{-1}^{1} \int_{-1}^{1} \frac{1}{\bar{\rho}^{\alpha}(x)\bar{\rho}^{\beta}(y)|x-y|^{1-2s}} \, \mathrm{d}x \, \mathrm{d}y \leq \int_{-1}^{1} \bar{\rho}^{2s-\alpha-\beta}(y) \, \mathrm{d}y.$$

Therefore we conclude (B.3) by (B.1).

$\square$

As we mentioned in Remark 2.3, we give the following estimates on $\|\bar{\rho}^{-\frac{1}{2}}\Phi(t)\|_{L^2}$.

**Lemma B.4.** *For $\gamma < 1 + 2s/3$, there exists a constant $C > 0$ such that*

$$\|\bar{\rho}^{-\frac{1}{2}}\Phi(t)\|_{L^2} \leq C(\|\partial_x \eta(t)\|_{L^{\infty}}^2 + \|\bar{\rho}^{\frac{2\gamma-1}{2}}\partial_{xx}\eta(t)\|_{L^2}^2). \tag{B.6}$$

*Proof.* It holds that

$$\|\bar{\rho}^{-\frac{1}{2}}\Phi(t)\|_{L^2}^2 \leq \int_{-1}^{1} \frac{1}{\bar{\rho}(x)} \left( \int_{-1}^{1} \frac{\bar{\rho}(y)|\partial_x \eta(t,x) - \partial_y \eta(t,y)|}{|x-y|^{2-2s}} \, \mathrm{d}y \right)^2 \, \mathrm{d}x$$
$$+ \int_{-1}^{1} \frac{1}{\bar{\rho}(x)} \left( \int_{-1}^{1} \frac{|\bar{\rho}(x) - \bar{\rho}(y)| \cdot |\partial_y \eta(t,y)|}{|x-y|^{2-2s}} \, \mathrm{d}y \right)^2 \, \mathrm{d}x$$
$$+ \int_{-1}^{1} \frac{1}{\bar{\rho}(x)} \left( \int_{-1}^{1} \bar{\rho}(y)|\partial_x \eta(t,x) - \partial_y \eta(t,y)| \cdot |R_{1,0}(t,x,y)| \, \mathrm{d}y \right)^2 \, \mathrm{d}x$$
$$+ \int_{-1}^{1} \frac{1}{\bar{\rho}(x)} \left( \int_{-1}^{1} |\bar{\rho}(x) - \bar{\rho}(y)| \cdot |\partial_y \eta(t,y) R_{1,0}(t,x,y)| \, \mathrm{d}y \right)^2 \, \mathrm{d}x.$$

By the estimate (3.51) for $R_{1,0}$ and choosing $\delta \in (1-2s, 4-3\gamma)$, we have

$$\|\bar{\rho}^{-\frac{1}{2}}\Phi(t)\|_{L^2}^2 \leq \int_{-1}^1 \frac{1}{\bar{\rho}(x)}\left(\int_{-1}^1 \frac{\bar{\rho}(y)|\partial_x\eta(t,x) - \partial_y\eta(t,y)|}{|x-y|^{2-2s}}\,\mathrm{d}y\right)^2 \mathrm{d}x$$
$$+ \int_{-1}^1 \frac{1}{\bar{\rho}(x)}\left(\int_{-1}^1 \frac{|\bar{\rho}(x) - \bar{\rho}(y)| \cdot |\partial_y\eta(t,y)|}{|x-y|^{2-2s}}\,\mathrm{d}y\right)^2 \mathrm{d}x$$
$$\leq \underbrace{\int_{-1}^1 \frac{1}{\bar{\rho}(x)}\left(\int_{-1}^1 \frac{\bar{\rho}^{\frac{-1+\delta}{2}}(y)|\bar{\rho}^{\frac{3-\delta}{2}}(x)\partial_x\eta(t,x) - \bar{\rho}^{\frac{3-\delta}{2}}(y)\partial_y\eta(t,y)|}{|x-y|^{2-2s}}\,\mathrm{d}y\right)^2 \mathrm{d}x}_{I_1}$$
$$+ \underbrace{\int_{-1}^1 \frac{1}{\bar{\rho}(x)}\left(\int_{-1}^1 \frac{\bar{\rho}^{\frac{-1+\delta}{2}}(y)|\bar{\rho}^{\frac{3-\delta}{2}}(x) - \bar{\rho}^{\frac{3-\delta}{2}}(y)| \cdot |\partial_x\eta(t,x)|}{|x-y|^{2-2s}}\,\mathrm{d}y\right)^2 \mathrm{d}x}_{I_2}$$
$$+ \underbrace{\int_{-1}^1 \frac{1}{\bar{\rho}(x)}\left(\int_{-1}^1 \frac{|\bar{\rho}(x) - \bar{\rho}(y)| \cdot |\partial_y\eta(t,y)|}{|x-y|^{2-2s}}\,\mathrm{d}y\right)^2 \mathrm{d}x}_{I_3}.$$

For $I_1$, by Hölder's inequality, Lemma (B.3) and the choice of $\delta$ such that $\delta > 1-2s$, we have

$$I_1 \leq \int_{-1}^1 \int_{-1}^1 \frac{1}{\bar{\rho}(x)\bar{\rho}^{1-\delta}(y)|x-y|^{1-2s}}\,\mathrm{d}y\,\mathrm{d}x$$
$$\times \int_{-1}^1 \int_{-1}^1 \frac{|\bar{\rho}^{\frac{3-\delta}{2}}(x)\partial_x\eta(t,x) - \bar{\rho}^{\frac{3-\delta}{2}}(y)\partial_y\eta(t,y)|^2}{|x-y|^{3-2s}}\,\mathrm{d}y\,\mathrm{d}x$$
$$\leq C\|\bar{\rho}^{\frac{3-\delta}{2}}\partial_x\eta(t)\|_{H^{1-s}}^2.$$

By Gagliardo-Nirenberg's inequality,

$$\|\bar{\rho}^{\frac{3-\delta}{2}}\partial_x\eta(t)\|_{H^{1-s}} \leq C\|\bar{\rho}^{\frac{3-\delta}{2}}\partial_x\eta(t)\|_{L^2}^s(\|\bar{\rho}^{\frac{1-\delta}{2}}\partial_x\eta(t)\|_{L^2} + \|\bar{\rho}^{\frac{3-\delta}{2}}\partial_{xx}\eta(t)\|_{L^2})^{1-s}.$$

Then by the choice of $\delta$ such that $\delta < 4-3\gamma$, we have

$$I_1 \leq C(\|\partial_x\eta(t)\|_{L^\infty}^2 + \|\bar{\rho}^{\frac{2\gamma-1}{2}}\partial_{xx}\eta(t)\|_{L^2}^2). \tag{B.7}$$

For $I_2$ and $I_3$, by Hölder's inequality, we have

$$I_2 \leq \|\partial_x\eta(t)\|_{L^\infty}^2 \int_{-1}^1 \int_{-1}^1 \frac{1}{\bar{\rho}(x)\bar{\rho}^{1-\delta}(y)|x-y|^{1-2s}}\,\mathrm{d}y\,\mathrm{d}x$$
$$\times \int_{-1}^1 \int_{-1}^1 \frac{|\bar{\rho}^{\frac{3-\delta}{2}}(x) - \bar{\rho}^{\frac{3-\delta}{2}}(y)|^2}{|x-y|^{3-2s}}\,\mathrm{d}y\,\mathrm{d}x,$$

and

$$I_3 \leq \|\partial_x\eta(t)\|_{L^\infty}^2 \int_{-1}^1 \int_{-1}^1 \frac{1}{\bar{\rho}(x)|x-y|^{1-2s}}\,\mathrm{d}y\,\mathrm{d}x \int_{-1}^1 \int_{-1}^1 \frac{|\bar{\rho}(x) - \bar{\rho}(y)|^2}{|x-y|^{3-2s}}\,\mathrm{d}y\,\mathrm{d}x.$$

Then following the estimate for $I_1$, by Lemma B.3 and Gagliardo-Nirenberg's inequality,

$$I_2, I_3 \leq C\|\partial_x\eta(t)\|_{L^\infty}^2. \tag{B.8}$$

Therefore, (B.6) is proved by (B.7) and (B.8).                    $\square$

**Acknowledgment:** JAC was supported by the Advanced Grant Nonlocal-CPD (Nonlocal PDEs for Complex Particle dynamics: Phase Transitions, Patterns and Synchronization) of the European Research Council Executive Agency (ERC) under the European Union's Horizon 2020 research and innovation programme (grant agreement No. 883363). JAC was also partially supported by the "Maria de Maeztu" Excellence Unit IMAG, reference CEX2020-001105-M, funded by the Spanish ministry of Science MCIN/AEI/10.13039/501100011033/ and the EPSRC grant numbers EP/T022132/1 and EP/V051121/1. RJD was partially supported by the General Research Fund (Project No. 14303523) from RGC of Hong Kong and also by the grant from the National Natural Science Foundation of China (Project No. 12425109). AWK was partially supported by National Science Centre Poland grant Sonata Bis no. 2020/38/E/ST1/00469. RJD and JHZ thank Professor Tao Luo at City University of Hong Kong for his helpful discussions on this work. RJD and JAC thank the Kan Ton Po fellowship scheme of the Royal Society that allowed us to initiate this research.

**Data availability:** The manuscript contains no associated data.

**Conflict of Interest:** The authors declare that they have no conflict of interest.

(J. A. Carrillo) Mathematical Institute, University of Oxford, Woodstock Road, Oxford, OX2 6GG, United Kingdom.

*Email address*: `jose.carrillo@maths.ox.ac.uk`

(R.-J. Duan) Department of Mathematics, The Chinese University of Hong Kong, Shatin, Hong Kong, P.R. China

*Email address*: `rjduan@math.cuhk.edu.hk`

(A. Wróblewska-Kamińska) Institute of Mathematics, Polish Academy of Sciences, Śniadeckich 8, 00-656 Warszawa, Poland

*Email address*: `awrob@impan.pl`

(J.-H. Zhang) Department of Mathematics, The Chinese University of Hong Kong, Shatin, Hong Kong, P.R. China

*Email address*: `jhzhang@math.cuhk.edu.hk`